%% file: BroussousSchneiderFinal.tex
\input defps

\parindent=0pt

\def\s{\sigma}
\def\lra{\longrightarrow}
\def\VV{{\cal V}}
\def\CC{{\cal C}}
\def\XE{{\rm sd}(X_E )}

\def\XX{{\rm sd}(X)}
\def\RR{{\cal R}}
\def\FF{{\cal F}}
\def\HH{{\cal H}}
\def\KK{{\cal K}}
\def\SS{{\cal S}}
\def\boundary{\buildrel \partial\over\longrightarrow}
\def\augmentation{\buildrel \epsilon\over\longrightarrow}
\def\psiv{\buildrel \vee\over\psi}
\def\fv{\buildrel \vee\over{f}}
\def\Rfr{\cal R}
\def\St{\rm St}
\def\WW{{\cal W}}
\def\UU{{\cal U}}
\def\Jm{J_{\rm max}}
\def\kmax{\kappa_{\rm max}}
\def\Bfrm{\Bfr_{\rm max}}
\def\lm{\lambda_{\rm max}}

\def\em{e_{\rm max}}
\def\Ap{{\cal A}}
\def\AL{{\cal A}_{L}}
\def\uVV{\underline{\cal V}}

\def\augmentationL{\buildrel \epsilon_{L}\over\longrightarrow}
~

\medskip

\centerline{\bf Simple characters and coefficient systems on the
building}
\bigskip

\centerline{P. Broussous and P. Schneider}
\medskip

\centerline{Version : february 2014}

\vskip1cm

{\it Abstract.}  Let $F$ be a non-archimedean local field and $G$ be the group ${\rm GL}(N,F)$ for some integer
$N\geq 2$. Let $\pi$ be a smooth complex representation of $G$ lying in the Bernstein block ${\cal B}(\pi )$ of  some
simple type in the sense of Bushnell and Kutzko. Refining the approach of the second author and U. Stuhler,
we canonically attach to $\pi$ a subset $X_\pi$ of the Bruhat-Tits building $X$ of $G$, as well as  a $G$-equivariant
coefficient system $\CC [\pi ]$ on $X_\pi$.
 Roughly speaking the coefficient system is obtained by taking isotypic components
of $\pi$ according to some representations constructed from the Bushnell and Kutzko type of $\pi$. We conjecture
that when $\pi$ has central character, the augmented chain complex associate to $\CC (\pi )$ is a projective 
resolution of $\pi$
in the category ${\cal B}(\pi )$. Moreover we reduce this conjecture to a technical lemma of representation
theoretic nature. We prove this lemma when $\pi$ is an irreducible discrete series of $G$. We then attach to any
 irreducible discrete series $\pi$ of $G$ an explicit pseudo-coefficient $f_\pi$ and obtain
a Lefschetz type  formula for the value of the Harish-Chandra character of $\pi$ at a regular elliptic element.
In contrast to that obtained by U. Stuhler and the second author,
 this formula allows explicit character value computations.

{\it R\'esum\'e}. Soient $F$ un corps local non archim\'edien et $G$ le groupe ${\rm GL}(N,F)$, pour un entier
$N\geq 2$. Soit $\pi$ une repr\'esentation lisse complexe de $G$ appartenant au block de Bernstein ${\cal B}(\pi )$
d'un type simple au sens de Bushnell et Kutzko. En affinant l'approche que proposent le second auteur et U.
Stuhler, nous attachons canoniquement \`a $\pi$ un sous-ensemble $X_\pi$ de l'immeuble de Bruhat-Tits $X$
de $G$, ainsi qu'un syst\`eme de coefficients $G$-\'equivariant $\CC [\pi ]$ sur $X_\pi$. Grossi\`erement parlant, 
le syst\`eme de coefficients est construit en prenant des composantes isotypiques de $\pi$ selon des repr\'esentations
construites \`a partir du type de Bushnell et Kutzko de $\pi$. Nous conjecturons que lorsque $\pi$ poss\`ede un 
caract\`ere central, le complexe de cha\^{\i}nes augment\'e associ\'e \`a $\CC (\pi )$ est une r\'esolution de $\pi$
dans la cat\'egorie ${\cal B}(\pi )$. De plus nous r\'eduisons cette conjecture \`a un lemme technique
en th\'eorie des repr\'esentations. Nous d\'emontrons ce lemme  lorsque $\pi$ est une repr\'esentation
irr\'eductible de la s\'erie discr\`ete de $G$. Nous attachons ensuite \`a toute
repr\'esentation irr\'eductible $\pi$ de la s\'erie discr\`ete de $G$ un pseudo-coefficient explicite $f_\pi$ et
obtenons une formule de type Lefschetz pour la valeur du caract\`ere de Harish-Chandra de $\pi$ en un \'el\'ement
elliptique r\'egulier. Contrairement à celle obtenue par U. Stuhler et le second auteur,
 notre formule permet des calculs explicites. 
 
\vskip2cm

\centerline{\bf Introduction}
\vskip1cm

 Let $F$ be a non-archimedean local field and, for some integer $N\geq 2$, let 
$G$ denote the locally compact group ${\rm GL}(N,F)$ and $X$ its Bruhat-Tits building
The aim of this work is to refine
the construction of [SS] (also see [SS2]) to attach to certain representations of $G$
new equivariant coefficient systems on the Bruhat-Tits building. These representations
belong to the Bernstein blocks  of the category of smooth complex representations of
$G$ corresponding to {\it simple types} in the sense of Bushnell and Kutzko [BK1].
 Let $(\pi ,\VV )$ be a smooth complex representation
of $G$. In [SS] an equivariant  coefficient system  $\CC (\pi )$ is  constructed by
attaching to each simplex $\sigma$ of
$X$ the  space of vectors fixed by  a certain congruence subgroup of level $e$ of the
parahoric subgroup of $G$ fixing $\sigma$. Here the integer $e$ is such that $\VV$ is
generated as a $G$-module by its vectors fixed by the principal congruence subgroup
of level $e$ of some maximal compact subgroup of $G$.  In [SS] it is proved that the augmented
chain complex $C_{\bullet} (X, \CC (\pi ))\lra \VV$ of $X$ with coefficients in $\CC (\pi )$
is exact. If one moreover  assumes that $(\pi ,\VV )$ admits a
central character $\chi$, then  $C_{\bullet} (X, \CC (\pi ))\lra \VV$ is a projective resolution
of $(\pi ,\VV)$ in the category of smooth representations of $G$ with central character $\chi$.
In [Br], the first author gave another proof of this fact for Iwahori-spherical representations.
In [SS2], the second author and  U. Stuhler draw some important consequences concerning
the harmonic analysis on $G$ as well as the homological algebra of the category of smooth
representations of $G$.  Among other things they prove that these projective
resolutions give rise to pseudo-coefficients for discrete series representations
(generalizing the pseudo-coefficient constructed by Kottwitz
in [Kot] for the Steinberg representation) as well as a Lefschetz type character formula for
the Harish-Chandra character of any  smoooth representation.  Note that if the construction of
[SS] is restricted to the group $G$, [SS2] gives a generalization to any connected reductive
$F$-group $\Gdss$  and most of its results are  valid  without restriction on $G$
(but  sometimes $F$ is assumed to have characteritic $0$, and $\Gdss (F)$ to
 have compact center).

 If the construction and results of [SS], [SS2] have important theoretic consequences,
they do not allow explicit calculations. Indeed in general the coefficient system $C(\pi )$
cannot be explicitely computed (except may be  in the {\it level $0$ case}, but this is
nowhere written). Indeed the only explicit way to be given an irreducible
smooth representation of $G$ is to specify its Bushnell and Kutzko type. This is why
it is natural to seek for a refinement of [SS] based on Bushnell and Kutzko theory.

 In this paper, for technical reasons, we restrict to  representations belonging to
Bernstein blocks of $G$ attached to simple types.  These Bernstein blocks are exactly 
those containing discrete series representations. We fix a simple type $(J,\lambda )$
and denote by $\RR_{\lambda}(G)$ the category of smooth representations of $J$ that
are generated by their $\lambda$-isotypic component. We fix a smooth representation
$(\pi ,\VV )$ of $G$ lying in $\RR_\lambda (G)$. To the datum $(J
,\lambda )$, in a non canonical way, one may associate a field
extension $E/F$ of degree dividing $N$ whose multiplicative group
$E^\times$ is embedded in $G$. The centralizer $G_E$ of $E^\times$ in
$G$ is isomorphic to ${\rm GL}(N/[E:F], E)$. Using a result of the first author and
B. Lemaire [BL], we may view the Bruhat-Tits building $X_E$ of $G_E$
as being embedded in $X$ in a $G_E$-equivariant way. We show that {\it hidden}
in the properties of {\it Heisenberg representations} constructed in [BK1]{\S}(5.1)
and in the {\it mobility} of simple characters established in {\it loc. cit.} {\S}(3.6),
there is a {\it geometric structure} allowing to attach to $\pi$ a $G_E$-equivariant
coefficient system $\CC_E [\pi ]$ on the first barycentric subdivision ${\rm sd}(X_E )$ of $X_E$.
More precisely, in a non canonical way, we attach to $(J,\lambda )$ a collection of
pairs $(J^1 (\sigma ,\tau ), \eta (\sigma ,\tau ))_{\sigma\subset \tau}$, where
$\sigma$ and $\tau$ run over the simplices of $X_E$ satisfying $\sigma\subset \tau$. Here
$J^1 (\sigma ,\tau )$ is some compact open subgroup of $G$ and $\eta (\sigma ,\tau )$
a Heisenberg representation of $J^1 (\sigma ,\tau )$ as considered in {\it loc. cit.}
 (5.1.14) (but Bushnell and Kutzko do not use this language nor this notation). Moreover
 the collection $(J^1 (\sigma ,\tau ), \eta (\sigma ,\tau ))_{\sigma\subset \tau}$ is 
$G_E$-equivariant. Exploiting the compatibility relations among the various 
$\eta (\sigma ,\tau )$ proved in {\it loc. cit.} {\S}(5.1),
and by taking isotypic components of $\VV$ according to the Heisenberg representations
$\eta (\sigma ,\tau )$, we construct our equivariant coefficient system $\CC_ E [\pi ]$.

 We then show that the subset $X[E]$ of $X$ obtained by taking the
 union of the $g.X_E$, where $g$ runs over $G$ has the structure of
 a $G_E$-simplicial complex containing $X_E$ as a subcomplex. We
 naturally attach to $\CC_E [ \pi ]$ a $G$-equivariant coefficient
 system $\CC [\pi ]$ on the first barycentric subdivision of $X[E]$
 and show that it actually derives from a coefficient complex on
 $X[E]$, still denoted by $\CC [\pi ]$.  We prove that
 the simplicial complex $X[E]$ and the coefficient system $\CC [\pi ]$
 are actually independent of any choice made in their construction:
 these are objects canonically attached to $\pi$. Moreover the support $X_\pi$
 of $\CC [\pi ]$ maybe explicitely determined. In [BK1]{\S}5, the Hecke
 algebra of $(J,\lambda )$ is described using a non canonical
 unramified  field extension $L/E$. It gives rise to a general linear group
 $G_L\subset G_E\subset G$, to a Bruhat-Tits building $X_L\subset
 X_E \subset X$ and to a simplicial complex
$$
X[L]=\bigcup_{g\in G}g.X_L \subset X[E]\ .
$$
Then the support of $\CC [\pi ]$ is $X_\pi = X[L]$.

We then consider the augmented chain complex of $X_\pi$ with
coefficients in $\CC [\pi ]$:
$$
C_{\bullet} (X_\pi , \CC [\pi ])\lra \VV\ .
\leqno{\hbox{($*$)}}
$$
We show that this  complex lies  in the category $\RR_\lambda (G)$. We
cannot in general prove its exactness  that we consider as a
conjecture. However we propose a strategy to tackle this exactness that generalizes the
approach that the first author uses in  [Br]. Indeed if $(\pi ,\VV )$ has level $0$ then 
$X[L]=X$ and the coefficient system $\CC [\pi ]$ coincides with that constructed in [SS].
 In [Br], for Iwahori-spherical representations (they have level $0$),  one proves 
the exactness of ($*$) using type theory and an argument of geometric nature.

 Let us explain how this generalized approach works. Let $\HH (G)$ be the Hecke algebra
 of locally constant complex functions with compact support on $G$. It is equipped with 
the convolution product $\star$ coming from a fixed Haar measure on $G$. Let $e_\lambda$ 
be the idempotent of $\HH (G)$ attached to $\lambda$ so that
for any smooth complex representation $\cal W$ of $G$, 
$e_\lambda \star {\cal W}={\cal W}^{\lambda}$ is the $\lambda$-isotypic
component of $\cal W$. One basic fact of type theory is that  the functor
$$
\RR_\lambda (G)\lra e_\lambda \star \HH (G) \star e_\lambda -{\rm Mod} \ ,
 \ {\cal W}\lra {\cal W}^{\lambda}
$$
induces an equivalence of categories. It follows that in order to prove the exactness of  ($*$),
 we are reduced to proving the exactness of the chain complex ($**$) in 
$e_\lambda .\HH (G).e_\lambda$-${\rm Mod}$ obtained 
 from ($*$) by applying the functor ${\cal W}\lra {\cal W}^{\lambda}$:
$$
 C_{\bullet} (X_\pi , \CC [\pi ])^{\lambda} \lra \VV^{\lambda}\ .
\leqno{\hbox{($**$)}}
$$

 In fact we shall not work with the type $\lambda$, but with an equivalent type $\lambda '$
 defining the same Bernstein block; to make things simpler we ignore this difficulty
in the introduction. Then generalizing [Br]  we prove that modulo a conjectural technical
hypothesis (Conjecture (X.4.1)), as a complex of $\Cdss$-vector spaces, ($**$) is canonically
isomorphic to the augmented chain complex of a certain apartment ${\cal A}_L$ of $X_L$ with
{\it constant} coefficients in $\VV^{\lambda}$. Of course ${\cal A}_L$ being a finite dimensional
euclidean space, it is a contractible topological space, and its augmented chain complex with
 constant coefficients in any abelian group is exact.

 We prove Conjecture (X.4.1), whence the exactness of ($*$), when the representation $\pi$
belongs to the discrete series of $G$. Indeed in that case we are able to entirely compute
the coefficient system $\CC [\pi ]$ by using some technical lemmas proved by the second author
and Zink in [SZ]. We actually  prove that there exists a $G$-equivariant
 collection  of pairs $(G_\sigma ,\lambda_\sigma )$ such that the coefficient system is given by
$\CC [\pi ]_\sigma =\VV^{\lambda_\sigma}$ (isotypic component), where $\sigma$ runs over the simplices
of $X_\pi$, $G_\sigma$ denotes the stabilizer of $\sigma$ in $G$, and $\lambda_\sigma$ is an
irreducible smooth representation of $G_\sigma$.  Moreover for any  simplex $\sigma$ of $X[L]$,
the restriction of $\lambda_\sigma$ to the maximal compact subgroup of $G_\sigma$ only depends
on $(J,\lambda )$ but not on $\pi$.

 Closely following  [SS2],  we attach to the coefficient system $\CC [\pi ]$ an \break
Euler-Poincar\'e function  $f_{\rm EP}^\pi$ on $G$ and prove that it is a pseudo-coefficient of 
$\pi$. This pseudo-ceofficient should be very close to  that constructed in [Br2] by the 
first author using an entirely different approach (but also based on Bushnell and Kutzko type 
theory), however the comparison has to be done. In contrast with that of [Br2], the pseudo 
coefficient $f_{\rm EP}^\pi$ is given  by a formula adapted to explicit computations. 
In particular by computing certain orbital integrals,   we derive a Lefschetz type character
 formula  for the value of the Harish-Chandra character $\Theta_\pi$ of $\pi$ at  a regular 
elliptic element $\gamma$ of $G$. This formula takes the form:
$$
\Theta_\pi (\gamma )={\rm Tr}\, \left(\gamma\  ,\  {\rm EP}\, 
H^* (X_\pi^\gamma ,\CC [\pi ])\ \right)
\leqno{\hbox{($***$)}}
$$
where  ${\rm EP}\, H^* (X_\pi^\gamma ,\CC  [\pi ])$ denotes the homology Euler-Poincar\'e module of
 the restriction of $\CC (\pi )$ to the fixed point set $X_\pi^\gamma$ of $\gamma$ in $X_\pi$.
 We cannot expect to make formula ($***$) entirely explicit. Indeed if $\gamma$ is an element of $G$
there is no known easy description of the  fixed  point set $X^\gamma$. Nevertheless when the elliptic
regular element $\gamma$ is {\it minimal over $F$} in the sense of Bushnell and Kutzko,
then $X_\pi^\gamma$ is either empty or reduced to a point. In that case the  Lefschetz formula for
$\Theta_\pi (\gamma )$ takes a striking simple form and allows explicit computations. In particular,
in that case we recover the two character formulas obtained in [Br2]. However our approach
gives a much more general result.

The paper is organized as follows. In section I we  establish some crucial
properties of the embedding $X_E \lra X$, where $E/F$ is a field extension such that $E^\times$
embeds in $G$. In sections II and III we review the main properties of simple characters
and of their endo-classes.  The construction of the $G_E$-equivariant coefficient system $\CC_E [\pi ]$ 
on $X_E$ is given in section IV and its extension  $\CC [\pi ]$ to a $G$-equivariant coefficient system
is done in sections V and VI. The canonicity of the coefficient
complex $\CC [\pi ]$ is studied in sections VII and VIII. To state this
result the right language is that of {\it endo-classes} (Propositions
(VII.2) and (VIII.1.2)). The support of $\CC [\pi ]$ is described in
Proposition (VIII.2.6). In section IX we prove that the chain complex
attached to $\CC [\pi ]$ actually lies in the Bernstein block of $\pi$
(Proposition (IX.2)). In section X we reduce the acyclicity of the
augmented chain complex attached to $\CC [\pi ]$ to a technical lemma
(Conjecture (X.4.1)). For an irreducible discrete series
representation, the conjecture is proved in section XI (Theorem
XI.2.7). The last section XII is devoted to applications. We first
construct an explicit pseudo-coefficient for any irreducible discrete
series representations (Theorem (XII.2.3)) and then derive an explicit
character formula for the Harish-Chandra character of such a
representation (Theorem (XII.3.2)). For elliptic minimal element the
formula simplifies a lot (Proposition (XII.4.4)) and give a new proof
of formulas already obtained in [Br2].

We shall assume that the reader is familiar with the formalism of
[BK1]. Indeed this work may be somehow viewed as a geometric reformulation of
Bushnell and Kutzko's construction of the discrete series of $G$.

We want to thank Shaun Stevens for his help. Proposition (XI.1.2) and
its proof are due to him as well as the proof of  Lemma (X.4.4).

 This work has a long story. Both authors started to collaborate as 
  the first one was in post-doctoral stay in Muenster in 2000/2001.
Results from sections I to  IX  where obtained already in 2004. 
\vfill\eject

\centerline{\bf Table of contents}
\vskip1cm

I. Field extensions and centralizers.\hfill\break
\hbox{}\hskip0.5cm I.1. Vector spaces and orders.\hfill\break
\hbox{}\hskip0.5cm I.2. Buildings.\hfill\break
\hbox{}\hskip0.5cm  I.3. Some properties of the embedding ${\rm sd}(X_E )\rightarrow {\rm
  sd}(X_F )$.\hfill\break
II. Simple characters and their endo-classes.\hfill\break
\hbox{}\hskip0.5cm  II.1. Simple pairs and their realizations.\hfill\break
\hbox{}\hskip0.5cm  II.2. Potential simple characters.\hfill\break
\hbox{}\hskip0.5cm  II.3. Endo-classes of ps-characters.\hfill\break
III. Ps-characters and pairs of orders.\hfill\break
\hbox{}\hskip0.5cm  III.1. Extensions to mixed groups.\hfill\break
\hbox{}\hskip0.5cm  III.2. Extensions to $1$-units of orders.\hfill\break
\hbox{}\hskip0.5cm  III.3. The degenerate case.\hfill\break
IV. The coefficient system on ${\rm sd}(X_E )$.\hfill\break
V. The coefficient system on ${\rm sd}(X)$.\hfill\break
VI. The degenerate case.\hfill\break
VII. Dependence on the endo-class. \hfill\break
VIII. On the support of ${\cal C}(\Theta ,V,{\cal V})$.\hfill\break
\hbox{}\hskip0.5cm   VIII.1. Endo-classes.\hfill\break
\hbox{}\hskip0.5cm   VIII.2. The support of ${\cal C}[\Theta_o ,V,{\cal V}]$.\hfill\break
IX. The chain complex attached to ${\cal C}_{(J,\lambda )}({\cal V})$· \hfill\break
X. Acyclicity of the chain complex: a strategy.\hfill\break
\hbox{}\hskip0.5cm   X.1. Some lemmas on $\lambda_{\rm max}$-isotypic components.\hfill\break
\hbox{}\hskip0.5cm   X.2. Orientation of $X[L]$.\hfill\break
\hbox{}\hskip0.5cm   X.3. Another chain complex. \hfill\break
\hbox{}\hskip0.5cm   X.4. $J_{\max}$-orbits of simplices.\hfill\break
\hbox{}\hskip0.5cm   X.5. Comparison of chain complexes. \hfill\break
XI. Acyclicity in the case of a discrete series representation. \hfill\break
\hbox{}\hskip0.5cm   XI.1. Determination of the chain complex. \hfill\break
\hbox{}\hskip0.5cm XI.2. Proof of conjecture  (X.4.1) for irreducible
discrete series representations.\hfill\break
XII. Explicit pseudo-coefficients for discrete series representations.\hfill\break
\hbox{}\hskip0.5cm   XII.1. The coefficient system ${\cal C}(\pi )$. \hfill\break
 \hbox{}\hskip0.5cm   XII.2. Euler-Poincar\'e functions.\hfill\break
\hbox{}\hskip0.5cm   XII.3. An explicit character formula.\hfill\break
 \hbox{}\hskip0.5cm XII.4. The character of discrete series representations at minimal elements.\hfill\break

\vfill\eject

\centerline {\bf  I. Field extensions and centralizers.}

\medskip

{\bf I.1 Vector spaces and orders}.

If $K$ is a non-archimedean local field we shall denote by  $\ofr_K$
its ring of integers and by $\pfr_K$ the maximal ideal of $\ofr_K$.
Once for all we fix such a field $F$.

 Let $E/F$ be a finite field extension and $V$ a finite
 dimensional $E$-vector space. Then $V$ is naturally an $F$-vector
 space. We write $A = {\rm End}_{F}V$, $G  = {\rm
 Aut}_F V$, $B = {\rm End}_{E}V$ and $G_E = {\rm
 Aut}_{E}V$. We have a natural inclusion of $F$-algebras $B\subseteq A$
and the group $G_E$ is naturally a subgroup of $G$. As an $F$-algebra
 $E$ embeds canonically in $A$ and its centralizer is $B$. Similarly,
 the left action of $E$ on $V$ allows us to see $E^{\times}$ as a
 subgroup of $G$; its centralizer is $G_E$.

 Let ${\rm Her}(A)$ (resp. ${\rm Her}(B)$) denote the set of hereditary
 $\ofr_{F}$-orders in $A$ (resp. hereditary $\ofr_E$-orders in
 $B$). These sets are posets (for inclusion) and $G$ and $G_E$
 respectively act on them by conjugation. We have a natural map
 $j_{\rm order}$: ${\rm Her}(B)\longrightarrow {\rm Her}(A)$, defined
 as follows. If $\Bfr$ is in ${\rm Her}(B)$, it is the stabilizer in $B$ of
 an $\ofr_{E}$-lattice chain ${\cal L}$ in $V$; this lattice chain
 may be seen as an $\ofr_F$-lattice chain in $V$ and $j_{\rm
 order}(\Bfr )$ is the attached order in $A$. We shall use the
 notation $j_{\rm order}(\Bfr )  = \Afr (\Bfr )$. The map $j_{\rm
 order}$ is $G_E$-equivariant and, by [BK1] (1.2.1), its
 image  consists of those orders in ${\rm Her}(A)$
 that are stabilized by $E^{\times}$.

\smallskip

{\bf I.2 Buildings.}

 We keep the notation as in (I.1). Let $X$ (resp. $X_E$) denote
 the semisimple affine building of $G$ (resp. $G_E$). The following
 fact will be crucial for our construction.

\smallskip

{\bf (I.2.1) Theorem}. ([BL] Theorem 1.1).  {\it There exists a unique affine
and $G_E$-equivariant map
$$
j_E : X_E \longrightarrow X .
$$
It induces a bijection $X_E \longrightarrow X^{E^{\times}}$.}

\smallskip

 We are going to give a more precise version of this theorem. Recall
 that the building $X$ is triangulated in a canonical way: it is the
 geometric realization of a $G$-simplicial complex that we still
 denote  by $X$. Let $F(X)$ be the set of simplices of $X$. It is a
 poset for inclusion and is equipped with an action of $G$ via poset
 isomorphisms. It is a standard result (compare [BT] Cor. 2.15) that we have an
 anti-isomorphism of posets, compatible with the $G$-actions:
$$
\matrix{
{\rm Her}(A)^{\rm opp} & \longrightarrow & F(X)\cr
\hfill\Afr & \longmapsto & F(\Afr)\hfill}
$$
where $F(\Afr )$ is the unique simplex stabilized by the
 normalizer of $\Afr$ in  $G$. Similarly, we have an
 anti-isomorphism of posets, compatible with the $G_E$-actions:
$$
\matrix{
{\rm Her}(B)^{\rm opp} & \longrightarrow & F(X_E )\cr
\hfill\Bfr & \longmapsto & F(\Bfr )\hfill }
$$
where the notation is obvious. We write $j_{\rm simp}$ for the
morphism of $G_E$-posets $F(X_E )\longrightarrow F(X)$ obtained from
$j_{\rm order}$ through the two previous isomorphisms.

Let ${\rm sd}(X)$ (resp. ${\rm sd}(X_E )$) be the first barycentric
subdivision of $X$ (resp. of $X_E$). This is the flag complex attached
to the poset $F(X)$ (resp. $F(X_E )$). Since $j_{\rm simp}$ is
increasing, it induces a $G_E$-equivariant simplicial map ${\rm
sd}(X_E )\longrightarrow {\rm sd}(X)$.

\smallskip

{\bf (I.2.2) Proposition}. {\it The map $j_{\rm simp}$: ${\rm sd}(X_E
) \longrightarrow {\rm sd}(X)$ induces $j_E$ on the geometric
realizations.}

{\it Proof.}  Let us  denote by $j_{\rm sd}$ the map $X_{E}
\longrightarrow X$ induced by $j_{\rm simp}$ on the geometric realizations
 (constructed with standard affine simplices). By construction $j_{\rm
 sd}$ is affine and $G_E$-equivariant. By unicity in Theorem (I.2.1),
 it must coincide with $j_E$.

\smallskip

In the sequel we shall use the language of hereditary orders instead
of simplices. In particular a $q$-simplex $\sigma$ in ${\rm sd}(X)$ is
a strictly decreasing sequence of orders $\sigma = (\Afr_0 \supset
\Afr_1 \supset\dots\supset \Afr_q )$. The map $j_{E} = j_{\rm simp}$ is
then given by
$$
j_E ( \Bfr_0 \supset \Bfr_1 \supset \dots
 \supset \Bfr_q ) =(\Afr (\Bfr_0 ) \supset \Afr (\Bfr_1 ) \supset \dots
 \supset \Afr (\Bfr_q )) .
$$
We shall also see ${\rm sd}(X_E )$ as being embedded in ${\rm sd}(X)$
: $j_E$ is now an inclusion.

\smallskip

 The map $j_E$ enjoys another property that is not proved in
 [BL]. Recall that $X_E$ and $X$ have invariant metrics which are
 unique up to a $>0$ factor. Since $G$ (resp. $G_E$) acts transitively
 on the apartments of $X$ (resp. of $X_E$) fixing a metric on $X$
 (resp. on $X_E$) amounts to fixing it on one of its apartments.
\smallskip

{\bf (I.2.3) Proposition}. {\it There exist normalizations of metrics
on $X_E$ and $X$ such that the map $j_E$ is an isometry.}
\smallskip

{\it Proof.} By invariance it suffices to prove that the restriction
of $j_E$ to some apartement ${\cal A}_E$ of $X_E$ is an isometry. By
[BL](5.1), $j_E ({\cal A}_E)$ is contained in an apartment $\cal A$ of
$X_E$. Set $n={\rm Dim}_E\, (V)$ and consider $\Rdss^n$ and $\Rdss^{n/[E:F]}$ equipped with
their standard euclidean structures. Then by the proof of Lemma (4.1)
of [BL], one may choose the apartment $\cal A$ and metrics on $X_E$
and $X$ such that :

 -- $\cal A$ identifies to the orthogonal of $(1,1,...,1)$ in $\Rdss^n$

 -- ${\cal A}_E$ identifies to  the orthogonal of $(1,1,..,1)$ in
    $\Rdss^{n/[E:F]}$

 -- the map $j_E$ is given by the restriction of the following linear
    map:
$$
J \ : \ \Rdss^{n/[E:F]}\lra  \Rdss^n\ , \ (x_1
,...,x_{n/[E:F]})\mapsto (x_i /e +\mu_j )_{i=1,...,n/[E:F],\
j=1,...,[E:F]}
$$
where $e$ is the ramification index of $E/F$ and the $\mu_i$ are some
real constants. It is clear that up to a scalar $J$ is an
isometry. Our result follows.

{\bf I.3 Some properties of the embedding ${\rm sd}(X_E )
\longrightarrow {\rm sd}(X)$.}

We keep the notation as in (I.1) and (I.2). We need first some more
notation and facts on orders. If $\Afr$ is a hereditary $\ofr_F$-order
in $A$, then its multiplicative group is a compact open subgroup of
$G$ that we denote by $U(\Afr )$ (this is indeed a parahoric subgroup
of $G$). Let $\Pfr$ be the Jacobson radical of $\Afr$. Then the
quotient $\Afr /\Pfr$ is a semisimple $\Fdss$-algebra, where $\Fdss$
is the residue field of $F$. In particular the multiplicative group
$(\Afr /\Pfr )^{\times}$ is the group of $\Fdss$-points of a product
of general linear groups defined over $\Fdss$. The subgroup
$U^{1}(\Afr ) = 1+\Pfr $ of $1$-units is a normal subgroup of $U(\Afr
)$ and the quotient canonically identifies with $(\Afr /\Pfr
)^{\times}$.

For $\Bfr$ a hereditary order in $B$, the symbol ${\cal N}(\Bfr )$
denotes the normalizer of $\Bfr$ in $G_{E}$, while if $\Afr$ is a
hereditary order in $A$, ${\cal N}(\Afr )$ denotes the normalizer of
$\Afr$ in $G$.

\smallskip

{\bf (I.3.1) Lemma.} {\it For any hereditary order $\Bfr$ in $B$, we
have}
$$
{\cal N}(\Afr (\Bfr )) = {\cal N}(\Bfr ) U(\Afr (\Bfr )) .
$$

{\it Proof.} Let $(L_{k})_{k\in \Zdss}$ be an $\ofr_{E}$-lattice chain
in $ V$  defining $\Bfr$. Let $v_{\Afr(\Bfr )} : A \lra \Zdss$ be the
the valuation map given by
$$
v_{\Afr}(a) = m \ {\rm iff}\  a\in \Pfr^{m} \backslash \Pfr^{m+1}\ ,
 \  m\in \Zdss
$$
where $\Pfr$ is the radical of $\Afr (\Bfr )$. Write $v_{\Bfr}$ for
the similar map $B\lra \Zdss$ defined by the powers of the radical of
$\Bfr$. From [BK1]{\S}1, we have

\smallskip

{\bf (I.3.2)} $(v_{\Afr })_{|B} = v_{\Bfr}$ and ${\cal N}(\Afr (\Bfr
))\cap G_{E} = {\cal N}(\Bfr )$.

\smallskip

Let $t\Zdss$, $t>0$,  be the image of the group homomorphism
$$
v_{\Afr} \ : {\cal N}(\Afr (\Bfr ))\lra {\Zdss} .
$$
Then ${\cal N}(\Afr (\Bfr )) = z^{\Zdss}U(\Afr (\Bfr ))$ for any $z$
in ${\cal N}(\Afr (\Bfr ))$ with  $\Afr$-valuation $t$. A similar
 statement holds for ${\cal N}(\Bfr )$. Now from [BF]
 one knows that $t$ is the smallest positive period of the map
$k\mapsto {\rm dim}_{{\Fdss}}L_{k}/L_{k+1}$.
  So $t$ is also the  smallest positive period of
 $k\mapsto {\rm dim}_{{\Fdss}_{E}}L_{k}/L_{k+1}$,
 where ${\Fdss}_E$ is the residue class field of $E$. Together with
 (I.3.2) this  implies that we can actually choose $z$ in ${\cal
 N}(\Bfr )$  and the result follows.

\smallskip

{\bf (I.3.3) Lemma.} {\it Let $\s = (\Bfr_{0}^{\s} \supset \dots
 \supset \Bfr_{q}^{\s})$  and $\tau = (\Bfr_{0}^{\tau} \supset \dots
 \supset \Bfr_{q}^{\tau})$ be two $q$-simplices in ${\rm sd}(X_E )$.
 Assume that $\s = g\tau$ for some $g\in G$. Then there exists $g_E$
in $G_E$ such that $\s = g_E \tau$. In particular any $g$ as above can
be written $g = g_E g_{\tau}$ with $g_{E}\in G_E$ and
 $g_{\tau}\in {\rm Stab}_{G}(\tau )$. }

{\it  Proof}. First we need to recall the classification of conjugacy
classes of hereditary orders  in $A$ (cf. [BF] or [Rei]).
 Let $\Afr$ be such an order and let $(L_{k})_{k\in \Zdss}$ be a
 lattice chain in $V$ defining $\Afr$. To $\Afr$ we attach the
sequence of integers $d(\Afr )_{k} = {\rm dim}_{\Fdss}L_{k}/L_{k+1}$,
 $k\in \Zdss$. Then two hereditary orders $\Afr_{1}$, $\Afr_{2}$ are
 conjugate if and only if the sequences $d(\Afr_{1})$ and
$d(\Afr_{2})$ coincide up to a translation of the indexing.
 We use the notation $d_{E}$ for the sequences attached to hereditary
orders  in $B$. If $\Bfr$ is such an order, attached to an
$\ofr_{E}$-lattice  chain $(L_{k})_{k\in \Zdss}$ in $V$, we have:
$$
d(\Afr (\Bfr ))_{k} = [{\Fdss}_E:{\Fdss}] d_{E}(\Bfr )_{k} \ ,
\   k\in {\Zdss} .
$$
We deduce:

\smallskip

{\bf (I.3.4)} {\it Let $\Bfr_{1}$ and $\Bfr_{2}$ be hereditary orders
in $B$. Then they are $G_{E}$-conjugate if and only if the orders
$\Afr (\Bfr_{1})$ and $\Afr (\Bfr_{2} )$ are $G$-conjugate. In other
words Lemma (I.3.3) holds when $q=0$.}

\smallskip

Now let us turn to the general case. By using (I.3.4), we may
replace $\tau$ by a conjugate under $G_E$ so that
$\Bfr_{q}^{\sigma} = \Bfr_{q}^{\tau}=:\Bfr_q$. By assumption there
exists a $g\in G$ such that $\Afr (\Bfr_{i}^{\sigma})  = \Afr
(\Bfr_{i}^{\tau})^{g}$ for $i=0,\dots ,q$. In particular  $\Afr
(\Bfr_{q})  = \Afr (\Bfr_{q})^{g}$, and, thanks to (I.3.1), we
may, by replacing $\tau$ by a $G_E$-conjugate, assume that
$$
\sigma = g\tau\ \ , \ \ \Bfr_{q}^{\sigma}= \Bfr_{q}^{\tau}
=\Bfr_q,\ \hbox{and}\ g\in U(\Afr (\Bfr_q )).
$$
But then $g \in U(\Afr (\Bfr_i ))$ for any $i=0,\dots ,q$ which
means that $g$ fixes $\tau$, i.e., that $\sigma = \tau$.

\smallskip

It is not possible to characterize the image of $\XE$ using numerical
invariants attached to simplices. But we are going to give a criterion
for a simplex of $\XX$ to belong to:
$$
X(E) :=  \bigcup_{g\in G} g \XE \ .
$$
Here $\XE$ is of course seen as being embedded in $\XX$.

Let $(L_{k})_{k\in \Zdss}$ be an $\ofr_F$-lattice chain in $V$
 and $\Afr$ be the attached order in $A$. Write $e=e(\Afr )$
 for the period of $\Afr$. The sequence of positive integers
 defined  by $d(\Afr )_{k} = {\rm dim}_{\Fdss}L_{k}/L_{k+1}$ is
 $e$-periodic and we have the partition:
$$
n = {\rm dim}{V} = d(\Afr )_{0} +\dots +d(\Afr )_{e-1}\ .
$$
We denote by  $p(\Afr )$  the least positive period of $(d(\Afr
)_{k})_{k\in \Zdss}$. We can rephrase [BK2] Prop. (1.2) as follows.

\smallskip

{\bf (I.3.5) Proposition}. {\it The order $\Afr$ has a conjugate
normalized  by $E^{\times}$ if and only if the following assertions
hold:

  i) $f(E/F)$ divides $d(\Afr )_{k}$ for all $k\in \Zdss$;

  ii) $e(E/F)$ divides $e(\Afr )/p(\Afr )$.

In other words the vertices of ${\rm sd}(X)$ which are in $X(E)$ are
exactly those vertices which correspond to hereditary orders $\Afr$
satisfying conditions (i) and (ii).}

\smallskip

We remark that the simplicial complex $X(E)$ is not simply
connected in general. For instance take $G = {\rm GL}(4,F)$ and
$E/F$ quadratic unramified. Then $\XE$ is the building of $G_E =
{\rm GL}(2,F)$ which is $1$-dimensional. Using the criterion of
(I.3.5), we get that any vertex of $X$ belongs to $X(E)$. On the
other hand the barycenter of an edge in $X$ attached to a
$2$-periodical $\ofr_F$-lattice chain $(L_k )_{k\in \Zdss}$ in $V$
lies in $X(E)$ if and only if ${\rm dim}_{\Fdss}(L_k /L_{k+1}) =2$
for all $k$. Any given chamber of $X$ therefore has exactly two
opposite edges $\sigma_0$ and $\sigma_1$ that lie in $X(E)$. If we
consider all chambers in an apartment of $X$ which contain
$\sigma_0$ then the corresponding edges opposite to $\sigma_0$
form a cycle in $X(E)$.

\bigskip

\centerline{\bf II. Simple characters and their endo-classes}

\medskip

Here we recall some basic facts about simple characters. References
are to be found in [BK1] and [BH]. We continue to use the notation of
(I).

\smallskip

{\bf II.1 Simple pairs and their realizations.}

Recall that a simple pair $[0,\beta ]$ ([BH](1.5)) is a finite field
extension $E/F$, equipped with a generator $\beta$ (i.e. $E=F(\beta)$)
and satisfying the following conditions:

 (SP1) $ \beta\not\in\ofr_{E}$,

 (SP2) $k_o (\beta , \Afr (E) )< 0$ (cf. [BK1]{\S}1).

For each finite dimensional $E$-vector space $V$, and for each $\Bfr
 \in {\rm Her}(B)$, we have a simple stratum $[\Afr (\Bfr ), n_{\Bfr}
 ,0,\beta ]$ in $A$, called a {\it realization} of $[0,\beta ]$ ([BH]
 p. 133). Here $n_{\Bfr}$ is the valuation of $\beta\in A$ with
 respect to $\Afr (\Bfr )$.

Attached to $[\Afr (\Bfr ),n_{\Bfr}, 0,\beta ]$ (so to $[0,\beta]$,
 $V$ and $\Bfr$), we have the following data:

-- Two {\it open compact subgroups} of $G$:  $U^1(\Bfr)
\subseteq H^{1}(\Bfr )\subseteq
 J^{1}(\Bfr )\subseteq U^{1}(\Afr (\Bfr ))$; they are both normalized by
 ${\cal N}(\Bfr )$.

-- A finite set of {\it simple characters}  ${\cal C}(\Bfr ) = {\cal
 C}(\Afr (\Bfr ),0,\beta )$ of $H^{1}(\Bfr )$; each character in
 ${\cal C}(\Bfr )$ having a $G$-intertwining given by $J^{1}(\Bfr )G_E
 J^{1}(\Bfr )$.

-- Moreover, for each $\theta\in {\cal C}(\Bfr )$, there exists (up to
 isomorphism) a unique irreducible representation $\eta (\theta )$ of
 $J^{1}(\Bfr )$ such that $\eta (\theta )_{| H^{1}(\Bfr )}$ contains
 $\theta$. The intertwining of $\eta (\theta )$ is again  $J^{1}(\Bfr )G_E
 J^{1}(\Bfr )$ and the representation ${\rm Ind}_{H^{1}(\Bfr
 )}^{J^{1}(\Bfr )} \theta $ is a multiple of $\eta (\theta )$.

In addition we need the {\it degenerate} simple characters ([BK1] p.
184). To have a uniform notation we in this case set $E := F$ and $B
:= A$; for any $\Bfr \in {\rm Her}(B) = {\rm Her}(A)$ we let
$H^1(\Bfr) := J^1(\Bfr) := U^1(\Bfr)$ and let ${\cal C}(\Bfr )$ denote
the one element set consisting of the trivial representation ${\bf
1}_{H^1(\Bfr)}$ of $H^1(\Bfr)$.

If we need to keep track of the $E$-vector space $V$ we some times
write $H^1(V,\Bfr)$, $J^1(V,\Bfr), {\cal C}(V,\Bfr)$ instead of
$H^1(\Bfr), J^1(\Bfr), {\cal C}(\Bfr)$, respectively.

\smallskip

{\bf II.2 Potential simple characters} ( cf. [BH] {\S}8).

Let $[0,\beta ]$ be a simple pair and $V_1$, $V_2$ be two finite
 dimensional $E$-vector spaces. Write $B_i = {\rm End}_E V_i$, $A_i =
 {\rm End}_F V_i$, $i=1,2$. For $i=1,2$, fix a hereditary order
 $\Bfr_i\in {\rm Her}(B_i )$. Recall ([BK1] (3.6)) that we have a
 canonical bijection (called a transfer map):
$$
{\bf \tau}_{\Bfr_1 ,\Bfr_2 ,\beta}: {\cal C}(V_1,\Bfr_1 )
\longrightarrow {\cal C}(V_2,\Bfr_2 ).
$$
These transfer maps satisfy the properties:
$$
{\bf \tau}_{\Bfr_1 ,\Bfr_2 ,\beta} = {\bf \tau}_{\Bfr_2 ,\Bfr_1
,\beta}^{-1}\ , \ {\bf \tau}_{\Bfr_1 ,\Bfr_3 ,\beta} = {\bf
\tau}_{\Bfr_1 ,\Bfr_2 ,\beta}
 \circ  {\bf \tau}_{\Bfr_2 ,\Bfr_3 ,\beta}.
$$
Write $\Rfr [0,\beta ] =\bigcup {\cal C}(V,\Bfr )$, where ${\cal
C}(V,\Bfr )$ runs over the sets of simple characters attached to all
possible realizations of $[0,\beta ]$. We say that $\theta_1$,
$\theta_2 \in \Rfr [0,\beta ]$, attached to $(V_i ,\Bfr_i )$, $i=1,2$
are {\it equivalent} if $\theta_2 = {\bf \tau}_{\Bfr_1 ,\Bfr_2
,\beta}\theta_1$. This is indeed an equivalence relation and the
equivalence classes are called {\it potential simple characters} (or
{\it ps-characters}) supported by $[0,\beta ]$.

In addition we let all possible degenerate simple characters form a
single class which will be called the {\it degenerate ps-character}.

\smallskip

{\it Remarks.} (i) To be given a ps-characters amounts to fixing some
$\theta\in {\cal C}(V,\Bfr )$ in some realization.

(ii) A ps-character $\Theta$ may be seen as a {\it function} of the
pairs $(V, \Bfr )$: to $(V,\Bfr )$ we attach the simple character
$\theta \in \Theta$ that lies in ${\cal C}(V,\Bfr )$. We shall also
say that $\Theta (V, \Bfr )$ is a {\it realization of $\Theta$
associated to $(V,\Bfr )$}.

\smallskip

{\bf II.3 Endo-classes of ps-characters} (cf. [BH]{\S}8).

Let $\Theta_i$, for $i=1,2$,  be two ps-characters. Then each
$\Theta_i$ is either supported by a simple pair  $[0,\beta_i]$ (with
$E_i := F(\beta_i)$) or is degenerate (with $E := F$). We say that two
realizations $\Theta_1 (V_1 ,\Bfr_1 )$ and $\Theta_2 (V_2 ,\Bfr_2 )$
are {\it simultaneous} if $[E_1 :F] =[E_2 :F]$ and if the $F$-vector
spaces $V_1$ and $V_2$ are the same.

\smallskip

{\bf (II.3.1) Definition} ([BH](8.6)). {\it Two ps-characters
$\Theta_1$ and $\Theta_2$ are called endo-equivalent, denoted
$\Theta_1 \simeq \Theta_2$, if there exist simultaneous realizations
$\Theta_1  (V_1 ,\Bfr_1 )$ and $\Theta_2 (V_2 ,\Bfr_2 )$ that
intertwine in ${\rm Aut}_F V$, where $V=V_1 =V_2$. We shall summarize
this condition by saying that $\Theta_1$ and $\Theta_2$ intertwine in
some simultaneous realization.}

\smallskip

The following proposition shows that $\simeq$ is indeed an equivalence
relation.

\smallskip

{\bf (II.3.2) Proposition} (cf. [BK1](3.6) and [BH] pp. 154-157).
 (i) {\it If $\Theta$ is a ps-character then any pair of simultaneous
  realizations of $\Theta$
 intertwine.}

(ii) {\it If $\Theta_1$ and $\Theta_2$ are ps-characters, they
intertwine in some simultaneous realization if and only if they
intertwine in any  simultaneous realization.}

\smallskip

A class for $\simeq$ is called an {\it endo-class of ps-characters}.

We shall need the following two facts.

\smallskip

{\bf (II.3.3) Proposition} ([BH] (8.11)). {\it Let $\bf \Theta$ be an
endo-class of non-degenerate ps-characters and $\Theta \in {\bf
\Theta}$ supported by $[0,\beta ]$. Then the  following integers only
depend on $\bf
\Theta$: $k_o (\beta ,\Afr (E))$, $v_{E}(\beta )$, $e(E/F)$
(ramification index)  and $f(E/F)$ (inertial degree).}

\smallskip

{\bf (II.3.4) Proposition} ([BK1] (3.5.11)). {\it Let ${\bf
\Theta}$ be an endo-class of ps-characters and $\Theta_1$, $\Theta_2
\in {\bf \Theta}$. Let $\theta_1 = \Theta_1 (V_1 ,\Bfr_1 )$ and
$\theta_2 = \Theta_2 (V_2 ,\Bfr_2 )$ be simultaneous realizations.
Assume that $\Afr (\Bfr_1 ) =  \Afr (\Bfr_2 ) =:
\Afr$. Then there exists $x\in U(\Afr )$ such that $\theta_2
=\theta_1^{x}$. }

\bigskip

\centerline{\bf III. Ps-characters and pairs of orders}

\medskip

{\bf III.1 Extensions to mixed groups.}

We fix  a simple pair $[0,\beta ]$, a  ps-character $\Theta$ supported
by $[0,\beta ]$, as well as a finite dimensional $E$-vector space $V$.
We keep the notation as in (I) and (II).

The ps-character $\Theta$ gives rise to a function $\theta$; it maps
an order $\Bfr\in {\rm Her}(B)$ to the simple character
 $\theta (\Bfr ) = \Theta (V, \Bfr )$ of ${\cal C}(\Bfr )$. For each
$\Bfr\in {\rm Her}(B)$, let $\eta (\Bfr ) = \eta (V,\Bfr )$ be the
Heisenberg representation of $J^{1} (\Bfr )$ which contains
$\theta(\Bfr)$ when restricted to $H^{1}(\Bfr )$.

For each pair of hereditary orders $\Bfr_1 \subseteq \Bfr_2$ in ${\rm
 Her}(B)$, we have $U (\Bfr_1 )\subseteq U(\Bfr_2 )$ and $U^{1}(\Bfr_2
 )\subseteq U^{1}(\Bfr_1 )$. Since $U^{1}(\Bfr_1 ) \subseteq U(\Bfr_2 )$
 and $U(\Bfr_2 )$ normalizes $J^{1}(\Bfr_2 )$, one may form the
 group
$$
J^{1}(\Bfr_1 ,\Bfr_2 ):= U^{1}(\Bfr_1 )J^{1}(\Bfr_2 ) .
$$

\smallskip

{\bf (III.1.1) Proposition} ([BK1](5.1.14-16), (5.1.18), (5.1.19)).
{\it There exists a unique family of irreducible representations $\{
(J^{1}(\Bfr_{1} ,\Bfr_{2}), \eta (\Bfr_{1} ,
\Bfr_{2}))\}_{\Bfr_{1}\subseteq \Bfr_{2}}$ (determined up to isomorphism)
which extends $\{ \eta (\Bfr )\}_{\Bfr }$ in the following  sense:

  (i)  $\eta (\Bfr , \Bfr ) = \eta (\Bfr )$ for any $\Bfr$ in ${\rm Her}(B)$;

 (ii)  $ \eta (\Bfr_{1} , \Bfr_{2})_{|J^{1}(\Bfr_{2})} \simeq
\eta (\Bfr_{2})$ for all $\Bfr_{1}\subseteq \Bfr_{2}$ in ${\rm Her}(B)$;

 (iii) the following induced representations are irreducible and equivalent:
$$
{\rm Ind}_{J^{1}(\Bfr_{1})}^{U^{1}(\Afr (\Bfr_{1}))} \eta (\Bfr_{1})
\simeq {\rm Ind}_{J^{1}(\Bfr_{1},\Bfr_{2})}^{U^{1}(\Afr (\Bfr_{1}))}
\eta (\Bfr_{1} ,\Bfr_{2})\ .
$$
Moreover we have:

 (iv) The compatibility condition:
 $\eta (\Bfr_{1},\Bfr_{3})_{|J^{1}(\Bfr_{2},\Bfr_{3})}
 = \eta (\Bfr_{2} ,\Bfr_{3})$,
 for any triple $\Bfr_{1}\subseteq \Bfr_{2} \subseteq \Bfr_{3}$
 in ${\rm Her}(B)$;

 (v) the intertwining formula:
$$
{\cal I}_{G}(\eta (\Bfr_{1} ,\Bfr_{2})) = J^{1}(\Bfr_{2})G_E
J^{1}(\Bfr_{2})\ .
$$}

\smallskip

Note that the representation
$$
\eta (\Afr (\Bfr )) := \eta (V,\Afr (\Bfr )) :=
{\rm Ind}_{J^{1}(\Bfr )}^{U^{1}(\Afr (\Bfr ))} \eta (\Bfr )
$$
is irreducible for all $\Bfr$. Its intertwining is given by
$$
{\cal I}_{G}(\eta (\Afr (\Bfr )) = U^{1}(\Afr (\Bfr )) G_E  U^{1}(\Afr
(\Bfr ))\ .
$$

\smallskip

{\bf (III.1.2) Proposition}. {\it For all $g$ in $G_E$ and
$\Bfr_{1}\subseteq \Bfr_{2}$, we have $[J^{1} (\Bfr_{1},
\Bfr_{2})]^{g} = J^{1}(\Bfr_{1}^{g}, \Bfr_{2}^{g})$ and the
representations $\eta (\Bfr_{1} ,\Bfr_{2})^g$ and $\eta (\Bfr_{1}^g ,
\Bfr_{2}^g)$ are isomorphic.}

\smallskip

First we need the following result. Let $V'$ denote the $F$-vector
space $V$ equipped with a possibly different $E$-vector space
structure. We then find an element $x \in G$ such that $B' := {\rm
End}_E V'$ satisfies $B' = B^x := xBx^{-1}$ and hence ${\rm Her}(B') =
\{\Bfr^x : \bfr \in {\rm Her}(B)\}$.

\smallskip

{\bf (III.1.3) Lemma}. {\it a)  For any $\Bfr \in {\rm Her}(B)$ we
have:

 i) $H^{1} (V',\Bfr^{x}) = H^{1}(V,\Bfr)^{x}$ and
 $J^{1} (V',\Bfr^{x}) = J^{1}(V,\Bfr)^{x}$.

 ii) If $\theta \in \Cscr (V,\Bfr)$, then
$\theta^{x} \in \Cscr (V',\Bfr^x)$.

b) For $g \in G_{E}$ and $\Bfr \in {\rm Her}(B)$ we have $\theta (\Bfr
)^{g} = \theta (\Bfr^g )$. }

{\it Proof.} The point a) follows immediately from the inductive
definition of simple characters and groups (cf. [BK1] {\S}3).

We need to recall the characterization of the transfer maps
 $\tau_{\Bfr_1 ,\Bfr_2 ,\beta}$  for a pair of orders $\Bfr_i$,
 $i=1,2$, in ${\rm Her}(B)$ ([BK1](3.6)): If $\theta_i \in {\cal
 C}(\Bfr_i )$, $i=1, 2$, then $\theta_2 = \tau_{\Bfr_1 ,\Bfr_2 ,\beta}
\theta_1$ if and only if $1\in G_E$ intertwine $\theta_1$ and
 $\theta_2$.

Consider the two characters $\theta (\Bfr^{g})$ and $\theta (\Bfr )^g$
of $H^{1}(\Bfr^{g})$.  Since $g$ intertwines $\theta (\Bfr ) $,
 we must have $\theta (\Bfr )^{g}_{| H^{1}(\Bfr )\cap H^{1}(\Bfr^g )}
= \theta (\Bfr )_{| H^{1}(\Bfr )\cap H^{1}(\Bfr^g )}$.
 So $\theta (\Bfr )^g \in {\cal C}(\Bfr^g )$ coincides with
$\tau_{\Bfr ,\Bfr^{g},\beta}(\theta (\Bfr ))$, that is with $\theta
(\Bfr^g)$ by definition of $\Theta$.

\smallskip

Turning to the proof of Proposition (III.1.2) the $G_E$-equivariance
of the family $\{ (J^{1}(\Bfr_{1} ,\Bfr_{2}), \eta (\Bfr_{1}
,\Bfr_{2}))\}_{\Bfr_{1}\subseteq \Bfr_{2}}$ follows now from that of
$\{\theta (\Bfr )\}_{\Bfr \in {\rm Her}(B)}$ by a unicity argument.

\smallskip

{\bf III.2 Extensions to $1$-units of orders.}

We now quote some properties of the representations $\eta (\Afr (\Bfr
 ))$. In the following we abbreviate $\Afr (\Bfr_{*}) = \Afr_{*}$
 for any subscript ``$*$". Let $\Bfr_{1}\subseteq \Bfr_{2}$ be hereditary orders
in $B$.

We first note that
$$
J^{1}(\Bfr_{1} ,\Bfr_{2}) \subseteq U^{1}(\Bfr_{1})U^{1}(\Afr_{2})
\subseteq U^{1}(\Afr_{1})\ .
$$
So we can consider the irreducible representation
$$
\eta (\Afr_{1} ,\Afr_{2}) :=
 {\rm Ind}_{J^{1}(\Bfr_{1}
 ,\Bfr_{2})}^{U^{1}(\Bfr_{1})U^{1}(\Afr_{2})}
\eta (\Bfr_{1} ,\Bfr_{2})\ .
$$

\smallskip

{\bf (III.2.1) Proposition} {\it a) The representation $\eta (\Afr_{1}
,\Afr_{2} )$  satisfies

 (i) $\eta (\Afr_{1} ,\Afr_{2})_{|U^{1}(\Afr_{2})} \simeq \eta (\Afr_{2})$;

 (ii) ${\rm Ind}_{U^{1}(\Bfr_{1})U^{1}(\Afr_{2})}^{U^{1}(\Afr_{1})}
\eta (\Afr_{1} , \Afr_{2}) \simeq \eta (\Afr_{1} )$.

b) Moreover for any triple of hereditary orders
 $\Bfr_{1}\subseteq \Bfr_{2} \subseteq \Bfr_{3}$, we have
$$
\eta (\Afr_{1} ,\Afr_{3})_{|U^{1}(\Bfr_{2})U^{1}(\Afr_{3})}
 = \eta (\Afr_{2} ,\Afr_{3})\ .
$$ }

\smallskip

{\it  Proof.} Assertion a) (ii) is a consequence of Proposition
 (III.1.1)(iii). We must prove b). By Mackey's restriction formula
and since the double quotient
$$
U^{1}(\Bfr_{2})U^{1}(\Afr_{3})\backslash
U^{1}(\Bfr_{1})U^{1}(\Afr_{3})/U^{1}(\Bfr_{1}) J^{1}(\Bfr_{3})
$$
is reduced to one element, we get that the restriction in b) is
$$
{\rm Ind}_{U^{1}(\Bfr_{1})J^{1}(\Bfr_{3})
\cap U^{1}(\Bfr_{2})U^{1}(\Afr_{3})}^{U^{1}(\Bfr_{2})U^{1}(\Afr_{3})}
 \eta (\Bfr_{1} ,\Bfr_{3}) = {\rm
Ind}_{U^{1}(\Bfr_{2})J^{1}(\Bfr_{3})}^{U^{1}
(\Bfr_{2})U^{1}(\Afr_{3})} \eta (\Bfr_{1} ,\Bfr_{3})\ .
$$
Now the result follows from Proposition (III.1.1)(iv).

\smallskip

{\bf III.3 The degenerate case.}

The constructions of (III.1) and (III.2) trivially extend to the case
where $\Theta$ is the degenerate ps-character. Indeed, in that case,
we set $E=F$ and for all pairs of orders $\Bfr_1 \subset \Bfr_2$ in
${\rm Her}(B)$, we set:
\smallskip

-- $J^{1}(\Bfr_1 ,\Bfr_2 )= U^{1}(\Bfr_1 ) U^{1}(\Bfr_2 ) =
   U^{1}(\Afr_1 )U^{1}(\Afr_2 )= U^{1} (\Afr_1 )$;

-- $\eta (\Bfr_1 ,\Bfr_2 ) = \eta (\Afr_1 ,\Afr_2 ) ={\bf 1}_{U^1
   (\Afr_1 )}$.

\bigskip

\centerline{\bf IV. The coefficient system on ${\rm sd}(X_E)$}

\medskip

As in the previous section, we fix a
 ps-character $\Theta$. It is either degenerate ($E := F$)
 or supported by a simple pair $[0, \beta]$ ($E := F(\beta)$).
 We also fix a finite dimensional $E$-vector space $V$.

Let $\VV$ be a smooth complex
 representation of $G = {\rm Aut}_{F}V$. In a first step, we are going
 to construct a $G_E$-equivariant coefficient system $\CC_o (\VV ) =
 \CC_o (\Theta ,V,\VV )$ on ${\rm sd}(X_E )$. We shall first construct
 this coefficient system on the stars of the vertices of $X_E$ and then
 extend it to any simplex.

We call a simplex $\sigma = (\Bfr_0 \supset \dots \supset \Bfr_q )$
 {\it semistandard} if it belongs to the star of some vertex in $X_E$,
 that is if $\Bfr_0$ is a maximal order.

\smallskip

{\bf (IV.1) Definition}. {\it  i) For any semistandard simplex
 $\sigma = (\Bfr_0 \supset \dots \supset \Bfr_q )$ of ${\rm sd}(X_E
 )$, we set
$$
\VV(\sigma ) = \VV^{\eta (\Bfr_q ,\Bfr_0 )}\ ,
$$
the $\eta (\Bfr_q ,\Bfr_0 )$-isotypic component of $\VV$.

 ii) For an arbitrary simplex $\sigma$ of ${\rm sd}(X_E )$, we set
$$
\VV (\sigma ) = \sum_{\tau\ {\rm semistandard},\tau \supseteq
 \sigma} \VV (\tau )\  .
$$}

\smallskip

{\bf (IV.2) Proposition}. {\it i) The previous definition is
consistent.

 ii) For any pair of simplices $\sigma, \tau$ of $\XE$ with
$\sigma \subseteq \tau$, we have $\VV (\tau )\subseteq \VV (\sigma
)$.}

{\it Proof.} We only need to prove the second assertion in the case of
semistandard simplices $\sigma,\tau$. Suppose therefore that $\sigma =
(\Bfr_0 \supset \dots \supset \Bfr_q )$ is semistandard. The stars of
two distinct vertices being disjoint, the simplex $\tau$ must then
have the form $\tau = (\Bfr_0 \supset \dots \supset
\Bfr_r )$ containing $(\Bfr_0 \supset \dots \supset \Bfr_q )$ as a subflag.
By (III.1.1)(iv), we have $\eta (\Bfr_r ,\Bfr_0 ){| J^{1}(\Bfr_q
,\Bfr_0 )} = \eta (\Bfr_q ,\Bfr_0 )$. So $\VV^{\eta (\Bfr_r ,
\Bfr_0 )}\subseteq  \VV^{\eta (\Bfr_q ,\Bfr_0 )}$ and the result
follows.

\smallskip

By taking inclusions as transition maps, the family $\CC_o (\VV ) =
(\VV (\sigma ))_{\sigma }$, $\sigma$ running over the simplices of
$\XE$, is then a coefficient system of $\Cdss$-vector spaces over
$\XE$.

\smallskip

{\bf (IV.3) Proposition.} {\it For the obvious action of $G_E$ on the
$\VV (\sigma )$, $\sigma$ simplex of $\XE$, $\CC_o (\VV )$ is
naturally endowed with a structure of $G_E$-equivariant coefficient
system.}

{\it Proof}. We must prove that $g\VV (\sigma ) = \VV (g\sigma )$, for
all $g\in G_E$ and $\sigma$ simplex of $\XE$. Also we may clearly
reduce to the case where $\sigma$ and $\tau$ are semi-standard.

Let $\sigma = (\Bfr_0 \supset \dots \supset \Bfr_q ) $ be semistandard
 and $g$ be in $G_E$. Then $g\sigma = (\Bfr_{0}^{g} \supset \dots
 \supset \Bfr_{q}^{g})$ and
$$
g \VV (\sigma ) =g\VV^{\eta (\Bfr_q ,\Bfr_0 )}\ {\rm and } \
 \VV (g\sigma )  = \VV^{\eta (\Bfr_{q}^g ,\Bfr_{0}^g )} \ .
$$
 By (III.1.2), this last vector space is $\VV^{ \eta (\Bfr_q ,\Bfr_0
)^g}$. Now our result follows from the following observation. Let
$(K,\rho )$ be a smooth irreducible representation of a compact open
subgroup $K$ of $G$ and let $g\in G$. Then $g\VV^{\rho} =
\VV^{\rho^g}$, where $\rho^g$ is the representation of $K^g = gKg^{-1}$
given by $\rho^g (k) = \rho (g^{-1}kg)$.

\bigskip

\centerline{\bf V. The coefficient system on ${\rm sd}(X)$}

\medskip

We keep the notations from the previous sections. As in (I) we see
$\XE$ as a subcomplex of $\XX$. We are now going to construct a
coefficient system $\CC (\VV ) = \CC (\Theta , V, \VV )$ on $\XX$.

For any subscript ``$*$'', we shall write $\Afr_{*}$ for $\Afr
 (\Bfr_* )$. In particular, if $(\Bfr_0 \supset \dots\supset \Bfr_q )$ is a
 flag of orders in ${\rm Her}(B)$ then $(\Bfr_0 \supset \dots\supset \Bfr_q
 )$, $(\Afr (\Bfr_0 ) \supset \dots\supset \Afr (\Bfr_q ) )$ and $(\Afr_0
 \supset \dots \supset \Afr_q )$ denote the same object, i.e. a
 simplex of $\XE$ seen as a simplex of $\XX$.

We shall need the two following lemmas.

\smallskip

{\bf (V.1) Lemma}. {\it Let $(\rho,W)$ be a smooth irreducible
representation of some compact open subgroup $K \subseteq G$ and let
$\VV^{\rho}$ denote the $\rho$-isotypic component of $\VV$. Then
$\VV^{\rho}$ is invariant under any subgroup of $N_G(K)$ which
intertwines $\rho$.}

{\it Proof.} Let $v\in \VV^{\rho}$ and $g\in G$ be an element
normalizing $K$ and intertwining $\rho$.  By definition, there exist
$\varphi$ in ${\rm Hom}_{\rho}(W,\VV )$ and $w\in W$ such that $v
=\varphi (w)$. Since  $gKg^{-1}=K$ and $\rho^g \simeq \rho$, and since
$\rho$  is irreducible, there must exist an intertwining operator
$\psi\in {\rm Aut}_{\Cdss}(W)$ such that $\rho^g (k ) = \psi^{-1}\circ
\rho (k) \circ \psi$, for all $k\in K$. It easily follows that $g\varphi
\psi^{-1}$ belongs to ${\rm Hom}_{\rho}(W, \VV )$. So $gv = [g\varphi
\psi^{-1}] (\psi (w))$ and $gv\in \VV^{\rho}$, as required.

\smallskip

{\bf (V.2) Lemma.} {\it  Let $H\subseteq K$ be compact open subgroups
of $G$. Let $\eta_{H}$ be an irreducible smooth representation of $H$
and assume that $\eta_K := {\rm Ind}_{H}^{K}\eta_H$ is irreducible as
well. Then $\VV^{\eta_K} = K\VV^{\eta_{H}}$.}

{\it Proof}. Let $\Phi\in {\rm Hom}(\eta_K ,\VV )$. We may  see
$\eta_H$ as an $H$-submodule of $\eta_K$ so that
$$
\eta_K =\bigoplus_{k\in K/H} k \eta_H \ .
$$
So
$$
\Phi (\eta_K ) = \sum_{k\in K/H} k\Phi (\eta_H )\ ,
$$
with $\Phi (\eta_H )$ contained in $\VV^{\eta_H}$ since $\Phi$ is
$H$-equivariant. This gives the inclusion $\VV^{\eta_K} \subseteq
K\VV^{\eta_{H}}$ Conversely, since the smooth representations of $H$
are semisimple, $\VV^{\eta_H}$ decomposes into a direct sum
$$
\VV^{\eta_H} = \bigoplus_{i\in I}\VV_i\ ,
$$
each $\VV_i$ being isomorphic to $\eta_H$ as an $H$-module. Now each
$K.\VV_i \subseteq \VV$ is isomorphic to $\eta_K$ as a $K$-module and
the opposite inclusion follows.

\smallskip

{\bf (V.3) Definition}. {\it For $\sigma$ a semistandard simplex in
$\XE$, we set
$$
\VV_{\sigma} = \sum_{g\in {\rm Stab}_{G}(\sigma )} g\VV (\sigma
)\subseteq \VV \ .
$$}

\smallskip

{\bf (V.4) Lemma.} {\it Let  $\sigma =(\Bfr_0 \supset \dots
\supset \Bfr_q )$ be a semistandard simplex of $\XE$. Then
$$
{\rm Stab}_{G}(\sigma ) = E^{\times}U(\Afr_q ) \ .
$$}
{\it Proof.} The group $E^{\times}U(\Afr_q )$ certainly normalizes
$(\Bfr_0 \supset \dots
\supset \Bfr_q )$ $=$ $(\Afr_0 \supset \dots
\supset \Afr_q )$ and lies in ${\rm Stab}_{G}(\sigma )$. Conversely if
$g\in {\rm Stab}_{G}(\sigma )$, then $g$ normalizes the principal
order $\Afr_0$ and must lie in its stabilizer which by (I.3.1) is
equal to $E^{\times}U(\Afr_0 )$. Write $g =\lambda h$, with
$\lambda\in E^{\times}$ and $h\in U(\Afr_0 )$. Since $\lambda$ is in
${\cal N} (\Afr_q ) = {\cal N}(\Bfr_q )U(\Afr_q )$, so is $h$. Now $h$
must be in the maximal compact subgroup of $ {\cal N}(\Bfr_q )U(\Afr_q
)$, that is  $U(\Afr_q )$, and the lemma follows.

\smallskip

{\bf (V.5) Proposition.} {\it Let $\sigma =(\Bfr_0 \supset \dots
\supset \Bfr_q )$ be a semistandard simplex of $\XE$. Then:
$$
\VV_{\sigma} = \sum_{g\in U(\Afr_q )/U(\Bfr_q )J^{1}(\Bfr_0
)}g\VV (\sigma ) = \sum_{g\in U( \Afr_q )/U^{1}(\Afr_q )}g\VV^{\eta
(\Afr_q )}\ .
$$}
{\it Proof}.  The subgroup $U(\Bfr_q)J^{1}(\Bfr_0 )$ normalizes
$J^{1}(\Bfr_q ,\Bfr_0 ) = U^{1}(\Bfr_q )J^{1}(\Bfr_0 )$. Moreover it
intertwines $\eta (\Bfr_q ,\Bfr_0 )$ by (III.1.1)(v). As a consequence
of (V.1), (V.4), and the definition of $\VV (\sigma )$ we therefore
obtain the first equality in
$$
\VV_{\sigma} = \sum_{E^{\times}U(\Afr_q )/U(\Bfr_q )J^{1}(\Bfr_0)}
g\VV(\sigma ) = \sum_{U(\Afr_q )/U^{1}(\Bfr_q )J^{1}(\Bfr_0)}
g\VV(\sigma ) \ .
$$
The second one is immediate from the fact that $E^{\times}$ stabilizes
$\VV(\sigma )$. Now, using (III.1.1)(iii), we may apply (V.2) with $H
=J^{1}(\Bfr_q ,\Bfr_0 )$, $K = U^{1}(\Afr_q )$, $\eta_H = \eta (\Bfr_q
,\Bfr_0 )$, $\eta_K =
\eta (\Afr_q )$ to get the second equality in the proposition.

\smallskip

{\bf (V.6) Proposition}. {\it Let $\sigma$ and $\tau$ be semistandard
simplices of $\XE$ with $\sigma \subseteq \tau$. Then $\VV_{\tau}
\subseteq \VV_{\sigma}$.}

{\it Proof}. Write $\sigma = (\Bfr_0 \supset \dots \supset \Bfr_q )$
and $\tau = (\Bfr_0 \supset \dots \supset \Bfr_r )$, with $\Bfr_0$
maximal. By (IV.2)(ii), we have $\VV (\tau )\subseteq \VV (\sigma )$.
Moreover $U(\Afr_r )\subseteq U(\Afr_q )$. Our inclusion follows now
from the first equality in (V.5).

\smallskip

{\bf (V.7)  Proposition}. {\it Let $\sigma$ and $\tau$ be semistandard
simplices of $\XE$, and assume that $\tau = g\sigma$ for some $g\in
G$. Then $g\VV_{ \sigma} = \VV_{\tau}$.}

{\it Proof.} Write $\s = (\Bfr_{0}^{\sigma} \supset \dots \supset
 \Bfr_{q}^{\sigma})$ and $\tau = (\Bfr_{0}^{\tau} \supset
\dots \supset \Bfr_{q}^{\tau})$. Using (I.3.3), we can decompose $g$ as
 $g_{E}g_{\sigma}$, $g_{E}\in G_{E}$, $g_{\sigma}\in {\rm Stab}_{G}(\sigma
 )$.
 By construction we have $g_{\sigma}V_{\sigma} = V_{\sigma}$.
 So $g V_{\sigma} = g_{E}V_{\sigma}$.  We get:
$$
gV_{\sigma} = \sum_{g\in U(\Afr (\Bfr_{q}^{\sigma}))  } g_E g
g_{E}^{-1} g_{E}V(\sigma )\ .
$$
By (IV.3), we have $g_{E}V(\sigma ) =V(\tau )$, and it follows that
$$
g V_{\sigma} =  \sum_{h\in U(\Afr (\Bfr_{q}^{\sigma}))^{g_{E}}}
hV(\tau )\ .
$$
Now the result follows from the $G_E$-equivariance of the map $\Bfr
\mapsto U(\Afr (\Bfr ))$ and from the definition of $V_{\tau}$.

\smallskip

{\bf (V.8) Definition}. {\it A simplex of $\XX$ is called}
$E$-semistandard {\it if it is conjugate to a semistandard simplex of
$\XE$. We define a vector space $\VV_{\sigma}$, for each simplex
$\sigma$ of $\XX$, as follows:

i) If $\sigma = g\tau$, for $\tau$ semistandard in $\XE$ and $g\in
G$ , then $\VV_{\sigma}= g \VV_{\tau}$;

ii) If $\sigma$ is an arbitrary simplex of $\XX$, then
$$
\VV_{\sigma} = \sum_{\tau \ E{\rm -semistandard},\tau\supseteq \sigma}
\VV_{\tau} \ .
$$}

\smallskip

{\bf (V.9) Proposition}. {\it  i) The previous definition is
consistent.

ii) For any pair of simplices $\sigma \subseteq \tau$ of $\XX$, we
have $\VV_{\tau} \subseteq \VV_{\sigma}$. In particular, by taking
inclusions as transition maps, the collection $\CC (\VV ):=
(\VV_{\sigma} )_{\sigma}$  is  a coefficient system of $\Cdss$-vector
spaces on $\XX$.

iii) For the obvious action of $G$, the coefficient system $\CC (\VV
)$ is equivariant.}

{\it Proof.} i) By (V.7) the definition of $\VV_{\sigma}$  in (V.8)(i)
does not depend on the choice of $g$.

To prove that the definition (V.8)(ii) is consistent, we must prove
that if $\sigma$ and $\tau$ are $E$-semistandard simplices of $\XX$
satisfying $\tau \supseteq \sigma$, then $\VV_{\tau}\subseteq
\VV_{\sigma}$. Write $\sigma = g\sigma_o$, $\tau =h\tau_o$, with
$g,h\in G$ and $\sigma_o = (\Bfr_0\supset \dots
\supset \Bfr_q )$, $\tau_o = (\Cfr_0 \supset \dots \supset \Cfr_r)$
semistandard in $\XE$. By definition $\VV_{\sigma}= g\VV_{\sigma_o}$
and $\VV_{\tau}= h\VV_{\tau_o}$.

The hypothesis $\tau\supseteq \sigma$ implies $\Bfr_i =
g^{-1}h\Cfr_{j(i)}$ for some $j(0)=0 < \ldots < j(q) \leq r$. By (V.6)
and (V.7), we have
$$
\VV_{\tau} = h\VV_{\tau_o} = h\VV_{(\Cfr_0 \supset\dots \supset \Cfr_r)}
\subseteq h\VV_{(\Cfr_{j(0)} \supset\dots \supset \Cfr_{j(q)})} =
hh^{-1}g\VV_{\sigma_o} = g\VV_{\sigma_o} = \VV_{\sigma} \ .
$$

ii) We can obviously reduce to the case where $\sigma$ and $\tau$ are
$E$-semistandard, and the inclusion has just been proved in i).

iii) We must simply prove that $g\VV_{\sigma} = \VV_{g\sigma}$, for
any simplex $\sigma$ of $\XX$ and $g\in G$. We may reduce to the case
where $\sigma$ is $E$-semistandard where the result follows trivially
from the definition of $\CC (\VV)$.

\smallskip

By construction the coefficient system $\CC (\VV )$ is supported
on
$$
X(E) = \bigcup_{g\in G} g \XE
$$
viewed naturally as a simplicial subcomplex of $\XX$. But $\CC
(\VV )$ has the following additional constancy property.

\smallskip

{\bf (V.10) Proposition}. {\it Let the vertex $\sigma_o = (\Bfr)$
in $\XE$ be the barycenter of a simplex $\widetilde{\sigma}$ of
$X_E$; then $\VV_{\sigma} = \VV_{\sigma_o}$ for any simplex
$\sigma$ in $\XE$ such that $\sigma_o \subseteq \sigma \subseteq
\widetilde{\sigma}$.}

{\it Proof.} Put $\Afr := \Afr(\Bfr)$ and
$$
\VV_o := \sum_{g\in U( \Afr)/U^{1}(\Afr)}g\VV^{\eta (\Afr)}\ .
$$
If $\sigma$ is semistandard then $\VV_{\sigma} = \VV_o$ by (V.5).
Consider therefore the case that $\sigma = (\Bfr_0 \supset \dots
\supset \Bfr_q )$ with $\Bfr_q = \Bfr$ is not semistandard, and
let $\tau$ be any $E$-semistandard simplex in $\XX$ such that
$\tau \supseteq \sigma$. We have to show that $\VV_{\tau}
\subseteq \VV_o$. Write $\tau = g\tau_o$ with $g \in G$ and
$\tau_o = (\Cfr_0 \supset \dots \supset \Cfr_r)$ semistandard in
$\XE$. By (I.3.3) we may assume that $\tau_0 \supseteq
g^{-1}\sigma = \sigma$. We then have $\Bfr_i = \Cfr_{j(i)}$ for
some $ 0 \leq j(0) < \ldots < j(q) \leq r$. Since $\tau_0$ is
semistandard whereas $\sigma$ is not the order $\Cfr_0$ is maximal
but $\Bfr_0$ is not. This means that $0 < j(0)$. It follows that
$\tau_1 := (\Cfr_0 \supset \Bfr_0 \supset \dots \supset \Bfr_q )$
is a semistandard simplex in $\XE$ such that $\tau_0 \supseteq
\tau_1 \supseteq \sigma$. Hence $\tau = g\tau_0 \supseteq g\tau_1
\supseteq g\sigma = \sigma$. Since $\tau \supseteq g\tau_1$ both
are $E$-semistandard we know from the proof of (V.9)(i) that
$\VV_{\tau} \subseteq \VV_{g\tau_1}$. On the other hand, by
(I.3.1) we may write $g = hg'$ with $h \in \Nscr(\Bfr) \subseteq
G_E$ and $g' \in U(\Afr)$. We obtain that $g\tau_1 = h\tau_1$ in
fact is semistandard in $\XE$. Using (V.5) we conclude that
$\VV_{\tau} \subseteq \VV_{g\tau_1} = \VV_{h\tau_1} = \VV_o$.

\smallskip

This result is best expressed in the following way. The simplicial
structure of $X_E$ (before subdivision) can be described in terms
of the simplicial structure of $X$ as follows: The interiors
$\widetilde{\sigma}^{\rm o}$ of simplices $\widetilde{\sigma}$ of
$X_E$ are precisely the (nonempty) subsets of the form
$\widetilde{\tau}^{\rm o} \cap X_E$ for some simplex
$\widetilde{\tau}$ of $X$.

Suppose now that $\widetilde{\sigma}_1,\widetilde{\sigma}_2$ are
simplices of $X_E$ such that
$$
g(\widetilde{\sigma}_1)^{\rm o} \cap (\widetilde{\sigma}_2)^{\rm
o} \neq \emptyset \qquad \hbox{for some}\ g \in G\ .
$$
Write $(\widetilde{\sigma}_i)^{\rm o} = (\widetilde{\tau}_i)^{\rm
o} \cap X_E$ with simplices $\widetilde{\tau}_i$ of $X$. Then
$g(\widetilde{\tau}_1)^{\rm o} \cap (\widetilde{\tau}_2)^{\rm o}
\neq \emptyset$ and hence $g\widetilde{\tau}_1 =
\widetilde{\tau}_2$ since the $G$-action on $X$ is simplicial. In
particular, $g$ maps the barycenter of $\widetilde{\tau}_1$ into
the barycenter of $\widetilde{\tau}_2$. Since both barycenters lie
in $X_E$ we conclude from (I.3.3) that there also is an element
$g_E \in G_E$ which maps the first barycenter into the second one.
It follows that $g_E\widetilde{\sigma}_1 = \widetilde{\sigma}_2$
and $g_E\widetilde{\tau}_1 = \widetilde{\tau}_2$. Hence
$gg_E^{-1}$ fixes $\widetilde{\tau}_2$ and, by (I.3.1) applied to
its barycenter, can be written $gg_E^{-1} = h_Eh$ with $h_E \in
G_E$ fixing $\widetilde{\sigma}_2$ and $h \in G$ fixing
$\widetilde{\tau}_2$ pointwise. We obtain
$$
g(\widetilde{\sigma}_1)^{\rm o} =
gg_E^{-1}(\widetilde{\sigma}_2)^{\rm o} =
h_Eh(\widetilde{\sigma}_2)^{\rm o} =
h_E(\widetilde{\sigma}_2)^{\rm o} = (\widetilde{\sigma}_2)^{\rm
o}\ .
$$
>From this fact one deduces in a straightforward way that $X(E)$
carries a simplicial structure where the simplices are the subsets of
the form $g\widetilde{\sigma}$ for $\widetilde{\sigma}$ a simplex of
$X_E$ and $g \in G$. We write $X[E]$ for $X(E)$ equipped with this
simplicial structure. Similarly as for $X_E$ the interiors of
simplices of $X[E]$ are the nonempty intersections
$\widetilde{\tau}^{\rm o} \cap X[E]$ for $\widetilde{\tau}$ running
over the simplices of $X$. The barycentric subdivision of $X[E]$ is
$X(E)$.

The Proposition (V.10) (together with $G$-equivariance) says that $\CC
(\VV )$ in an obvious way derives from an equivariant coefficient
system $\CC[\VV] = \CC[\Theta,V,\VV] := (\VV[\sigma])_\sigma$ on
$X[E]$ given by
$$
\VV[\sigma] := \sum_{g\in U(\Afr )/U^{1}(\Afr)} hg\VV^{\eta (\Afr )}
$$
for any simplex $\sigma$ of $X[E]$ where $\sigma$ is the image under
some $h\in G$ of a simplex of $X_E$ with barycenter $(\Bfr )$ and
where $\Afr :=\Afr (\Bfr )$.

\bigskip

\centerline{\bf VI. The degenerate case}

\medskip

We recall that in the degenerate case we have $E = F$, $B = A$,
$H^1(\Afr) = J^1(\Afr) = U^1(\Afr)$, and $\theta(\Afr) =
\eta(\Afr) = {\bf 1}$. The coefficient system $\CC
(\VV) = (\VV_\sigma)_\sigma$ on $X$ associated with a smooth
complex representation $\VV$ of $G$ is given by the fixed vectors
$$
\VV_\sigma = \VV^{U^1(\Afr)} \qquad \hbox{for}\ \sigma = F(\Afr).
$$
This is precisely one of the coefficient systems considered in [SS]
(namely the one corresponding to ``level'' $n = 1$). From loc. cit. we
therefore have the following result.

\smallskip

{\bf (VI.1) Theorem.} {\it The oriented chain complex of $\CC (\VV
)$ is an exact resolution in the category of all smooth complex
$G$-representations of the subrepresentation of $\VV$ generated
(as a $G$-representation) by $\VV^{U^1(\Afr_0)}$ for some vertex
$F(\Afr_0)$. Moreover, if the center of $G$ acts on $\VV$ through
a character $\chi$ then this resolution is a projective resolution
in the category of smooth $G$-representations with central
character $\chi$.}

\bigskip

\centerline{\bf VII. Dependence on the endo-class}

\medskip

 We fix a finite dimensional $F$-vector space $V$. Let $\Theta_i$,
 for $i=1,2$, be two ps-characters with simultaneous realizations in
 $A={\rm Aut}_F V$. So for each $i$ we are in one of the following cases:

 1) The ps-character $\Theta_i$ is supported by a simple pair
    $[0,\beta_i ]$ and $V$ is an $E_i$-vector
    space, where $E_i = F(\beta_i)$, and as such will be denoted by $V_i$.
    Following previous notations we have the
    centralizer $B_i$ of $E_i$ in $A$, the centralizer $G_{E_i}$ of
    $E_{i}^{\times}$ in $G = A^{\times}$, the affine
    building $X_{E_i}$ of $G_{E_i}$, and $\displaystyle X[E_i ]=\bigcup_{g\in G}
    gX_{E_i}$ equipped with the simplicial structure defined in (V).

 2) The degenerate case.

\smallskip

{\bf (VII.1) Lemma}. {\it  If $\Theta_1$ and $\Theta_2$ are
endo-equivalent then $X[E_1 ]=X[E_2 ]$ as sets and simplicial
complexes.}

{\it Proof}. It suffices to prove that the barycentric subdivisions
$X(E_1)$ and $X(E_2)$ coincide as simplicial subcomplexes of ${\rm
sd}(X)$. Being the fixed points sets of groups acting simplicially on
${\rm sd}(X)$ they are full subcomplexes. Hence they coincide if they
have the same vertices. For $i=1,2$, a vertex $(\Afr )$ of ${\rm
sd}(X)$ lies in $X(E_i )$ if and only if the order $\Afr\in {\rm
Her}(A)$ satisfies the numerical criterion of (I.3.5). But by
(II.3.3), since $\Theta_1$ and $\Theta_2$ are endo-equivalent, we have
$f(E_1 /F)=f(E_2 /F)$, $e(E_1 /F) =e(E_2 /F)$ and the numerical
criteria for $X(E_1)$ and $X(E_2)$ are the same.

\smallskip

{\bf (VII.2) Proposition}. {\it Let $\VV$ be a smooth representation
of $G$. If $\Theta_1$ and $\Theta_2$ are endo-equivalent, then the two
coefficient systems $\CC [\Theta_i,V_i,\VV]$ coincide.}

{\it Proof}. If $\Theta_1$ and $\Theta_2$ are endo-equivalent, then
they are both supported by simple pairs or are both degenerate. In
this latter case the coefficient systems in question coincide for
trivial reasons. So we may assume in the following that $\Theta_i$,
for $i=1,2$, is supported by a simple stratum $[0,\beta_i ]$ and use
the notations of case 1). Recall that the ps-character $\Theta_i$
gives rise to the simple character valued function $\theta_i :=
\Theta_i(V_i,.)$ on ${\rm Her}(B_i)$. We write $\eta_i(V_i,.)$ for the
representations corresponding to $\theta_i$ which were introduced in
(III.1).

We now fix a simplex $\sigma$ of $X[E_1 ]=X[E_2 ]$ and write $\sigma
=h_i \sigma_i$ with $h_i\in G$ and $\sigma_i$ a simplex of $X_{E_i}$.
Let $(\Bfr_i)$, for $\Bfr_i \in {\rm Her}(B_i)$, be the barycenter of
$\sigma_i$ and put $\Afr_i := \Afr(\Bfr_i) \in {\rm Her}(A)$. We have
to show that the identity
$$
h_1(\sum_{g\in U(\Afr_1)/U^{1}(\Afr_1)} g\VV^{\eta_1(V_1,\Afr_1)}) =
h_2(\sum_{g\in U(\Afr_2)/U^{1}(\Afr_2)} g\VV^{\eta_2(V_2,\Afr_2)})
$$
holds true. Since $\Afr_1^{h_1}$ and $\Afr_2^{h_2}$ both correspond to
the barycenter of $\sigma$ they are equal. Hence setting $h :=
h_1^{-1}h_2$ the above identity can equivalently be written as
$$
\sum_{g\in U(\Afr_1)/U^{1}(\Afr_1)} g\VV^{\eta_1(V_1,\Afr_1)} =
\sum_{g\in U(\Afr_1)/U^{1}(\Afr_1)} gh\VV^{\eta_2(V_2,\Afr_2)}\ .
$$
It therefore suffices to find an $x \in U(\Afr_1)$ such that
$$
\eta_1(V_1,\Afr_1)^x \cong \eta_2(V_2,\Afr_2)^h\ .
$$
For this in turn it certainly is sufficient to show that
$$
\theta_1(\Bfr_1)^x = \theta_2(\Bfr_2)^h\ .
$$
Let $V_2^h$ be the $F$-vector space $V_2 = V$ with the new
$E_2$-vector space structure given by $E_2 \hookrightarrow {\rm
End}_{E_2} V_2 \mathop{\longrightarrow}\limits^{h.h^{-1}} {\rm End}_F
V$. By (III.1.3)(a) we have $\Theta_2(V_2,\Bfr_2)^h \in
\CC(V_2^h,\Bfr_2^h)$. Hence there is a unique ps-character
$\Theta_2^h$ supported by $[0,\beta_2]$ such that
$\Theta_2^h(V_2^h,\Bfr_2^h) = \Theta_2(V_2,\Bfr_2)^h$. Obviously
$\Theta_2(V_2,\Bfr_2)$ and $\Theta_2^h(V_2^h,\Bfr_2^h)$ are
simultaneous realizations which intertwine in $G$. Therefore
$\Theta_2$ and $\Theta_2^h$ and hence $\Theta_1$ and $\Theta_2^h$ are
endo-equivalent. Since $\Afr(\Bfr_2^h) = \Afr_2^h = \Afr_1 =
\Afr(\Bfr_1)$ we may apply (II.3.4) to $\Theta_1$ and $\Theta_2^h$ and
obtain an $x \in U(\Afr_1)$ such that
$$
\Theta_1(V_1,\Bfr_1)^x = \Theta_2^h(V_2^h,\Bfr_2^h) =
\Theta_2(V_2,\Bfr_2)^h\ .
$$

\bigskip

\centerline{\bf VIII. On the support of $\CC (\Theta,V,\VV)$.}

\medskip

We fix a simple type $(J,\lambda)$ in the sense of [BK1](5.5.10).
Recall that this means one of the following two possibilities.

(a) There are given a simple pair $[0,\beta]$, an $E$-vector space $V$
where $E := F(\beta)$, and a {\it principal} order $\Bfr_o$ in ${\rm
Her}(B)$ where $B := {\rm End}_{E} V$. The representation $\lambda$ of
the group $J := J(\Bfr_o) := J^1(\Bfr_o)\cdot U(\Bfr_o)$ is of the
form $\kappa\otimes\rho$, where $\kappa$ is a $\beta$-extension to $J$
of a simple character $\theta_o \in \CC(V,\Bfr_o)$ (cf. [BK1](5.2.1))
and $\rho$ is the inflation of an irreducible cuspidal representation
of $J/J^{1}(\Bfr_o)$ of the following kind. Recall that
$J/J^{1}(\Bfr_o )$ identifies with ${\rm GL}_{n/e}(\Fdss_{E})^{\times
e}$, where $n := {\rm dim}_{E} V$, $e$ is the period of the
$\ofr_{E}$-lattice chain corresponding to $\Bfr_o$, and $\Fdss_{E}$
denotes the residue class field of $E$. Then the condition on $\rho$
is that, as a representation of ${\rm GL}_{n/e}(\Fdss_{E})^{\times
e}$, it is of the form $\rho_o^{\otimes e}$, where $\rho_o$ is an
irreducible cuspidal representation of ${\rm GL}_{n/e} (\Fdss_{E})$.
We let $\Theta_o$ denote the unique ps-character supported by
$[0,\beta]$ such that $\Theta_o(V,\Bfr_o) = \theta_o$.

(b) We are in the degenerate case. There is given an $F$-vector space
$V$ and a {\it principal} order $\Afr_o$ in $A := {\rm End}_F V$. The
representation $\lambda$ of the group $J := J(\Afr_o) := U(\Afr_o )$
is the inflation of an irreducible cuspidal representation of
$U(\Afr_o )/U^1 (\Afr_o )$ of the following kind. We have $U(\Afr_o
)/U^1 (\Afr_o ) \cong {\rm GL}_{n/e}(\Fdss_{F})^{\times e}$, where $n
:= {\rm dim}_F V$ and $e$ is the period of the $o_F$-lattice chain
corresponding to $\Afr_o$. Then $\lambda \cong \rho_{o}^{\otimes e}$
for some irreducible cuspidal representation $\rho_o$ of ${\rm
GL}_{n/e}(\Fdss_F )$. In order to have a notation consistent with the
case a), we set $E :=F$, $B := A$, $\Bfr_o := \Afr_o$, $\kappa := {\bf
1}_J$, $\theta_o := {\bf 1}_{J^1(\Afr_o)}$, $\rho := \lambda$, and we
let $\Theta_o$ denote the corresponding degenerate ps-character.

\smallskip

Let $\RR (G)$ denote the category of smooth complex representations
 of $G := {\rm Aut}_F V$. Recall that the full subcategory $\RR_{(J,\lambda)}(G)$
 whose objects are the representations generated by their
 $\lambda$-isotypic component is stable under the formation of
 subquotients. It coincides with a
 Bernstein component of $\RR (G)$ attached to a single point in the
 Bernstein spectrum of $G$ (cf. [BK1] and [BK3](9.3)).

Throughout this section we fix a nonzero representation $\VV$ in
$\RR_{(J ,\lambda)}(G)$.

\smallskip

{\it Remark.} The coefficient systems $\CC_o(\Theta_o,V,\VV)$ and
$\CC(\Theta_o,V,\VV)$ are nonzero.

{\it Proof.} Choose a maximal order $\Bfr$ in ${\rm Her}(B)$
containing $\Bfr_o$ so that $\sigma := (\Bfr \supset \Bfr_o)$ is a
semistandard simplex of $\XE$. It is a consequence of (III.1.1)(iii)
that $\VV(\sigma)$ and $\VV^{\eta(\theta_o)}$ generate the same
$U^1(\Afr(\Bfr_o))$-invariant subspace of $\VV$. But the latter
contains the isotypic component $\VV^{\lambda}$ which is nonzero by
assumption.

\smallskip

{\bf VIII.1 Endo-classes.}

Let $\Theta'$ be an arbitrary ps-character which can be realized in an
$E'$-vector space $V'$ which as an $F$-vector space coincides with
$V$. We assume that $\Theta '$ is either degenerate or supported by a
simple pair $[0,\beta ']$; in particular $E' = F$ or $E' = F(\beta')$.
Let $B' := {\rm End}_{E'}(V')$ and let $\eta (\ldots) =
\eta(V',\ldots)$ denote the various representations attached to
$\Theta'$ and $V'$ as introduced in (III.1).

\smallskip

{\bf (VIII.1.1) Lemma}. {\it Assume that there exists $\Cfr_0\in {\rm
Her}(B')$ such that $\VV$ contains the simple character $\Theta'
(V',\Cfr_0)$. Then there exist a $\Cfr\in {\rm Her}(B')$, a
$\beta'$-extension $\kappa'$ of the Heisenberg representation $\eta
(\Cfr )$ attached to $\Theta' (V',\Cfr )$, and an irreducible cuspidal
representation $\rho'$ of the $\Fdss_{E'}$-reductive  group $U(\Cfr
)/U^1 (\Cfr )$ such that $\VV$ contains the representation
$\kappa'\otimes \rho'$ of $J(\Cfr )$.}

{\it Proof}. (This fact actually is a consequence of the proof of
[BK1](8.1.5), p. 268/269. We shall nevertheless give the argument, for
the context of {\it loc. cit.} is slightly different.)\hfill\break
Take $\Cfr$ minimal among the orders of ${\rm Her}(B')$ such that
$\VV$ contains $\Theta' (V',\Cfr )$. Then $\VV$ must contain the
Heisenberg representation $\eta (\Cfr )$ associated to $\Theta'
(V',\Cfr )$ and a fortiori an irreducible representation $\lambda'$ of
$J(\Cfr )$ such that $\lambda'_{\vert J^1 (\Cfr )}$  contains $\eta
(\Cfr )$. By [BK1](5.2.2) such a representation $\lambda'$ is of the
form $\lambda' =\kappa'\otimes \rho'$, where $\kappa'$ is a
$\beta'$-extension of $\eta (\Cfr )$ and $\rho'$ is (the inflation of)
an irreducible representation of $J(\Cfr )/J^1 (\Cfr ) =
U(\Cfr)/U^1(\Cfr)$. It remains to prove that the minimality condition
on $\Cfr$ implies that $\rho'$ is cuspidal. Assume therefore that
$\rho'$ is not cuspidal. Then there exists a proper parabolic subgroup
$\Pdss$ of $U(\Cfr )/U^1 (\Cfr )$, with unipotent radical $\Udss$,
such that $\rho'_{\vert
\Udss}$ contains the trivial character. There is a uniquely determined
hereditary order $\Cfr_1\subseteq \Cfr$ such that $\Pdss$ (resp.
$\Udss$) is the image of $U(\Cfr_1)$ (resp. $U^1 (\Cfr_1)$) in the
quotient $U(\Cfr )/U^1 (\Cfr )$. Since $\Pdss$ is proper, the
containment $\Cfr_1\subset \Cfr$ is strict. Let $\eta (\Cfr_1)$ be the
Heisenberg representation associated to $\Theta' (V',\Cfr_1)$. We show
that $\VV$ contains $\eta (\Cfr_1)$, hence also $\Theta' (V',\Cfr_1)$,
which contradicts the minimality assumption on $\Cfr$.

The representation $\VV$ contains
$$
(\kappa' \otimes \rho')_{\vert U^1 (\Cfr_1)J^1 (\Cfr )}= \kappa'_{
\vert U^1 (\Cfr_1)J^1 (\Cfr )}\otimes \rho'_{\vert U^1 (\Cfr_1)J^1
(\Cfr )}
$$
which contains
$$
\kappa'_{\vert U^1 (\Cfr_1)J^1 (\Cfr )}\otimes
{\bf 1}_{\vert U^1 (\Cfr_1)J^1 (\Cfr )}
 = \kappa'_{\vert U^1 (\Cfr_1)J^1 (\Cfr )}\ .
$$
Hence our claim follows from [BK1](8.1.6) which says that the
representations of $U^1 (\Afr (\Cfr_1))$ induced by $\eta (\Cfr_1)$
and $\kappa'_{\vert U^1(\Cfr_1)J^1 (\Cfr )}$ are irreducible and
equivalent to each other.

\smallskip

{\bf (VIII.1.2) Proposition}. {\it If the coefficient system $\CC(
\Theta',V',\VV)$ is nonzero then the ps-characters $\Theta'$ and $\Theta_o$
are endo-equivalent and $\CC[\Theta',V',\VV] = \CC[\Theta_o,V,\VV]$.}

{\it Proof}.  The second part of the assertion is a consequence of the
first part by (VII.2). If  $\CC (\Theta',V',\VV)$ is nonzero then
there is a vertex $(\Cfr)$ of ${\rm sd}(X_{E'})$ such that
$\VV^{\eta(\Afr(\Cfr))} \neq 0$. Let $v$ be a nonzero vector in this
isotypic component, let $\VV_0$ be the $G$-subrepresentation of $\VV$
generated by $v$, and let $\VV_1 \subseteq \VV_0$ be the largest
$G$-subrepresentation which does not contain $v$. Then $\VV_0/\VV_1$
is an irreducible $G$-representation in the category
$\RR_{(J,\lambda)}(G)$. By construction we have
$(\VV_0/\VV_1)^{\eta(\Afr(\Cfr))} \neq 0$ and hence $\CC
(\Theta',V',\VV_0/\VV_1) \neq 0$. In order to show that $\Theta'$ and
$\Theta_o$ are endo-equivalent, we may therefore assume in the
following without loss of generality that $\VV$ is irreducible. By
definition if $\CC (\Theta',V',\VV)$ is nonzero then $\CC_o
(\Theta',V',\VV)$ is nonzero, too. In particular there exists a semi-
standard simplex $(\Cfr_0 \supset \dots \supset
\Cfr_q )$ of ${\rm sd}(X_{E'} )$ such that $\VV^{\eta (\Cfr_q ,\Cfr_0
)}\not= 0$.  Since $\eta (\Cfr_q,\Cfr_0 )_{\vert J^{1}(\Cfr_0 )}\simeq
\eta (\Cfr_0 )$, we have $\VV^{\eta (\Cfr_0 )}\not= 0$. On the other
hand the condition that $\VV$ lies in $\RR_{(J,\lambda)}(G)$ implies
that $\VV^{\eta(\theta_o)}\not=0 $. Since $\VV$ is irreducible it
follows that the representations $(J^{1}(\Bfr_o ), \eta(\theta_o) )$
and $(J^{1}(\Cfr_0 ),\eta (\Cfr_0 ))$ intertwine in $G$. Moreover
$\eta(\theta_o)_{\vert H^{1}(\Bfr_o )}$ (resp. $
\eta(\Cfr_0 )_{\vert H^{1}(\Cfr_0 )}$) is a multiple of $\theta_o =
\Theta_o(V,\Bfr_o)$ (resp. of $\Theta' (V',\Cfr_0 )$); so these simple
characters must intertwine as well. Hence, for the endo-equivalence of
$\Theta_o$ and $\Theta'$, it remains to establish the equality $[E:F]
= [E':F]$.

Applying (VIII.1.1) we have that $\VV$ contains a pair $(J(\Cfr),
\kappa' \otimes \rho')$, where $\Cfr\in {\rm Her}(B' )$, $\kappa'$
is a $\beta'$-extension of the Heisenberg representation $\eta(\Cfr)$
attached to $\Theta'(V',\Cfr)$, and $\rho'$ is an irreducible cuspidal
representation of $U (\Cfr)/U^1 (\Cfr)$. Write
$$
U(\Cfr)/U^1 (\Cfr)\simeq \prod_{i=1}^{e'} {\rm GL}_{m_i}(\Fdss_{E'})\
,
$$
where $e' :=e(\Cfr/\ofr_{E'})$ and the $m_i$ are integers $\geq 1$.
Then $\rho'$ writes $\rho_{1}' \otimes \cdots \otimes
\rho_{e'}'$, where, for $i=1,\dots ,e'$,  $\rho_{i}'$ is an
irreducible cuspidal representation of ${\rm GL}_{m_i}(\Fdss_{E'})$.
The pair $(J(\Cfr), \kappa'\otimes \rho')$ is either a simple type or
a split type in the sense of [BK1](8.1), according to whether the
$\rho'_i$ are equivalent to each other or are not. When it is a split
type it has level $(0,0)$ ([BK1](8.1.2)) or level $(-v_{\Afr
(\Cfr)}(\beta')/e(\Afr (\Cfr)/\ofr_{F}),0)$ ([BK1](8.1.4)), according
to whether $\Theta'$ is degenerate or not.

Assume first that $(J(\Cfr), \kappa'\otimes \rho')$ is a simple type.
Since $\VV$ is irreducible, it then is a type for the same Bernstein
component ${\cal R}_{(J,\lambda )}(G)$ of ${\cal R}(G)$. By
[BK1](7.3.17) the pairs $(J,\lambda)$ and $(J(\Cfr),
\kappa'\otimes \rho')$ must be conjugate in $G$, and in particular $e
= e(\Bfr_o/\ofr_{E}) = e(\Cfr/\ofr_{E'}) = e'$ (cf. the proof of
loc.cit.). Setting $n := {\rm dim}_E V$ and $n' := {\rm dim}_{E'} V$
it follows that
$$
GL_{n/e}(\Fdss_E)^{\times e} \cong J/J(\Bfr_o) \cong J(\Cfr)/J^1(\Cfr)
\cong GL_{n'/e'}(\Fdss_{E'})^{\times e'} = GL_{n'/e}(\Fdss_{E'})^{\times
e}
$$
and hence that $n = n'$, i.e., that $[E:F] = [E':F]$.

Assume now that $(J(\Cfr), \kappa'\otimes \rho')$ is a split type. In
this case we need to consider the cuspidal support of the irreducible
representation $\VV$. From the point of view of the simple type
$(J,\lambda)$ in $\VV$ we know from [BK1](7.3.12) that the cuspidal
support of $\VV$ is of the form $(M,\mu)$ where $\mu = \mu_1\otimes
 \ldots \otimes \mu_e$ is a supercuspidal representation of the Levi
subgroup $M={\rm Aut}_{F}(W)^{\times e}$ corresponding to a
decomposition $V=W\oplus \ldots \oplus W$ ($e$ factors) of the
$E$-vector space $V$. Moreover each supercuspidal representation
$\mu_i$ contains ``the maximal type'' $(\widetilde{J},
\widetilde{\lambda})$ attached to $(J,\lambda)$ ([BK1](7.2.18)(iii)).
We do not repeat the definition of $(\widetilde{J},
\widetilde{\lambda})$ but only recall that its underlying simple pair
still is $[0,\beta]$.

On the other hand, from the point of view of the split type $(J(\Cfr),
\kappa'\otimes \rho')$ in $\VV$ we deduce from [BK1](8.3.3) and (6.2.2)
that the cuspidal support of $\VV$ must be of the form $(M',\mu')$
where $\mu' = \mu'_1 \otimes\ldots\otimes \mu'_{e'}$ is a
supercuspidal representation of the Levi subgroup
$M'=\prod_{i=1}^{e'}{\rm Aut}_{F}(W_i)$  corresponding to a
decomposition $V=W_1\oplus\ldots\oplus W_{e'}$ of the $E$-vector space
$V$. Moreover each supercuspidal representation $\mu'_i$ contains a
simple type $(J_i,\lambda_i)$ with underlying simple pair
$[0,\beta']$. By unicity of the cuspidal support, the pairs $(M,\mu )$
and $(M',\mu')$ are conjugate in $G$. So after conjugation, we may
reduce to the case where, e.g., the representation $\mu_1 \cong
\mu'_1$ contains two simple types with underlying simple pairs $[0,\beta]$
and $[0,\beta']$, respectively, and we may conclude as in the first
case.

\medskip

Let ${\rm Coeff}_G(\XX)$ denote the category of $G$-equivariant
coefficient systems on the simplicial complex $\XX$. The Proposition
(VIII.1.1) together with the above Remark imply that, given a simple
type $(J,\lambda)$, the functor
$$
\matrix{
  \CC_{(J,\lambda)} : \RR_{(J,\lambda)}(G) & \longrightarrow & {\rm
  Coeff}_G(\XX) \cr\cr
  \hfill \VV & \longmapsto & \CC(\Theta_o,V,\VV) \hfill
  }
$$
is independent of any additional choices. In order to be able later on
to show that this functor in fact is a fully faithful embedding we
first have to analyze the support of these coefficient systems.

\smallskip

{\bf VIII.2 The support of $\CC [\Theta_o,V,\VV]$.}

As in (IV) and (V) we write $\CC_o(\Theta_o,V,\VV) = (\VV(\s))_\s$ and
$\CC[\Theta_o,V,\VV] = (\VV[\s])_\s$. To start with we fix a simplex
$\sigma_0 = (\Bfr_{max} \supset \dots
\supset \Bfr_{min})$ of maximal dimension in $\XE$ such that
 $\Bfr_{min} \subseteq \Bfr_o \subseteq \Bfr_{max}$. Recall
 ([BK0](5.2.2-5)) that the $\beta$-extension $\kappa$
 gives rise to a compatible family of $\beta$-extensions $\kappa
 (\Bfr)$ where $(\Bfr)$ runs over the vertices of $\sigma_0$.
These $\kappa (\Bfr )$ are characterized as follows:

(a) The induced representations
$$
{\rm Ind}_{J(\Bfr )}^{U(\Bfr )U^{1}(\Afr (\Bfr ))} \kappa (\Bfr )
\quad\hbox{and}\quad {\rm Ind}_{U(\Bfr )J^1(\Bfr_{max} )}^{U(\Bfr )U^{1}(\Afr
(\Bfr ))} \kappa (\Bfr_{max} )
$$
are isomorphic (and irreducible);

(b) $\kappa (\Bfr_o ) = \kappa$.

\smallskip

 Set $\Gdss = U(\Bfr_{max} )/U^{1}(\Bfr_{max} )\simeq {\rm GL}_n(\Fdss_{E})$.
 Following [SZ]{\S}5, we define the $\Gdss$-module $\VV (\Bfr_{max}) := {\rm
 Hom}_{J^{1}(\Bfr_{max} )}(\kappa (\Bfr_{max} ),\VV )$, using the obvious action
 of $J(\Bfr_{max} )$  and the canonical identification
 $$
 J(\Bfr_{max} )/J^1 (\Bfr_{max}
 )\simeq U(\Bfr_{max} )/U^1 (\Bfr_{max} )\ .
 $$
 Recall (loc. cit.) that for
 $\Bfr_{min} \subseteq \Bfr \subseteq \Bfr_{max}$ the image of $U(\Bfr )J^1 (\Bfr_{max}
 )$ in  $J(\Bfr_{max} )/J^1 (\Bfr_{max} )$ is a parabolic subgroup
 $\Pdss_{\Bfr}$ of $\Gdss$ whose unipotent radical $\Udss_{\Bfr}$
 is the image of $U^{1}(\Bfr )J^{1}(\Bfr_{max} )$ and
 whose Levi quotient $\Ldss_{\Bfr }$  canonically identifies with
 $U(\Bfr )/U^1 (\Bfr )$.

\smallskip

{\bf (VIII.2.1) Proposition}. {\it For any $\Bfr_{min} \subseteq \Bfr
\subseteq \Bfr_{max}$ we have linear isomorphisms
$$
\VV (\Bfr_{max} )^{\Udss_{\Bfr}}\simeq {\rm Hom}_{J^{1}(\Bfr ,\Bfr_{max} )}(\eta
(\Bfr ,\Bfr_{max} ),\VV  )
\simeq  {\rm Hom}_{J^{1}(\Bfr )}(\eta (\Bfr ), \VV)\ ,
$$
where $\VV (\Bfr_{max} )^{\Udss_{\Bfr}}$ denotes the Jacquet module
with respect to the parabolic subgroup $\Pdss_{\Bfr}$.}

{\it Proof}. According to the proof of Lemma 2 in  [SZ]{\S}5, we have:
$$
\VV (\Bfr_{max} )^{\Udss_{\Bfr}}\simeq {\rm Hom}_{J^{1}(\Bfr ,\Bfr_{max}
)}(\kappa (\Bfr_{max} ),\VV  )
\simeq  {\rm Hom}_{J^{1}(\Bfr )}(\kappa (\Bfr ), \VV)\ .
$$
So we must simply  prove the isomorphisms:
$$
\kappa (\Bfr_{max} )_{\vert J^{1}(\Bfr ,\Bfr_{max} )} \simeq \eta(\Bfr ,\Bfr_{max}
)\ ,\  \kappa (\Bfr )_{\vert J^{1}(\Bfr )} \simeq \eta (\Bfr)\ .
$$
The second one is clear by definition of a $\beta$-extension. Write
$$
\eta_{max} := \kappa (\Bfr_{max} )_{\vert J^1 (\Bfr ,\Bfr_{max} )}\ .
$$
Using Mackey's restriction formula, the restrictions to $U^1 (\Afr
(\Bfr ))$ of the isomorphic representations ${\rm Ind}_{J(\Bfr
)}^{U(\Bfr )U^1 (\Afr (\Bfr ))}\kappa (\Bfr )$ and ${\rm Ind}_{U(\Bfr)
J^1(\Bfr_{max} )}^{U(\Bfr )U^1 (\Afr (\Bfr ))}\kappa (\Bfr_{max}  )$
are
$$
{\rm Ind}_{J^1 (\Bfr )}^{U^1 (\Afr (\Bfr ))}\eta (\Bfr )\simeq {\rm
Ind}_{J^{1}(\Bfr ,\Bfr_{max} )}^{U^1 (\Afr (\Bfr ))}\eta_{max} \ .
$$
Moreover by definition of a $\beta$-extension $\eta_{max\,\vert J^1
(\Bfr_{max} )} \simeq \eta (\Bfr_{max} )$.  So by definition of $\eta
(\Bfr ,\Bfr_{max} ) $ (cf. (III.1.1) and [BK1](5.1.16)), we have
$\eta_{max} \simeq \eta (\Bfr ,\Bfr_{max} ) $, as required.

\smallskip

 We shall also need the following fact from.

\smallskip

{\bf (VIII.2.2) Proposition} ([SZ]{\S}5 Prop. 3). {\it Any irreducible
constituent of the $\Gdss$-module $\VV (\Bfr_{max} )$ has cuspidal
support $(\Ldss_{\Bfr_o} ,\rho)$.}

\smallskip

As a corollary of the last two propositions we obtain the following
result.

\smallskip

{\bf (VIII.2.3) Proposition}. {\it Let $\sigma =(\Bfr_{max} \supset
\dots \supset \Bfr_q )$ be a semistandard simplex contained in $\sigma_0$.
Then $\VV(\sigma )=\VV^{\eta (\Bfr_q ,\Bfr_{max})}\not= 0$ if and only
if $\Ldss_{\Bfr_q}$ contains a Levi subgroup conjugate in $\Gdss$ to
$\Ldss_{\Bfr_o}$. In other words, if the invariant of the conjugacy
class of $\Ldss_{\Bfr_q}$ is the unordered partition $(n_1 ,\dots ,n_s
)$ of $n$, we have $\VV (\sigma )\not= 0$ if and only if $n/e$ divides
$n_i$ for any $i=1,\dots ,s$.}

\smallskip

As in (I.3.3) we introduce, for any $\Bfr \in {\rm Her}(B)$, the
integers
$$
d_E(\Bfr)_k := {\rm dim}_{\Fdss_E} L_k/L_{k+1}
$$
where $(L_k)_{k \in\Zdss}$ is an $o_E$-lattice chain in $V$
corresponding to $\Bfr$. The condition of the proposition can be read
off the sequence $[d_{E}(\Bfr_q )_k ]_{k\in \Zdss}$ and, of course,
only depends on the $G_E$-conjugacy class of $\Bfr_q$.

\smallskip

{\bf (VIII.2.4) Corollary}. {\it Let $\sigma = (\Bfr_0 \supset
\dots \supset \Bfr_q )$ be any simplex of ${\rm sd}(X_E )$.
Then $\VV (\sigma )\not= 0$ if and only if $n/e$ divides $d_{E}(\Bfr_q
)_k$ for any $k\in \Zdss$.}

{\it Proof}. By definition $\VV (\sigma )\not= 0$ if and only if their
exists a semistandard simplex $\tau$ containing $\sigma$ such that
$\VV (\tau )\not= 0$. So $\VV (\sigma )\not= 0$ if and only if their
exists $\Bfr\in {\rm Her(B)}$ such that $\Bfr \subseteq \Bfr_q$ and
$(n/e)\vert d_{E}(\Bfr )_k$ for any $k\in \Zdss$. But this implies
that $(n/e)\vert d_{E}(\Bfr_q )_k$ for any $k\in \Zdss$, since any
$\ofr_{E}$-lattice occurring in a lattice chain defining $\Bfr_q$
occurs in any lattice chain defining $\Bfr$.

\smallskip

{\bf (VIII.2.5) Corollary}.  {\it Let $\sigma$ be a simplex of $X[E]$;
write $\sigma$ as the image under some $h\in G$ of a simplex of
$X_{E}$ with barycenter $(\Bfr)$, $\Bfr\in {\rm Her}(B)$. Then
$\VV[\sigma]\not= 0$ if and only if $(n/e)\vert d_{E}(\Bfr )_k$ for
any $k\in \Zdss$.}

{\it Proof}. By $G$-equivariance, we may assume that $h=1$. We have
$$
\VV[\sigma] =\sum_{g\in U(\Afr )/U^1 (\Afr )}g\VV^{\eta (\Afr )}\ ,
$$
where $\Afr =\Afr (\Bfr )$. So if $\Bfr_{max}$ is any maximal order in
${\rm Her}(B)$ containing $\Bfr$, by (V.5), we have
$$
\VV[\sigma] =\sum_{g\in U(\Afr )}g \VV ((\Bfr_{max} \supset \Bfr))
$$
and the result follows easily.

\smallskip

We are now going to describe the support of $\CC [\Theta_o,V,\VV]$ in
terms of an auxiliary building. Thanks to (I.3.5), we find an
unramified extension $L$ of $E$ contained in $B$, of degree $n/e$, and
such that $L^{\times}$ normalizes $\Bfr_o$. Write $C:={\rm
End}_{L}V\simeq {\rm M }(e,L)$ for the centralizer of $L$ in $B$. From
(I.2.1), we have a canonical $G_L$-equivariant embedding $j_L$ of the
semisimple affine building $X_L$ of $G_L$ into $X_{E}$. Since $L/E$ is
unramified, this embedding is actually simplicial; indeed in that case
if $\Cfr \in {\rm Her}(C)$ is maximal (i.e. corresponds to a vertex in
$X_L$) then the corresponding order $\Bfr (\Cfr )\in {\rm Her}(B)$ is
maximal as well.  We see $X_L$ as a simplicial subcomplex of $X_{E}$
and ${\rm sd}(X_L)$ as a simplicial subcomplex of $X_E$. So as in (V),
we may consider the simplicial complex $X[L]$; this is a $G$-invariant
simplicial subcomplex of $X[E]$.

\smallskip

{\bf (VIII.2.6) Proposition}. {\it For any simplex $\sigma$ of $X[E]$,
we have $\VV[\sigma]\not= 0$ if and only if $\sigma$ lies in $X[L]$.}

{\it Proof}. By $G$-equivariance we may assume that $\sigma$ actually
lies in $\XE$. We then must prove that $\VV[\sigma]\not= 0$ if and
only if $\sigma \in G_E (X_L)$. By the criterion of (I.3.5) (applied
to ``$E/F$''$=$``$L/E$''), this latter condition is equivalent to
$f(L/E)\vert d_{E}(\Bfr )_{k}$ for any $k\in \Zdss$, where $\Bfr$ is
the barycenter of $\sigma$. But $f(L/E)=n/e$. So we are done using
(VIII.2.5).

\smallskip

We therefore may and will view the functor $\CC_{(J,\lambda)}$
introduced at the end of (VIII.1) as a functor
$$
  \CC_{(J,\lambda)} : \RR_{(J,\lambda)}(G) \longrightarrow {\rm
  Coeff}_G(X[L])\ .
$$

\bigskip

\centerline {\bf  IX. The chain complex attached to $\CC_{(J,\lambda)}(\VV)$}

\medskip

As in the previous section we fix a simple type $(J,\lambda )$ in $G =
{\rm Aut}_F V$ where $\lambda = \kappa\otimes\rho$, with ps-character
$\Theta_o$ having a realization in $V$, and a smooth complex
representation $\VV$ in ${\cal R}_{(J,\lambda )}(G)$. We also keep
most of the other notations introduced in (VIII). We consider the
$G$-equivariant coefficient system $\CC := \CC[\Theta_o,V,\VV] =
(\VV[\sigma])_\sigma$ that we see as a coefficient system on the
$G$-invariant simplicial subcomplex $X[L]$ of $X[E]$. This complex has
dimension
$$
d:={\rm dim}\,X[L] = {\rm dim}_F V/[L:F] - 1 = e - 1
$$
where $e$ is the divisor of ${\rm dim}_E V $ fixed in (VIII). We
denote by $X[L]_q$ the set of $q$-simplices of $X[L]$ for
$q=0,\dots,d$. The following considerations are copied from [SS].

An ordered $q$-simplex in $X[L]$ is a sequence $(\sigma_0 ,\dots,
\sigma_q )$ of vertices such that $\{ \sigma_0 ,\dots ,\sigma_q\}$
is a $q$-simplex. Two such ordered simplices are called equivalent if
they differ by an even permutation of the vertices; the corresponding
equivalence classes are called oriented $q$-simplices and are denoted
by $\langle \sigma_0 ,\dots ,\sigma_q \rangle$. We let $X[L]_{(q)}$ be
the set of oriented $q$-simplices of $X[L]$. The space $C_{q}^{\rm
or}(X[L],\CC )$ of oriented $q$-chains of $X[L]$ with coefficients in
$\CC$ is the $\Cdss$-vector space of all maps
$$
\omega~: X[L]_{(q)} \longrightarrow \VV
$$
such that:

-- $\omega$ has finite support,

-- $\omega (\langle \sigma_0 ,\dots ,\sigma_q \rangle )\in
   \VV[\{\sigma_0 ,\dots ,\sigma_q\}]$,

--  $\omega (\langle \sigma_{\iota (0)} ,\dots ,\sigma_{\iota (q)}
    \rangle ) = {\rm sgn}(\iota ) \cdot \omega (\langle \sigma_0 ,\dots
    ,\sigma_q \rangle )$ for any permutation $\iota$.

The group $G$ acts smoothly on  $C_{q}^{\rm or}(X[L],\CC )$ via
$$
(g\omega )(\langle \sigma_0 ,\dots ,\sigma_q \rangle ):= g(\omega
(\langle g^{-1} \sigma_0 ,\dots , g^{-1}\sigma_q \rangle ))\ .
$$
With respect to the $G$-equivariant boundary maps
$$
\matrix{
\partial~:\ C_{q+1}^{\rm or}(X[L],\CC ) & \rightarrow &
            C_{q}^{\rm or}(X[L],\CC) \hfill\cr\cr
   \hfill \omega & \mapsto & [\langle \sigma_0 ,\dots ,\sigma_q \rangle\mapsto
\mathop{\sum}\limits_{\{\sigma ,\sigma_0 ,\dots \sigma_q\}\in X[L]_{q+1}}\omega
(\langle \sigma ,\sigma_0 ,\dots ,\sigma_q \rangle )] }
$$
we then have the augmented chain complex in $\RR (G)$:
$$
C_d^{\rm or}(X[L],\CC )\boundary \cdots \boundary C_0^{\rm
or}(X[L],\CC )\augmentation \VV \leqno (\hbox{IX.1})
$$
where $\displaystyle \epsilon (\omega )=\sum_{\sigma\in X[L]_{(0)}}
\omega (\sigma )\in \VV$.

\smallskip

{\bf (IX.2) Proposition}. {\it For all $q=0,\dots ,d$, the $G$-module
$C_{q}^{\rm or}(X[L],\CC )$ lies in $\RR_{(J,\lambda )}(G)$. In
particular the complex (IX.1) is a chain complex in the category
$\RR_{(J,\lambda )}(G)$.}

\smallskip

To prepare for the proof we let $\sigma_{\Cfr}$, for any  $\Cfr\in
{\rm Her}(C)$, denote the simplex of $X_L$ with barycenter $(\Cfr )$.
Moreover let ${<}\sigma_{\Cfr}{>}$ be a fixed oriented simplex with
underlying simplex $\sigma_{\Cfr}$ and let
$\overline{{<}\sigma_{\Cfr}{>}}$ denote that oriented simplex with the
same underlying simplex $\sigma_{\Cfr}$ but with the reversed
orientation (for vertices we have
${<}\sigma_{\Cfr}{>}=\overline{{<}\sigma_{\Cfr}{>}}$). The order
$\Bfr_o$ defining $J=J(\Bfr_o)$ corresponds to a minimal order
$\Cfr_{min}$ of $C$. We write $\Bfr_{min} := \Bfr_o = \Bfr (\Cfr_{min}
)$ and put $\Afr_{min} :=\Afr (\Bfr_{min} )$. We fix a maximal order
$\Cfr_{max}$ of $C$ containing $\Cfr_{min}$ and put $\Bfr_{max} :=
\Bfr (\Cfr_{max} )$ and $\Afr_{max} :=
\Afr(\Bfr_{max})$. We have $\Afr_{min} \subset \Afr_{max}$ and
$\Bfr_{min} \subset \Bfr_{max}$ and since $L/E$ is unramified,
$\Bfr_{max}$ is a maximal order of $B$.  Note that $\Bfr_{min}$ in
general is not a minimal hereditary order of $B$.

Any simplex in $X[L]$ lies in the $G$-orbit of a simplex $\sigma_\Cfr$
with $\Cfr_{min} \subseteq \Cfr \subseteq \Cfr_{max}$. Hence
$$
C_{q}^{\rm or}(X[L] ,\CC ) =
\sum_{\Cfr_{min}\subseteq\Cfr\subseteq\Cfr_{max} ,{\rm
dim}\sigma_{\Cfr}=q} C_q^{\rm or}(\sigma_{\Cfr},\CC)
$$
where
$$
C_q^{\rm or}(\sigma_{\Cfr},\CC):=\{\omega\in C_q^{\rm or}(X[L],\CC):
\omega\ \hbox{has support in}\ G{<}\sigma_{\Cfr}{>}\cup\
G\overline{{<}\sigma_{\Cfr}{>}}\}
$$
and we are reduced to showing that the $G$-subrepresentations
$C_q^{\rm or}(\sigma_{\Cfr},\CC)$, for
$\Cfr_{min}\subseteq\Cfr\subseteq\Cfr_{max}$, are generated by their
$\lambda$-isotypic components. In the following we fix such an order
$\Cfr_{min}\subseteq\Cfr\subseteq\Cfr_{max}$ and put $\Bfr :=
\Bfr(\Cfr)$ and $\Afr := \Afr(\Bfr)$. We may embed
$\VV[\sigma_{\Cfr}]$ in a $U(\Afr)$-equivariant way into $C_q^{\rm
or}(\sigma_{\Cfr},\CC)$ by
$$
\matrix{
\VV[\sigma_{\Cfr}] & \longrightarrow & C_q^{\rm
or}(\sigma_{\Cfr},\CC)\cr \hfill v & \longmapsto & \omega_v \hfill }
$$
where
$$
\omega_v(.):=\left\{\matrix{ +v\hfill &
\hbox{if}\ \ .={<}\sigma_{\Cfr}{>}\ \ ,\hfill\cr
-v\hfill & \hbox{if}\ q>0\ \hbox{and}\ \
.=\overline{{<}\sigma_{\Cfr}{>}}\ \ ,\hfill\cr 0\hfill &
\hbox{otherwise}\ \ .\hfill\cr}\right.
$$
In the following we view this embedding as an inclusion. Clearly
$\VV[\sigma_{\Cfr}]$ generates $C_q^{\rm or}(\sigma_{\Cfr},\CC)$ as a
$G$-representation. According to (V.10) and (V.5), we have
$$
\VV[\sigma_{\Cfr}] = \VV_{(\Bfr_{max} \supset \Bfr)} =
\sum_{g\in U(\Afr )/U(\Bfr )J^{1}(\Bfr_{max}
)}g\VV^{\eta (\Bfr ,\Bfr_{max} )} \ .
$$
and hence
$$
C_q^{\rm or}(\sigma_{\Cfr},\CC) =
\sum_{g\in G/U(\Bfr)J^{1}(\Bfr_{max})}
    g\VV^{\eta (\Bfr ,\Bfr_{max} )}\ . \leqno (\hbox{IX.3})
$$
Having fixed a compatible family of $\beta$-extensions $\kappa(.)$ as
in (VIII.2) we in particular have $\kappa(\Bfr_{max})$ as a
representation of $J(\Bfr_{max}) = J^1(\Bfr_{max})\cdot
U(\Bfr_{max})$. We then may form the representation $\lambda_{max}
=\kappa (\Bfr_{max} )\otimes \rho$ of $U(\Bfr_{min} )J^{1}(\Bfr_{max}
)$. Both factors in this tensor product are irreducible, the second
factor by assumption and the first factor since it restricts to the
irreducible representation $\eta(\Bfr_{max})$ on $J^1(\Bfr_{max})$.
Therefore, by the argument in the proof of [BK1](5.3.2), the
representation $\lambda_{max}$ is irreducible.

\smallskip

{\bf (IX.4) Lemma}. {\it i) A smooth $G$-representation lies in
$\RR_{(J,\lambda)}(G)$ if and only if it is generated by its
$\lambda_{max}$-isotypic component.

ii) We have $\VV^{\eta (\Bfr_{min} ,\Bfr_{max} )}
=\VV^{\lambda_{max}}$.}

{\it Proof}. i) According to the proof of [BK1](5.5.13) we have the
isomorphism
$$
{\rm Ind}_{U(\Bfr_{min} )J^{1}(\Bfr_{max} )}^{U(\Bfr_{min}
)U^{1}(\Afr_{min} )}\lambda_{max} \simeq {\rm Ind}_{J(\Bfr_{min}
)}^{U(\Bfr_{min} )U^{1}(\Afr_{min} )}\lambda
$$
So by Frobenius reciprocity, the $U(\Bfr_{min} )U^{1}(\Afr_{min}
)$-submodules  in a smooth $G$-representation $\cal W$ generated by
${\cal W}^{\lambda_{max}}$ and ${\cal W}^{\lambda}$, respectively,
coincide.

ii) According to the proof of (VIII.2.1) we have
$$
\kappa (\Bfr_{max} )_{\vert J^{1}(\Bfr_{min},\Bfr_{max} )} \simeq
\eta(\Bfr_{min},\Bfr_{max})\ .
$$
Hence $\lambda_{max\, \vert J^{1}(\Bfr_{min},\Bfr_{max} )}$ is
$\eta(\Bfr_{min},\Bfr_{max})$-isotypic which shows that
$\VV^{\lambda_{max}} \subseteq \VV^{\eta (\Bfr_{min} ,\Bfr_{max} )}$.
But it also implies that $\VV^{\eta (\Bfr_{min} ,\Bfr_{max} )}$ is the
image of
$$
\kappa(\Bfr_{max}) \otimes {\rm Hom}_{J^1 (\Bfr_{min},\Bfr_{max} )} (\kappa
(\Bfr_{max} ),\VV )
$$
under the canonical map into $\VV$. For the reverse inclusion
$\VV^{\eta (\Bfr_{min} ,\Bfr_{max} )} \subseteq \VV^{\lambda_{max}}$
it therefore suffices to prove that ${\rm Hom}_{J^1
(\Bfr_{min},\Bfr_{max} )} (\kappa (\Bfr_{max} ),\VV )$ as a
$U(\Bfr_{min} )/U^1 (\Bfr_{min} )$-module is $\rho$-isotypic. But by
the first formula in the proof of (VIII.2.1) this latter module is the
Jacquet module $\VV(\Bfr_{max})^{\Udss_{\Bfr_{min}}}$ (notation of
(VIII)) of the $U(\Bfr_{max} )/U^1 (\Bfr_{max} )$-module ${\rm
Hom}_{J^1 (\Bfr_{max} )}(\kappa (\Bfr_{max} ), \VV )$. From (VIII.2.2)
we know that the latter representation has cuspidal support
$(\Ldss_{\Bfr_{min}},\rho)$. Since our $\rho$ is of the form $\rho
\simeq \rho_o^{\otimes e}$ it follows that the Jacquet module
 $\VV(\Bfr_{max})^{\Udss_{\Bfr_{min}}}$ indeed is $\rho$-isotypic.

\smallskip

In order to prove that the $G$-representation in (IX.3) is generated
by its $\lambda$-isotypic component, it suffices to prove that it is
generated by its $\lambda_{max}$-isotypic component. Since the right
hand version of this representation visibly is generated by $\VV^{\eta
(\Bfr ,\Bfr_{max} )}$ and since $\VV^{\eta (\Bfr_{min} ,\Bfr_{max} )}
\subseteq \VV^{\eta (\Bfr ,\Bfr_{max} )}$ by (III.1.1)(iv), we are finally reduced to
establishing the following fact.

\smallskip

{\bf (IX.5) Lemma}. {\it $\VV^{\eta (\Bfr ,\Bfr_{max} )}$, as a
$U(\Bfr )J^1 (\Bfr_{max} )$-module, is generated by\break $\VV^{\eta
(\Bfr_{min} ,\Bfr_{max} )}$.}

{\it Proof}. We first of all note that $\VV^{\eta (\Bfr ,\Bfr_{max}
)}$, by (III.1.1)(v) and (V.1), indeed is $U(\Bfr )J^1 (\Bfr_{max}
)$-invariant. In the proof of (IX.4)(ii) we have seen that $\VV^{\eta
(\Bfr_{min} ,\Bfr_{max} )}$ is the image of
$$
\kappa(\Bfr_{max}) \otimes {\rm Hom}_{J^1 (\Bfr_{min},\Bfr_{max} )} (\kappa
(\Bfr_{max} ),\VV )
$$
under the canonical map into $\VV$. Analogously $\VV^{\eta (\Bfr,
\Bfr_{max} )}$ is the image of
$$
\kappa(\Bfr_{max}) \otimes {\rm Hom}_{J^1 (\Bfr,\Bfr_{max} )} (\kappa
(\Bfr_{max} ),\VV )\ ,
$$
and this in fact in a $U(\Bfr )J^1 (\Bfr_{max} )$-equivariant way
since $U(\Bfr )J^1 (\Bfr_{max} )$ normalizes $J^1(\Bfr,\Bfr_{max})$.
We therefore are reduced to proving that
$$
{\rm Hom}_{J^1(\Bfr_{min},\Bfr_{max} )} (\kappa (\Bfr_{max} ),\VV )
$$
generates
$$
{\rm Hom}_{J^1(\Bfr,\Bfr_{max} )} (\kappa (\Bfr_{max} ),\VV )
$$
as a $U(\Bfr )J^1 (\Bfr_{max} )$-module. But in that proof we also
have seen (with the notations of (VIII)) that the former is the
Jacquet module $\VV(\Bfr_{max})^{\Udss_{\Bfr_{min}}}$ and the latter
is the Jacquet module $\VV(\Bfr_{max})^{\Udss_{\Bfr}}$ of the
$\Gdss$-module $\VV(\Bfr_{max}) = {\rm Hom}_{J^1 (\Bfr_{max} )}(\kappa
(\Bfr_{max} ),\VV )$. Hence we are further reduced to showing that the
module $\VV(\Bfr_{max})^{\Udss_{\Bfr}}$ for the Levi group
$\Ldss_{\Bfr}$ is generated by its Jacquet module
$\VV(\Bfr_{max})^{\Udss_{\Bfr_{min}}}$. For this it suffices that all
irreducible constituents of the $\Ldss_{\Bfr}$-module
$\VV(\Bfr_{max})^{\Udss_{\Bfr}}$ have cuspidal support on
$\Ldss_{\Bfr_{min}}$. Because of the special form of the group
$\Ldss_{\Bfr_{min}}$ this follows from the fact (cf. (VIII.2.2)) that
any irreducible constituent of the $\Gdss$-module $\VV(\Bfr_{max})$
has cuspidal support on $\Ldss_{\Bfr_{min}}$.

\smallskip

This finishes the proof of Proposition (IX.2).

\bigskip

\centerline{\bf  X. Acyclicity of the chain complex: a strategy}

\bigskip

In this section we consider the augmented complex (IX.1). We reduce its exactness to a technical hypothesis
(conjecture (X.4.1)) that we cannot prove. In the next section we shall prove this hypothesis
for irreducible discrete series representation.

\medskip

{\bf X.1 Some lemmas on $\lambda_{\rm max}$-isotypic components.}
\medskip

As in {\S}IX we fix a simple type $(J, \lambda )$ in $G$ and a smooth
complex representation $\VV$ in $\RR_{(J,\lambda )}(G)$. We keep the
same notation. We abbreviate $J_{\rm max}=U(\Bfr_{\rm
  min})J^{1}(\Bfr_{\rm max})$ and write $\Lambda$ for the representation
space of $\lm$.

 We fix a Haar measure $\mu$ on $G$ and let ${\cal H}(G)$ denote the
(convolution) Hecke algebra of locally constant
functions with compact support on $G$. For $\varphi\in \HH (G)$ and $g\in G$,
we also define ${}^g  \varphi\in \HH (G)$ by ${}^g\varphi (x)=\varphi (g^{-1}x)$.
We also recall the Schur orthogonality formula: if $(\rho,\WW)$ is an irreducible
representation of a compact subgroup $K$ of $G$, with contragredient representation
$(\check\rho,\check\WW)$, then
$$
\int_K \langle\rho(x^{-1})w,\check w\rangle\langle \rho(x)v,\check v\rangle dk =
{{\rm dim}(\rho) \over {\mu (K)}}
\langle w,\check v\rangle \langle v,\check w\rangle,
 \qquad v,w\in\WW,\ \check v,\check w\in\check\WW,
$$
where
$\langle -,-\rangle$~: $\WW\times {\check\WW}\longrightarrow \Cdss$
denotes the canonical pairing.

The irreducible representation
$\lm$ gives rise to an idempotent $e_{\rm max}$ of ${\cal H}(G)$ defined
  as follows: it has support $\Jm$ and is given by
$$
e_{\rm max}(j) = \mu (\Jm )^{-1}{\rm dim}(\lm ){\rm Tr}(\lm (j^{-1}))
$$
for $j\in \Jm$ (cf. [BK1] {\S}4.2)
Note that $e_{\rm max}$ may be considered as an idempotent of the
Hecke algebra ${\cal H}(\Jm):=\{ f\in {\cal H}(G)\ ; \ {\rm
  Support}(f)\subset \Jm\}$. If $(\zeta ,\UU )$ is a smooth
representation of $G$ (resp. of $\Jm$) then $(\zeta ,\UU )$ extends to a
representation of ${\cal H}(G)$ (resp. ${\cal H}(\Jm)$) on $\UU$, and  we
then have  $\zeta (e_{\rm max})\star \UU = \UU^{\lm}$ (the
$\lm$-isotypic component of $\UU$).

 For $x\in G$, we denote by $\lm^{x}$ the representation of
 $\Jm^{x}:=x\Jm x^{-1}$ in the space $\Lambda$ given by
 $\lm^{x}(xjx^{-1})=\lm (j)$, $j\in \Jm$.
\medskip

{\bf (X.1.1) Proposition.} {\it i) Any non-zero function $f$ in the
  scalar Hecke algebra $\em\star \HH (G)\star \em$ has support in the
  $G$-intertwining $I_{G}(\lm )$ of $\lm$.

ii) Let $x$ be an element of $G$ such that $x\not\in I_G (\lm )$ and let $(\pi,\VV)$ be a smooth representation of $G$.
 Then the linear map $p_x$~:
$\VV^{\lm}\lra \VV^{\lm}$, given by $p_x (v)=\pi (\em ) \circ \pi
 (x) \circ \pi (\em ) .v$  is zero.}
\medskip

{\it Remark}. These facts are certainly well known but we could not find a
reference.

\medskip

{\it Proof}. i) Let $\HH (G,\lm )$ be the Hecke algebra of
$\lm$-spherical functions on $G$ ([BK](4.1)). Recall that if
$(\check\lm,\check\Lambda)$ denotes the contragredient
representation of $(\lm ,\Lambda)$, then $\HH (G,\lm )$ is the
convolution algebra of compactly supported functions $\Phi$~: $G\lra
{\rm End}_{\Cdss}(\check\Lambda)$ satisfying~:
$$
\Phi (j_1 gj_2 )=\check\lm(j_1 )\circ \Phi (g) \circ \check\lm (j_2
)\ ,  \ j_i\in \Jm\ , \ g\in G\ .
$$

>From [BK](4.1.1), any non-zero $\Phi\in \HH (G,\lm )$ has support in $I_{G}(\lm
)$. Moreover by [BK], proposition (4.2.4), we have an algebra
isomorphism
$${\bf \Upsilon}~:\ \HH (G,\lm )\otimes_{\Cdss} {\rm End}_{\Cdss} (\Lambda
)\lra \em \star \HH (G)\star \em\ .
$$
Identifying ${\rm End}_{\Cdss} (\Lambda)$ with
$\Lambda\otimes_{\Cdss}\check\Lambda$,
$\bf \Upsilon$ is given by
$$
{\bf \Upsilon} (\Phi\otimes w\otimes \check w)(g)={\rm dim}(\lm ){\rm Tr}
(w\otimes \Phi (g)\check w)
$$
for $g\in G$, $\ w \in \Lambda,\ \check w\in
\check\Lambda, \ \Phi\in \HH (G,\lm ).$
In particular we have:
$$
{\rm Support}\big({\bf \Upsilon} (\Phi \otimes w\otimes \check w)\big)\subset {\rm
Support}(\Phi )\ , \ w \in \Lambda,\ \check w\in
\check\Lambda, \ \Phi\in \HH (G,\lm )\ .
$$
It follows that any non-zero element of $ \em \star \HH (G)\star \em$ has
support in $I_{G}(\lm )$ as required.

ii) Recall that, for $\varphi\in \HH (G)$ and $g\in G$, we write
${}^g  \varphi\in \HH (G)$ for the function ${}^g\varphi (x)=\varphi (g^{-1}x)$. Then
 straightforward computations show that $\pi (\em )\circ  \pi (x)\circ \pi
 (\em )=\pi (\em\star {}^x\em )$ and that $ \em\star
 {}^x\em \in \em \star \HH (G)\star\em$. Moreover $  \em\star
 {}^x\em$ clearly has support in $\Jm x\Jm$, whence is zero since
 $x\not\in I_G (\lm )$.  So $p_x
 =\pi (\em\star {}^x\em )$ is the zero map.

\bigskip

{\bf (X.1.2) Proposition}. {\it Let $x$ be a fixed element of $I_G
(\lm )$.

i) There exist $m\geq 1$, $u_1 ,...,u_m$, $v_1
,...,v_m\in \Jm$, $\gamma_1 ,...,\gamma_m\in \Cdss$, such that
$$
\sum_{i=1}^{m}\gamma_i \, \em\star {}^{u_i xv_i}\em
$$
is an invertible element of $\em\star\HH (G)\star\em$.

ii) There exist $m\geq 1$, $u_1 ,...,u_m$, $v_1
,...,v_m\in \Jm$, $\gamma_1 ,...,\gamma_m\in \Cdss$, such that,
for $(\pi ,\VV )$ any smooth representation of $G$
$$
\sum_{i=1}^{m}\gamma_i\,  \pi (\em )\circ \pi (u_i xv_i )\circ \pi
(\em )
$$
induces a $\Cdss$-linear isomorphism $\VV^\lm\lra \VV^\lm$.}

{\it Proof}. To make the notation lighter, we shall set $K=J_{\rm
max}$, $\rho =\lm$, $e=\em$.

 Since
$$
\sum_{i=1}^{m}\gamma_i\,  \pi (e )\circ \pi (u_i xv_i )\circ \pi
(e ) = \pi \left(\sum_{i=1}^{m}\gamma_i e \star {}^{u_i xv_i}e \right)
$$
assertion ii) is a consequence of i).

 Via ${\bf \Upsilon}^{-1}$~: $e\star\HH (G)\star e\longrightarrow \HH
 (G,\rho )\otimes_{\Cdss}{\rm End}_{\Cdss}(\Lambda )$, an element
 $\varphi\in e\star \HH (G)\star e$ corresponds to the element
 of  $\HH (G,\rho )\otimes_{\Cdss}{\rm End}_{\Cdss}(\Lambda )$ given
 as follows (see the proof of [BK] Proposition (4.2.4), pages
 149--150). Fix a basis $\{ w_1 ,...,w_n\}$ of $\Lambda$ and let $\{
 {\check w}_1 ,...,{\check w}_n\}$ be the corresponding dual basis of
 ${\check \Lambda}$, so that
 $\langle w_i ,{\check w}_j\rangle =\delta_{ij}$ (Kronecker's delta
 symbol). For each pair of indices $(i,j)$ and for $g\in G$, define an
 operator $\Phi_{ij}(g)\in {\rm End}_{\Cdss}({\check \Lambda})$ by the
 formula:
$$
\langle w,\Phi_{ij}(g){\check w}\rangle =\int_{K}\int_{K}\varphi
(kgl)\langle \rho (l)w_i ,{\check w}\rangle \langle
w,{\check \rho}(k^{-1}){\check w}_j\rangle dkdl\ ,\leqno{\hbox{(1)}}
$$
for all $w\in\Lambda$, $\check w\in\check\Lambda$.
Then the function $g\mapsto \Phi_{ij}(g)$ lies in $\HH
(G,\rho )$, and we have
$$
{\bf \Upsilon}^{-1}(\varphi )={{\rm dim}(\rho )\over \mu
(K)^2}\sum_{i,j=1}^n\Phi_{ij}\otimes w_j \otimes {\check w}_i \ .
$$

 Assume now that $\varphi\in e\star\HH (G)\star e$ has support in
 $KxK$. Then from formula (1), the $\Phi_{ij}$ have support in
 $KxK$. We need the following result.

{\bf (X.1.3) Lemma}. {\it i) The $\Cdss$-vector space
$$
\{ \Phi\in \HH (G,\lm )\ ; \ {\rm Support}(\Phi)\subset J_{\rm
max}xJ_{\rm max}\}
$$
has dimension $1$.

ii) Any non-zero $\Phi$ in $\HH (G,\lm )$ with support $J_{\rm
max}xJ_{\rm max}$ is invertible.}

{\it Proof}. By [BK1](5.5), the $G$-intertwining of $\lambda$ and
$\lambda_{\rm max}$ are $JG_L J$ and\break
 $J_{\rm max}G_LJ_{\rm max}$ respectively,
 where $L/E$ is the
unramified extension introduced in {\S}VIII and $G_L$ the centralizer of
$L$ in $G$. Moreover by [BK1](5.5.13), there is a canonical algebra isomorphism
$\HH (G,\lambda )\lra \HH (G,\lambda_{\rm max})$ which preserves supports in
the following sense: if $y\in G_L$ and $\varphi\in \HH (G,\lambda )$ has support $JyJ$, then its
image $\varphi '\in \HH (G,\lambda_{\rm max})$ has support $\Jm y\Jm$. Moreover consider the Iwahori
subgroup of $G_L$ given by $I_L =U(\Cfr_{\rm min}) =U(\Bfr_{\rm min}\cap G_L$ and let
 $\HH_0 =\HH(G_L ,I_L )$ be the corresponding affine Hecke algebra of type $A$ formed of (locally
constant) bi-$I_L$-invariant compactly supported functions on $G_L$. By Theorem (5.6.6) of
[BK1], the algebras $\HH (G ,\lambda )$ and $\HH_0$ are isomorphic in a support preserving way:
there is (a non-canonical) isomorphism of $\Cdss$-algebras $\Psi$: $\HH_0\lra \HH (G,\lambda )$
such that for all $y\in G_L$ and for all $\varphi \in \HH_0$ with support $I_L y I_L$,
$\Psi (\varphi )$ has support $JyJ$. As a consequence, there exists an algebra isomorphism
$\Psi '$: $\HH_0 \lra \HH (G,\lm )$ enjoying the same support preservation property.

 Now assertions i) and ii) of our lemma hold for the corresponding assertions hold true for the
standard affine Hecke algebra $\HH_0$. Indeed if $y\in G_L$, we have:
\medskip

i) $\{ \varphi \in \HH_0 \ ; \ {\rm Support}(\varphi )\subset I_L yI_L\}$ is the line spanned
by the characteristic function of $I_L yI_L$,
\smallskip

ii) it is a standard fact that any $\varphi\in \HH_0$ with support $I_L yI_L$ is invertible.
\bigskip

Let us fix a non-zero element $\Phi_0$ in $\HH (G,\rho )$ with support $KxK$. Then
$$
{\bf \Upsilon}^{-1}(e\star {}^x e)={{\rm dim}(\rho )\over \mu
(K)^2}\Phi_0\otimes (\sum_{i,j=1}^n\gamma_{ij} w_j\otimes{\check w}_i)
$$
where $\gamma_{ij}$ is defined by $\Phi_{ij} =\gamma_{ij}\Phi_0$,
$i,j\in \{ 1,...,n\}$.  For the same reason, for all $u,v\in K$, there
exists a vector $\zeta (u,v)\in \Lambda\otimes {\check \Lambda}$ such
that
$$
{\bf \Upsilon}^{-1}(e\star {}^{uxv}e)={{\rm dim}(\rho )\over \mu
(K)^2}\Phi_0 \otimes   \zeta (u,v)\ .
$$

{\bf (X.1.4) Lemma}. {\it For all $u,v\in K$, we have
$$
\zeta (u,v)=[\rho (u)\otimes {\check \rho}(v^{-1})]  \zeta (1,1) \ .
$$
}

 Take this last lemma for granted. Since the representation
 $\rho\otimes {\check\rho}$ of $K\times K$ in $\Lambda\otimes
 {\check \Lambda}$ is irreducible, it is generated by the non-zero vector
$\zeta (1,1)$. We may find $m\geq 1$, $u_i , v_i \in K$,
 $\gamma_i\in \Cdss$, $i=1,...,m$, such that
 $\displaystyle \sum_{i=1}^{m}\gamma_i \zeta (u_i ,v_i )$ is an
 arbitrary element of $\Lambda\otimes {\check \Lambda }\simeq {\rm
 End}_{\Cdss}(\Lambda )$. In particular we may choose this element
invertible in ${\rm End}_{\Cdss}(\Lambda )$.  It follows that
$$
{\bf \Upsilon}^{-1}(\sum_{i=1}{m}\gamma_i e\star {}^{u_i xv_i}e)
$$
is invertible. This finishes the proof of Proposition (X.1.2)(ii).

{\it Proof of Lemma (X.1.4)}. The proof is somewhat technical but
straightforward. It is inspired from the calculation of [BK], pages
232--233.

 Write $\Phi_{ij}^{uv}\in \HH (G,\rho )$ for the functions attached to
 $\varphi =e\star {}^{uxv}e$ via formula (1). For $g\in G$, we have
$$
\varphi (g)=\int_K e_{\rho}(y) e_{\rho}((uxv)^{-1} y^{-1}g)dy
\qquad\qquad\qquad\qquad\qquad
$$
$$
\qquad\qquad
={{\rm dim}(\rho )^2 \over \mu (K)^2}\sum_{b,c=1}^n\int_{K} \langle \rho
(y^{-1})w_b ,{\check w}_b\rangle . \langle \rho (g^{-1}yuxv)w_c
,{\check w}_c\rangle dy \ .
$$
So for $w\in \Lambda$ and ${\check w}\in {\check \Lambda}$, we have
$$
{\mu (K)^2\over {\rm dim}(\rho )^2}
\langle w, \Phi_{ij}^{uv}{\check w}\rangle =
$$
$$
\sum_{b,c=1}^n \!\!
\int_{K^3} \!\!\!
\langle  \rho (y^{-1})w_b ,{\check w}_b\rangle
\langle \rho (l^{-1}g^{-1}k^{-1}yuxv)w_c ,{\check w}_c\rangle
\langle \rho (l)w_i ,{\check w}\rangle
 \langle w,{\check \rho}(k^{-1}){\check w}_j\rangle
dkdldy
$$
Integrating with respect to $l$ and using the Schur
orthogonality relation, we obtain:
$$
{\mu (K)\over {\rm dim}(\rho )}
\langle w, \Phi_{ij}^{uv}{\check w}\rangle =
$$
$$
\sum_{b,c=1}^n
\int_{K^2}
\langle  \rho (y^{-1})w_b ,{\check w}_b\rangle
\langle \rho (g^{-1}k^{-1}yuxv)w_c ,{\check w}\rangle
\langle w_i ,{\check w}_c\rangle
\langle w,  {\check \rho}(k^{-1}){\check w}_j\rangle
dkdy
$$
$$
=\sum_{b=1}^n\int_{K^2}
\langle  \rho (y^{-1})w_b ,{\check w}_b\rangle
\langle \rho (g^{-1}k^{-1}yuxv)w_i ,{\check w}\rangle
\langle w,  {\check \rho}(k^{-1}){\check w}_j\rangle
dkdy
$$
We now make the change of variable
$(k')^{-1}=k^{-1}yu$ and this last expression becomes:
$$
\sum_{b=1}^n
\int_{K^2}
\langle  \rho (y^{-1})w_b ,{\check w}_b\rangle
\langle \rho (g^{-1}(k')^{-1}x)\rho(v)w_i ,{\check w}\rangle
\langle w,{\check \rho}((k')^{-1}u^{-1}y^{-1}){\check w}_j\rangle
dk'dy
$$
$$
=\sum_{b=1}^n\int_{K^2}
\langle  \rho (y^{-1})w_b ,{\check w}_b\rangle
\langle \rho (g^{-1}(k')^{-1}x)\rho(v)w_i ,{\check w}\rangle
\langle \rho(y)\rho (uk')w,{\check w}_j\rangle dk'dy \ .
$$
Using again the Schur orthogonality relation, we obtain:
$$
\langle w, \Phi_{ij}^{uv}(g) {\check w}\rangle
=\sum_{b=1}^n\int_{K}
\langle w_b , {\check w}_j\rangle
\langle \rho (uk')w, {\check w}_b\rangle
\langle\rho (g^{-1}(k')^{-1}x)w_i , {\check w}\rangle
dk'
$$
$$
\qquad\quad\ \,=\int_K
\langle \rho (k)w, {\check \rho}(u^{-1}){\check w}_j\rangle
\langle \rho (g^{-1}k^{-1}x)\rho (v)w_i , {\check w}\rangle dk
$$
Let $(V_{ij})$ (resp. $(U_{ij})$) be the matrix of $\rho (v)$
(resp. ${\check \rho}(u^{-1})$) in the basis $\{ w_i\}$ (resp. in the
basis $\{ {\check w}_i\}$).  We have
$$
\langle w, \Phi_{ij}^{uv}(g) {\check w}\rangle
=\sum_{\alpha ,\beta=1}^n V_{\alpha i}U_{\beta j}\int_K
\langle \rho (y)w, {\check w}_{\beta}\rangle
\langle \rho (g^{-1}y^{-1}x)w_{\alpha}, {\check w}\rangle
dk
$$
$$
=\sum_{\alpha ,\beta=1}^nV_{\alpha i}U_{\beta j} \langle
w, \Phi_{\alpha\beta}^{11}{\check w}\rangle\ .
$$
In other words, we have proved that
$$
\Phi_{ij}^{uv}=\sum_{\alpha ,\beta=1}^nV_{\alpha i}U_{\beta
j}\Phi_{\alpha , \beta}^{11}=
\left(\sum_{\alpha ,\beta=1}^n
V_{\alpha i}U_{\beta j}\gamma_{\alpha\beta}\right)
\Phi_0\ .
$$
Hence we obtain:
$$
{\bf \Upsilon}^{-1}(e\star {}^{(uxv)} e)={{\rm dim}(\rho )\over \mu
(K)^2}
\Phi_0\otimes \sum_{i,j=1}^n\sum_{\alpha ,\beta=1}^n
V_{\alpha i}U_{\beta j}\gamma_{\alpha\beta} w_j \otimes{\check w}_i
$$
$$
={{\rm dim}(\rho )\over \mu (K)^2}
\Phi_0\otimes \sum_{\alpha,\beta=1}^n\gamma_{\alpha\beta}
\left(\sum_{j=1}^n U_{\beta j}w_j \right)
\otimes\left(\sum_{i=1}^nV_{\alpha i}{\check w}_i\right)
$$
$$=
 {{\rm dim}(\rho )\over \mu (K)^2}
\Phi_0\otimes \sum_{\alpha,\beta=1}^n\gamma_{\alpha\beta}
\rho (u)(w_\beta )\otimes {\check \rho}(v^{-1})({\check w}_\alpha)
$$
$$
=
 {{\rm dim}(\rho )\over \mu (K)^2}
\Phi_0\otimes
[\rho (u)\otimes {\check \rho }(v^{-1})]\left(
\sum_{\alpha,\beta=1}^n\gamma_{\alpha\beta}
w_\beta \otimes {\check w}_\alpha\right)
$$
as required. (We have used that the matrix of $\rho(u)$ with respect to
the basis $\{w_i\}$ is the transpose of the matrix of $\check\rho(u^{-1})$
with respect to the dual basis $\{\check w_i\}$.) This finishes the proof of Lemma (X.1.4).

\bigskip

{\bf X.2. Orientation of $X[L]$.}
\bigskip

 In order to work with a simpler version of the chain complex of
 {\S}IX, we are going to show that, as a simplicial complex, $X[L]$ has
 a $G^{\circ}$-invariant labelling, where
$$
G^{\circ}= \{ g\in G\ ; {\rm Det}(g)\in \ofr_{F}^{\times}\}
$$

 Recall [Brown] that a labelling of a $d$-dimensional simplicial
 complex $Y$ is a simplicial map $\bf l$~: $Y\longrightarrow
 \Delta_{d}$, from $Y$ to the standard $d$-dimensional simplex, such
 that ${\rm dim} ({\bf l}(\sigma ))={\rm dim}(\sigma )$ for any
 simplex $\sigma$ of $Y$.

 Fix a chamber $C$ of $X$. It is classical that the action of
 $G^{\circ}$ on $X$ has the  following property: any simplex
 $\sigma$ of $X$ has a unique $G^{\circ}$-conjugate that lies in
 (the closure of) $C$. In particular the stabilizer of $C$ in
 $G^{\circ}$ fixes $C$ pointwise. Even though $X[L]$ is not a building
 in general,  its maximal simplices  have the same dimension and we
 call them {\it chambers}.

    We now fix the chamber $C$ so that $C\cap X_{L} \not= \emptyset$.
It is false in general that $G^{\circ}$ acts transitively on the
 chambers of $X[L]$. For instance, if $L/F$ is a maximal unramified
 extension of $F$ in $A$, then $X[L]$ is $0$-dimensional and consists
 of the vertices of $X$.  But the action of $G^{\circ}$ on the
 vertices of $X$ is not transitive.

 Let us notice that $C\cap X[L]$ is a sub-simplicial complex of $X[L]$.
Indeed, passing to the first barycentric subdivisions, we first have that
${\rm sd}(C)\cap {\rm sd}(X[L]) ={\rm sd}(C)\cap X(L)$ is a sub-simplicial
complex of $X(L)$. To get our assertion, it suffices to prove that if
${\rm sd}(C)\cap X(L)$ contains a vertex $x_\sigma$ corresponding to the
isobarycenter of a simplex $\sigma$ of $X[L]$, then $\sigma\subset C\cap X[L]$.
The interior $\sigma^\circ$ of $\sigma$ is of the form $\Sigma^\circ \cap X[L]$,
where $\Sigma$ is some simplex of $X$. We have
 $x_\sigma\in \sigma^\circ \subset \Sigma^\circ$ and $x_\sigma\in C\cap X[L]\subset C$. In
particular $\Sigma^\circ \cap C \not= \emptyset$ and this forces the containment
$\Sigma\subset C$. Therefore $\sigma\subset C$ as required.
\bigskip

 {\bf (X.2.1) Lemma}. {\it The simplicial subcomplex  $C\cap X[L]$ of
 $X[L]$ is a disjoint  union of $f(L/F)$  chambers of $X[L]$.}
\bigskip

{\it Proof}. First we prove that any vertex of $C\cap X[L]$ is
contained in a chamber of $X[L]$ which is itself contained in
$C$. Let $s$ be such a vertex. There exist a field extension
$L'/F\subset A$ and an order $\Afr\subset A$ such  that:

 -- $e(L'/F)=e(L/F)$ and $f(L'/F)=f(L/F)$,

 -- the order $\Afr$ lies in ${\rm Her}(A)^{L'^{\times}}$,

 -- $s$ is the vertex of $X_{L'}$ attached to the maximal order $\Afr
    \cap {\rm End}_{L'}(V)$.

 Let $(N_{k})_{k\in \Zdss}$ be a chain in $V$ corresponding to
 $\Afr$. It must have $\ofr_{L'}$-period $1$, whence it has
$\ofr_F$-period $e(L'/F)$. Assume that $C$ corresponds to a
 lattice chain $(L_k )_{k \in \Zdss}$ in V of $\ofr_{F}$-period $N$.
 There exists an integer $k_o$ such that:
$$
N_{k}=L_{k_o +kN/e(L/F)}\ , \ k \in \Zdss\ .
$$
Since $f(L/F)$ divides $N/e(L/F)$, we have the containments:
$$
\{ L_{k_o +k.N/e(E/F)}\ ; \ k\in \Zdss \}\subset \{ L_{k_o +kf(L/F)}\
; \ k \in \Zdss\}\subset \{ L_k\ ; \ k\in \Zdss\} \ .
$$
By the numerical  criterion of (I.3.5), the set of lattices $\{ L_{k_o +kf(L/F)}\
; \ k \in \Zdss\}$ corresponds to a chamber $C_L$ of $X[L]$ and the
previous containments mean that $s\in C_L \subset C$.

 Let $C_L$ be a chamber of $X[L]$. There exist a field
extension $L'/F\subset A$ and an order $\Afr \subset A$ such that:

 -- $e(L'/F)=e(L/F)$ and $f(L'/F)=f(L/F)$,

 -- the order $\Afr$ lies in ${\rm Her}(A)^{L'^{\times}}$,

 -- $C_L$ is the chamber of $X_{L'}$ attached to the order $\Afr\cap End_{L'}(V)$.

 Let $(M_k )_{k\in \Zdss}$ be a chain in $V$ corresponding to $\Afr$
 and $\Bfr$. It has $\ofr_{L'}$-period $N/[L:F]$. So it has
 $\ofr_F$-period $e(L/F). N/[L:F]=N/f(L/F)$. Moreover, for all $k\in
 \Zdss$, we have:
$$
{\rm dim}_{\Fdss_{F}}(M_k /M_{k+1})=f(L/F)  {\rm dim}_{\Fdss_{L'}}(M_k
/M_{k+1}) = f(L/F).1 =f(L/F)\ .
$$

 Assume now that $C_L$ lies in the  chamber $C$ of $X$.
According to the previous discussion, there exists a coset $\Gamma$ of
 $f(L/F)\Zdss /N \Zdss$ in $\Zdss /N\Zdss$, such that:

$$
\{ M_k\  ;\ k\in \Zdss\}= \{ L_{l}\ ; \ l\in \Zdss\ {\rm and}\  l\ {\rm mod} \
 N\Zdss \in \Gamma\}\ .
$$

Conversely, using proposition (I.3.5), we have that  for all such
coset $\Gamma$, the lattice chain whose lattice set is given by
$ \{ L_{l}\ ; \ l\in \Zdss\ {\rm and}\  l\ {\rm mod} \ N\Zdss \in
\Gamma\}$
correspond to a chamber of $X[L ]$ contained in $C$. Indeed if $F_{\Gamma}$ is the
simplex of $X$ corresponding
to the lattice set  $ \{ L_{l}\ ; \ l\in \Zdss\ {\rm and}\  l\ {\rm mod} \ N\Zdss \in
\Gamma\}$, then the corresponding (closed) chamber of $X[L]$ is
$F_{\Gamma}\cap X[L]$. When $\Gamma$ runs over the $f(L/F)$ cosets of
$f(L/F)\Zdss /N \Zdss$ in $\Zdss /N\Zdss$, the corresponding (closed) chambers
are disjoint, as required.

 Since the simplicial complex $X[L]\cap C$ is a disjoint union of
 (closed) chambers, it is trivially labelable. Let us fix a labelling
 ${\bf l}_C$~: $X[L]\cap C\longrightarrow \Delta_C$, where
 $\Delta_{C}$ is the standard simplex of dimension ${\rm
 dim}X[L]=N/[L:F]-1$. For any simplex $\sigma$ of $X[L]$, we define
 a simplex ${\bf l}(\sigma )$ of $\Delta_C$ by ${\bf l}(\sigma )={\bf
 l}_{C}(\sigma_C )$, where $\sigma_C$ is the unique simplex of
 $X[L]\cap C$ which is a conjugate of $\sigma$ under the action of
 $G^{\circ}$.

 {\bf (X.2.2) Lemma}. {\it The map $\bf l$~: $X[L]\longrightarrow
 \Delta_C$ is a labelling. It is invariant under the action of
 $G^{\circ}$}.

 {\it Proof}. Obvious from the properties of the action of $G^{\circ}$ on
 $X$.

 From now on, we fix the $G^{\circ}$-invariant labelling $\bf l$ of
 $X[L]$ (by fixing ${\bf l}_C$). It gives rise to a
 $G^{\circ}$-invariant orientation of the simplicial complex $X[L]$ as
 well as $G^{\circ}$-invariant incidence numbers $[\sigma :\tau ]$ for
 any pair of simplices $\tau\subset \sigma $ of $X[L]$ with $\tau$ of
 codimension $1$ in $\sigma$.
\bigskip

{\bf X.3 Another chain complex.}
\bigskip

 We fix a smooth complex representation $(\pi ,\VV)$ in   $\RR_{(J,\lambda
 )}G$ and consider the coefficient system $\CC = (\VV [\sigma])_{\sigma}
 =\CC_{(J,\lambda )}(\VV )$ of {\S}VIII.

 For $q=0,\dots ,N/[L:F] -1$, let $X[L]_q$ denote the set of
 $q$-simplices of $X[L]$. The space $C_{q} (X[L],\CC )$ of
 (unoriented) $q$-chains of $X[L]$ with coefficient in $\CC$ is the
 $\CC$-vector space of all maps $\omega$~: $X[L]_{q}\longrightarrow
 \VV$ such that $\omega$ has finite support and $\omega (\sigma )\in
 \VV [\sigma ]$, for all $\sigma\in X[L]_q$. The group $G$ acts
 smoothly on $C_{q}(X[L],\CC )$ via $(g\omega )(\sigma ):=g(\omega
 (g^{-1}\sigma ))$. The orientation of $X[L]$ gives rise to boundary
 maps:

$$
\matrix{
\partial~:\ C_{q+1}(X[L],\CC ) & \rightarrow &
            C_{q}(X[L],\CC) \hfill\cr\cr
   \hfill \omega & \mapsto & [\ \ \sigma\mapsto
\mathop{\sum}\limits_{\tau\in X[L]_{q}\ , \ \tau\subset\sigma} [\sigma
   :\tau] \omega (\tau)\ \ ] }
$$
We obtain an augmented chain complex of $G^{\circ}$-modules:
$$
C_{N/[L:F]-1} (X[L],\CC )\boundary \cdots \boundary C_0 (X[L],\CC )
\augmentation \VV \leqno (\hbox{X.3.1})
$$
where $\epsilon (\omega )=\sum_{\sigma\in X[L]_0}\omega (\sigma )\in
\VV$.

{\bf (X.3.2) Lemma}. {\it As augmented chain complexes of
$G^{\circ}$-modules,  the complexes (IX.1) and (X.3.1) are canonically
isomorphic.}

{\it Proof}. By standard arguments.
\bigskip

{\bf X.4 $\Jm$-orbits of simplices.}
\bigskip

 Fix $q\in \{ 0,\dots ,N[L:F]-1 \}$. For any subset $\Sigma$ of
 $X[L]_q$, we denote by $C_{q}(\Sigma ,\CC )$ the subspace of $C_q
 (X[L] ,\CC )$ formed of those $q$-chains with support in $\Sigma$.

 Let $\Omega_q$ be the set of orbits of $\Jm$ in $X[L]_q$. As a
 $\Jm$-module, $C_q (X[L],\CC )$ decomposes as
$$
C_q (X[L],\CC ) = \coprod_{\Sigma\in \Omega_q} C_q (\Sigma ,\CC )\ .
$$
Fix $\Sigma\in \Omega_q$. There exist $\Cfr\in {\rm Her}(C)$
satisfying $\Cfr_{\rm min}\subset \Cfr \subset \Cfr_{\rm max}$ and
$x\in G$ such that $\Sigma = \Jm x.\sigma_{\Cfr}$.We have the disjoint
union:
$$
\Sigma =\bigcup_{j\in\Jm  /\Jm\cap U(\Afr )^{x}} \{ jx\sigma_{\Cfr
}\}\ ,
$$
where $\Afr =\Afr (\Bfr )$ and $\Bfr =\Bfr (\Cfr )$, from which we
deduce the following isomorphisms of $\Jm$-modules:
$$
C_q (\Sigma ,\CC )=\coprod_{j\in \in\Jm  /\Jm\cap U(\Afr )^{x}} C_q
(jx\sigma_{\Cfr},\CC )= \coprod_{j\in \in\Jm  /\Jm\cap U(\Afr )^{x}} jxC_q
(\sigma_{\Cfr},\CC )\ .
$$
We have a natural $\Jm$-homomorphism $S_{\Sigma}$~: $C_{q}(\Sigma ,\CC
)\longrightarrow \VV$, given by
$$
S_{\Sigma}(\omega )=\sum_{\sigma\in\Sigma}\omega (\sigma )\ .
$$
In other words:
$$
S_{\Sigma}(\bigoplus_{j\in\Jm /\Jm\cap U(\Afr )^x} jx\omega_j )
=\sum_{j\in\Jm /\Jm\cap U(\Afr )^x} jx \omega_j (\sigma_{\Cfr})\ , \
\omega_j\in \VV [\sigma_{\Cfr}]\ .
$$
We set $K_{\Sigma}={\rm Ker}S_{\sigma}$. We have the following exact
sequences of $\Jm$-modules and $\Cdss$-vector spaces respectively:
$$
0\longrightarrow K_{\Sigma} \longrightarrow C_{q}(\Sigma, \CC
)\longrightarrow \sum_{j\in \Jm}jx\VV[\sigma_{\Cfr}]\longrightarrow 0
$$
$$
0\longrightarrow K_{\Sigma}^{\lm} \longrightarrow C_{q}(\Sigma, \CC
)^{\lm}\longrightarrow \big(\sum_{j\in \Jm}jx\VV[\sigma_{\Cfr}]\big)^{\lm}\longrightarrow 0
$$

Moreover, by lemmas (IX.4) and (IX.5), we have
$$
\VV [\sigma_{\Cfr}]=\sum_{g\in U(\Afr )/\Jm}g\VV^{\eta (\Bfr_{\rm
min},\Bfr_{\rm max})}=\sum_{g\in U(\Afr )/\Jm}g\VV^{\lm}\ .
$$
Therefore we have
$$
\sum_{j\in\Jm}jx\VV [\sigma_{\Cfr}]=\sum_{j\in\Jm}\sum_{g\in U(\Afr
)}jxg\VV^{\lm}\subset \VV\ .
$$
By proposition (X.1.1), for all $j\in \Jm$, $g\in U(\Afr )$, we have
$$
\em\star \{ jxg\VV^{\lm} \} =\left\{
\matrix{
\VV^{\lm} & {\rm if}\ jxg\in I_{G}(\lm )\cr
0 & {\rm otherwise}
}\right.
$$
We deduce that
$$
\big(\sum_{j\in \Jm}jx\VV [\sigma_{\Cfr}]\big)^{\lm}=
\left\{\matrix{
\VV^{\lm} & {\rm if} \ \exists \ g\in U(\Afr),\ j\in \Jm\ {\rm
s.t.}\ jxg\in I_{G}(\lm )\cr
0 & {\rm otherwise}
}\right.
$$
Since the $G$-intertwining of $\lm$ is $\Jm  G_L\Jm$, this may be
rewritten:
$$
\big(\sum_{j\in \Jm}jx\VV [\sigma_{\Cfr}]\big)^{\lm}=
\left\{\matrix{
\VV^{\lm} & {\rm if} \ x\in \Jm G_L U(\Afr) \cr
0 & {\rm otherwise}
}\right.
$$

{\bf (X.4.1) Conjecture}. {\it For any $\Sigma \in \Omega_q$,
we have $K_{\Sigma}^{\lm}=0$.}

{\bf (X.4.2) Corollary}. {\it Assume that conjecture (X.4.1) holds.

 i) If $\Sigma \cap X_L \not=\emptyset$, then
$S_{\Sigma}$ induces an isomorphism of $\Cdss$-vector spaces: \break
$C_q (\Sigma ,\CC )^{\lm}\longrightarrow \VV^{\lm}$.

ii) If $\Sigma\cap X_L =\emptyset$, then $C_{q}(\Sigma ,\CC )^{\lm} =0$.}

 Indeed we have $\Jm x\sigma_{\Cfr} \cap X_{L}\not=\emptyset$ if and
 only if there exist $\Cfr '\in {\rm Her}(C)$ and $j\in \Jm$ such that
 $jx\sigma_{\Cfr}=\sigma_{\Cfr '}$. By lemma (I.3.1) and (I.3.3), this
 is equivalent to the existence of $z\in G_L$ and $g\in U(\Afr )$ such
 that $jx =zg$, as required.

 From now on we fix an apartment $\AL$ of $X_L$ containing the
 chamber $\sigma_{\Cfr_{\rm min}}$.

{\bf (X.4.3) Lemma}. {\it Let $\Sigma\in \Omega_q$. Assume that
$\Sigma\cap X_{L}\not=\emptyset$. Then $\Sigma\cap \AL\not=\emptyset$
and the intersection $\Sigma \cap \AL$ is reduced to a single
simplex. Moreover $\Sigma\cap X_L$ is a single $U(\Cfr_{\rm min})$-orbit.}

If $\Sigma\cap X_L\not=\emptyset$, then $\Sigma=\Jm.\sigma_L$ for some
$\sigma_L\in (X_{L})_q$. Since $J_m$ contains the Iwahori subgroup
$U(\Cfr_{\rm min})$ (of ${\rm Aut}_L (V)$) and that $(\AL )_q$ is a
system of representatives of the $U(\Cfr_{\rm min})$-orbits in $(X_L
)_q$, we have $\Sigma\cap\AL\not=\emptyset$. At this stage we need the
following technical result.

{\bf (X.4.4) Lemma}. {\it Let $\sigma$, $\tau$ be simplices of
$X_L$. Then if they are conjugate under the action of $U(\Afr_{\rm
min})$, they are conjugate under the action of $U(\Cfr_{\rm min})$.}

Lemma (X.4.3) follows from the previous lemma by observing that
 $\Jm$ is contained in $U(\Cfr_{\rm min})U^{1}(\Afr_{\rm max})\subset
 U(\Afr_{\rm min})$.

{\it Proof of Lemma (X.4.4)}. Let $C_0$ be the chamber of $X_L$ fixed
by $U(\Cfr_{\rm min})$. Fix an apartment $\Ap_L$ of $X_L$ containing
$C_0$ and $\sigma$. Let $x_\sigma$ the barycenter of $\sigma$.
 Then there exists a point $x_0\in C_0^{\circ}$ such that the
geodesic segment $[x_0 ,x_\sigma )\subset \Ap_L$ does not intersect any simplex
of $X_L$ of codimension greater than or equal to $2$. Indeed consider the
subsets of $\Ap_L$ of the form $C^\bullet = {\rm Cvx}\{ x_\sigma
,F\}\backslash \{ x_\sigma\}$, where $F$ is a simplex of codimension
greater than or equal to $2$ in $\Ap_L$ and where $\rm Cvx$ denotes a convex
hull. The set of such subsets is countable. Moreover these subsets
 have empty interiors and by Baire's theorem their union has
empty interior. It follows that this union cannot contain $C^{\circ}_0$
as required.

 Let $\Gamma$ be the set of chambers $D$ in $\Ap_L$ such that $D\cap [x_0
 ,x_\sigma )=D^{\circ}\cap [x_0 ,x_\sigma )\not= \emptyset$. Then it is easy to see
 that there exists an indexation $\Gamma =\{ D_i\ ; \ i=0,...,r\}$ of
 the elements of $\Gamma$ such that $(D_i )_{i=0,..,r}$ is a gallery
 satisfying:  $D_0 =C_0$ and $D_r$ contains $x_\sigma$, whence contains
 $\sigma$.  We can be more precise: for
 $i=0,...,r-1$, $D_{i+1}$ is the unique chamber adjacent to $D_i$ and
 intersecting $[y,x]$, where $\displaystyle [x_0 ,y]=[x_0
 ,x]\cap (\cup_{j=0,...,i}D_j )$. Moreover, for $i=0,...,r-1$,
 let $H_i$ be the wall separating $D_i$ and $D_{i+1}$. It defines two
 roots $H_i^{\pm}$ (half-spaces with boundary $H_i$), such that $H_i^-$
 contains $x_0$ and $H_i^+$ contains $x$. Then the gallery $(D_i
 )_{i=0,...,r}$ is constructed in such a way that
 $\displaystyle \bigcup_{j=0}^{i} D_j \subset H_i^-$ and
 $\displaystyle \bigcup_{j=i+1}^r D_j \subset H_i^+$.

 Let $g\in U(\Cfr_{\rm min})$ be such that $g\sigma =\tau$. Then $g$
 fixes $C_0$ pointwise. Recall that by (I.2.3), there exist
 normalizations of metrics on $X_L$ and $X$ such that
 the embedding $X_L\subset X$ is
 isometric. It follows that the set $g[x_0 ,x_\sigma
 ]$ is the geodesic segment in $X_L$ linking $g.x_0 =x_0$ and
 $g.x_\sigma \in X_L$. Recall that $X_L$ is a simplicial subcomplex of
 ${\rm sd}(X)$. For $i=0,...,r$, $gD_i$ is a simplex of ${\rm sd}(X)$
 whose interior intersects $X_L$. So this simplex belongs to $X_L$.
It follows that $(gD_i )_{i=0,...,r}$ is a gallery in $X_L$
 satisfying $gD_0 =C_0$ and $gD_r \supset \tau$.

 We are going to prove by induction on $t\in \{ 0,...,r\}$ that there
 exists $g_t\in G_L$ such that $g_t D_i =gD_i$, $i=0,...,t$. We will
 then have $g_t\in U(\Cfr_{\rm min})$ and $g_t^{-1} gD_r =D_r$. Since
 $g_t^{-1}g\in U(\Cfr_{\rm min})$ is a compact element of $G_L$, it
 must fix $D_r$ pointwise. It will follow that $g_t^{-1}g\sigma
 =\sigma$, that is $g_t \sigma =\tau$, as required.

 The result is obvious when $t=0$. Assume $t\in \{ 0,..,r-1\}$ and
 that the result is proved for $t$. Replacing $\tau$ by
 $g_t^{-1}\tau$, $g$ by $g_t^{-1}g$, we may assume that $gD_i =D_i$,
 $i=0,...,t$. The chamber $D_{t+1}$ does not belong to $H_t^-$ and has
 a codimension $1$ face contained in $H_t$. The chamber $gD_{t+1}$ has
 a codimension $1$ face contained in $H_t$ and does not belong to
 $H_t^-$, otherwise this would contradict the fact that $g[x_0
 ,x_\sigma ]$ is a geodesic segment. Let $\varsigma$ be the
 codimension $1$ simplex $H_t\cap D_{t+1} =H_t \cap gD_{t+1}$.
 Then the pointwise fixator of $H_t^-$ in $G_L$ acts transitively on the
 set of chambers containing $\varsigma$ and not contained in $H_t^-$
 (an easy exercise left to the reader).
 It follows that there exists $g_{t+1}$ fixing $H_t^-$ pointwise such
 that $g_{t+1}D_{t+1}=gD_{t+1}$, as required.

\bigskip

{\bf X.5 Comparison of chain complexes}
\bigskip

 As in the previous section, we fix an apartment $\AL$ containing
 $\sigma_{\Cfr_{\rm min}}$. As a subcomplex of $\AL$, the topological
 space $\AL$ is equipped with its canonical triangulation. We denote
 by $\uVV^{\lm}$ the constant coefficient system on $\AL$ such that
 for any simplex $\sigma$, $\uVV^{\lm}[\sigma ]=\VV^{\lm}$. It gives
 rise to the chain complex $C_{\bullet}(\AL ,\uVV^{\lm})$, with an
 augmentation map: $C_{\bullet}(\AL , \uVV^{\lm})\augmentationL
 \VV^{\lm}$. This complex is exact since the topological space $\AL$
 is contractible (more precisely it is homeomorphic to a finite
 dimensional affine space).  We shall denote by $\partial_L$ the
 boundary maps of that complex.

 Denote by $C_{\bullet}(X[L], \CC )^{\lm}\augmentation \VV^{\lm}$ the
 augmented chain complex obtain by applying the functor of
 $\lm$-isotypic components to the augmented complex (X.3.1) (or
 equivalently to the augmented complex of (IX.1)). It lies in the
 category of left $\em\star \HH (G)\star\em$-modules. Since the
 functor
$$
\matrix{
\RR_{(J,\lambda )}(G) & \longrightarrow  & \em\star \HH (G)\star\em \ -\
 {\rm Mod}\cr
{\cal W} & \mapsto  & {\cal W}^{\lm}
}
$$
is an equivalence of categories, we have, using proposition (IX.2),
that the complex is exact if and only if
$C_{\bullet}(X[L], \CC )^{\lm}\augmentation \VV^{\lm}$ is exact.

{\bf (X.5.1) Proposition}. {\it Assume that the representation $(\pi ,\VV )$
satisfies  conjecture (X.4.1).
The augmented chains complexes
$C_{\bullet}(X[L], \CC )^{\lm}\augmentation \VV^{\lm}$ and $C_{\bullet}(\AL ,
\uVV^{\lm})\augmentation \VV^{\lm}$ are then  naturally isomorphic as
complexes of $\Cdss$-vector spaces.}

{\it Remark}. There is maybe a more precise result to prove. Indeed
there should be a natural  action of the scalar Hecke algebra on $C_{\bullet}(\AL ,
\uVV^{\lm})\augmentation \VV^{\lm}$ such that the complexes are
isomorphic as complexes of $\em\star \HH (G)\star\em$-modules.

 As a corollary, we have:

 {\bf (X.5.2) Theorem }. {\it Let $(J ,\lambda )$ be a simple type
 of $G$. Let  $(\pi ,\VV )$ be a smooth complex representation   in $\RR_{(J,\lambda
 )}(G)$ satisfying conjecture (X.4.1). Then  the augmented chain complex
$$
C_{\bullet}(X[L] , \CC_{(J,\lambda )}(\VV ))\augmentation \VV
$$
is a resolution of $\VV$ in the category  $\RR_{(J,\lambda
 )}(G)$. In particular, as a $G$-module, the space $\VV$ is given by
 the homology module $H_{0}(X[L] , \CC_{(J,\lambda )}(\VV ))$.}

 {\it Proof of proposition (X.5.1)}. We are going to construct a
 natural isomorphism of complexes from
$C_{\bullet}(X[L], \CC )^{\lm}\augmentation \VV^{\lm}$ to
$C_{\bullet}(\AL ,\uVV^{\lm})\augmentation \VV^{\lm}$. This  is a collection
 of isomorphisms : $[(\varphi_{q})_{q\geq 0},  \psi ]$, where
$$
\varphi_q \in {\rm Hom}_{\Cdss}\big( C_{q}(X[L], \CC )^{\lm}, C_{q}(\AL
 ,\uVV^{\lm})\big)\ , \ \psi\in {\rm Hom}_{\Cdss}(\VV^{\lm}, \VV^{\lm})\ ,
$$
and where the obvious square diagrams are commutative. We first take
$\psi$ to be the identity map of $\VV^{\lm}$. To define $\varphi_q$,
we note that
$$
C_{q}(X[L] , \CC )=\coprod_{\Sigma\in \Omega_q} C_q (\Sigma ,\CC )
$$
and that, by corollary (X.4.2)(ii), we have:
$$
C_{q}(X[L] , \CC )^{\lm}=\coprod_{\Sigma\in \Omega_q\ , \ \Sigma\cap
X_L\not=\emptyset}  C_q (\Sigma ,\CC )^{\lm}\ .
$$
For any simplex $\sigma$ of $\AL$, we let $\Sigma_{\sigma}$ denote
the $\Jm$-orbit of simplices through $\sigma$, so that:
$$
 C_{q}(X[L] , \CC )^{\lm}=\coprod_{\sigma\in {(\AL )q}}  C_q
 (\Sigma_{\sigma} ,\CC )^{\lm}\ .
$$
We now define $\varphi_{q}$~: $C_q (X[L], \CC )^{\lm}\longrightarrow
C_q (\AL ,\uVV^{\lm})$ by
$$
\varphi_q (\omega )(\sigma )=S_{\Sigma_{\sigma}}(\omega{\mid\Sigma_{\sigma}})\ , \sigma\in (\AL )_q\
.
$$
By corollary (X.4.2)(i), the map $\varphi_q$ is clearly an isomorphism
of $\Cdss$-vector spaces.

{\bf (X.5.3) Lemma}. {\it Under the assumptions of Theorem (X.5.2),
for \break
$q=1,$ $\dots ,N/[L:F]-1$, the following diagram is commutative:
$$\matrix{
   & C_{q}(X[L] ,\CC )^{\lm} &
   \buildrel{\partial}\over{\longrightarrow} &
C_{q-1} (X[L],\CC )\lm &   \cr
\varphi_q & \downarrow &      &  \downarrow & \varphi_{q-1} \cr
    & C_{q}(\AL ,\uVV^{\lm}) &
    \buildrel{\partial_L}\over\longrightarrow & C_{q-1}(\AL ,
    \uVV^{\lm}) &
}
$$
}
Fix $\omega\in C_q (X[L],\CC )^{\lm}$. We have
$$\varphi_q (\omega )(\beta )=\sum_{\tau\in \Sigma_{\beta}}\omega
(\tau ), \ \beta\in (\AL )_q\ ,
$$
and
$$
\partial_L (\varphi_q (\omega ))(\alpha )=\sum_{\beta\in (\AL )_q \ ,
\ \beta\supset \alpha}\big\{\sum_{\tau\in \Sigma_{\beta}}[\beta :\alpha
]\omega (\tau )\ ,\ \alpha \in (\AL )_{q-1}\ \big\} . \leqno{(E1) }
$$
On the other hand we have
$$
\partial\omega (\sigma )=\sum_{\theta\in X[L]_q ,\ ,
\theta\supset\sigma} [\theta :\sigma ] \omega (\theta
)\ , \sigma\in X[L]_q\ ,
$$
and
$$
\varphi_{q-1}(\partial\omega )(\alpha )=\sum_{\sigma\in
\Sigma_{\alpha}}\big\{\sum_{\theta \in X[L]_q\ , \ \theta\supset
\sigma}[\theta :\sigma ]\omega (\theta )\ , \ \alpha \in (\AL )_{q-1}\
\big\}
.\leqno{(E2)}
$$

 Fix $\alpha\in (\AL )_{q-1}$. The set $\Theta$ of $\theta\in X[L]_q$
 containing some $\sigma\in \Sigma_{\alpha}$ in general strictly contains
 the set of $\tau$ in $X[L]_q$ such that there exists $\beta\in (\AL
 )_q$, $\beta \supset \alpha$ and $\tau\in \Sigma_{\beta}$. However
 the first set $\Theta$ is stable under $\Jm$ and  splits into two
 disjoint subsets:

 -- the  subset $\Theta_1$ of those $\theta$ whose $\Jm$-orbits
    intersect $\AL$;

 -- the complementary subset $\Theta_2$.

 Let $\theta\in \Theta_1$ and $\sigma\in \Sigma_{\alpha}$ such that
 $\theta\supset\sigma$.  We have $\theta\in \Sigma_{\beta}$ for some
 simplex $\beta$ of $\AL$. The simplex $\beta$ contains a
 $\Jm$-conjugate of $\sigma$ lying in $\AL$. By unicity in lemma
 (X.4.3), that simplex must be $\alpha$. In other words $\theta$ lies
 in $\Sigma_{\beta}$ for some $\beta\in (\AL )_q$ containing $\alpha$
 and there is a unique $\sigma\in \Sigma_{\alpha}$ such that
 $\theta\supset \sigma$: if $\theta =j\beta$, $j\in \Jm$, then $\sigma
 =j\alpha$. Since the action of $\Jm$ preserves the incidence numbers,
 we must have $[\theta :\sigma ] =[\beta : \alpha ]$.

 From the previous discussion, we deduce:
$$
\varphi_{q-1}(\partial \omega )(\alpha
)=\partial_{L}(\varphi_{q}(\omega ))(\alpha )+\sum_{\sigma\in
\Sigma_{\alpha}}\big\{ \sum_{\theta\in \Theta_{2}\ , \
\theta\supset\sigma}[\theta :\sigma ]\omega (\theta )\big\}\ .
$$
Note that if $\theta$ is a simplex of $X[L]$, there is at most one
$\sigma\in \Sigma_{\alpha}$ such that $\theta\supset \sigma$. Indeed
two such simplices contained in $\theta$ must be equal since they have
the same label. In other words in the sum $\sigma$ depends in a
$\Jm$-equivariant way from $\sigma$; we shall write $\sigma =\sigma
(\theta )$. Let $\Omega_{q}(\Theta_2 )$ be the set of $\Jm$-orbits in
$\Theta_2$. We may write:
$$
\varphi_{q-1}(\partial \omega )(\alpha) -
\partial_{L}(\varphi_{q}(\omega ))(\alpha )
=\sum_{\Sigma\in \Omega_{q}(\Theta_2 )}\sum_{\theta\in \Sigma}[\theta
:\sigma (\theta )]\omega (\theta )
$$
$$
=\sum_{\Sigma\in \Omega_{q}(\Theta_2
)}\epsilon_{\Sigma}\sum_{\theta\in \Sigma} \omega (\theta )\ ,
$$
where $\epsilon$ is a sign depending only on $\Sigma$. For $\Sigma\in
\Omega_{q}(\Theta_2 )$, the restriction map:
$$
\matrix{
C_{q}(X[L],\CC ) & \longrightarrow  & C_{q}(\Sigma ,\CC ) \cr
\omega & \mapsto & \omega_{\mid \Sigma}
}
$$
is $\Jm$-equivariant and its restriction to $C_{q}(X[L],\CC )^{\lm}$
must have image in $C_{q}(\Sigma ,\CC )^{\lm}$. Since $\Sigma\cap \AL
=\emptyset$, by applying corollary (X.4.2)(ii), we obtain that
$C_{q}(\Sigma ,\CC )^{\lm} = 0$, whence
$\displaystyle \sum_{\theta\in \Sigma} \omega (\theta ) = 0$, for all
$\Sigma\in \Omega_{q}(\Theta_2 )$. Finally we get
$\varphi_{q-1}(\partial \omega )(\alpha) - \partial_{L}(\varphi_{q}(\omega ))(\alpha )
=0$ and the commutativity of the diagram.

 Using a quite similar proof we have the following result.

{\bf (X.5.4) Lemma}. {\it Under the assumptions of Theorem (X.5.2),  the
following diagram is commutative:
$$
\matrix{
   & C_{0}(X[L] ,\CC )^{\lm} &
   \augmentation &
\VV^{\lm}&   \cr
\varphi_0 & \downarrow &      &  \downarrow & {\rm Id} \cr
    & C_{0}(\AL ,\uVV^{\lm}) &
    \augmentationL  & \VV^{\lm} &
}
$$
}

This finishes the proof of proposition (X.5.1) and theorem (X.5.2).

\bigskip

\centerline{\bf XI. Acyclicity in the case of a discrete series
 representation.}
\bigskip

 The aim of this section is to prove conjecture (X.4.1) when
 the representation is irreducible and lies in the discrete series of
 $G$. More precisely we shall assume that our representation $(\pi
 ,\VV )$ is {\it an unramified twist of a (irreducible unitary)
 discrete series representation} of $G$ containing our fixed simple
 type $(J,\lambda )$.  In that case the chain complex attached to
 ${\cal C}_{(J,\lambda )}(\VV )$ may be entirely computed.
\bigskip

{\bf XI.1. Determination of the chain complex}.
\bigskip

 We keep the notation as in section IX. Let $\Cfr$ be a hereditary
 order of the $L$-algebra $C$ satisfying $\Cfr_{\rm
 min}\subset \Cfr \subset \Cfr_{\rm max}$, and let $\sigma_{\Cfr}$ be
 the corresponding simplex in $X_L \subset X[L]$. We want to
 understand the $U(\Afr )$-module structure of
$$
\VV [\sigma_\Cfr ]=\sum_{g\in U(\Afr )/U(\Bfr )J^1 (\Bfr_{\rm max} )}
 g. \VV^{\eta (\Bfr ,\Bfr_{\rm max})}\ .
$$
Let $\cal W$ an irreductible constituent of the $U(\Bfr )J^1 (\Bfrm
)$-module $\VV^{\eta (\Bfr ,\Bfrm )}$. By Frobenius reciprocity $\cal W$
embeds in a representation of the form $\kappa_{\rm max}\otimes\tau$,
where $\tau$ is an irreducible representation of $U(\Bfr )U^1 (\Bfr )$
seen as a representation of $U(\Bfr )J^1 (\Bfrm )$ trivial on
$U^1 (\Bfr )J^1 (\Bfrm )$. The following result implies that $\cal W$
actually has the form  $\kappa_{\rm max}\otimes\tau$.

{\bf (XI.1.1) Lemma}.  {\it If $\tau$ is an irreducible representation
of $U(\Bfr )/U^1 (\Bfr )$, then the  $U(\Bfr )J^1 (\Bfrm )$-module
$\kappa_{\rm max}\otimes\tau$ is irreducible.}

{\it Proof}. By Schur Lemma, it suffices to prove that ${\rm End}_{U(\Bfr )J^1
(\Bfrm )}\, \kmax\otimes  \tau$ is one-dimensional.
 For this we closely follow the proof of [BK]
Proposition (5.3.2)(ii), page 176. Write $X$ for the representation
space of $\kmax$ and $Y$ for the representation space of $\tau$.
Let $\varphi\in {\rm End}_{U(\Bfr )J^1 (\Bfrm )}\, \kmax\otimes \tau$
 that we write $\varphi =\sum_j S_j \otimes T_j$, where $S_j \in {\rm End}_\Cdss\, X$, $T_j\in {\rm
End}_\Cdss \, Y$, and where the $T_j$ are linearly independent.  For
$h\in J^1 (\Bfrm )$ we have
$$
(\kmax \otimes \tau )(h)\circ \varphi =\varphi\circ (\kmax \otimes \tau
)(h)\ .
$$
Since $J^1 (\Bfrm )\subset {\rm Ker}\, (\xi )$, we obtain
$$
\sum_j (\kmax (h)\circ S_j -S_j \kmax (h))\otimes T_j =0\ .
$$

Since the $T_j$ are linearly independent, we obtain that $S_j\in {\rm
End}_{J^1 (\Bfrm )} \, \eta (\Bfrm )$ for all $j$. But since
$\eta (\Bfrm )$ is irreducible, we have that ${\rm End}_{J^1 (\Bfrm )}\, \eta
(\Bfrm )$ and  ${\rm End}_{U(\Bfr )J^1 (\Bfrm )}\, \kmax$ are equal and
one-dimensional. So we may as well take $j=1$, so that $\varphi
=S\otimes T$, where $S\in {\rm End}_{U(\Bfr )J^1 (\Bfrm )}\, \kmax$
and $T\in {\rm End}_\Cdss \, Y$. Now any $h\in U (\Bfr )J^1 (\Bfrm )$
must satisfy
$$
(S\circ \kmax (h))\otimes (T\circ \xi (h)) =(\kmax (h)\circ S)\otimes
(T\circ \xi (h))=(\kmax (h)\circ S)\otimes (\xi (S)\circ T)\ .
$$
But this implies that $T\in {\rm End}_{U(\Bfr )}\, \tau$ and our result
follows from the irreducibility of $\tau$ and Schur Lemma.
\bigskip

First consider the case $\Cfr =\Cfr_{\rm max}$, so that we have
$U(\Bfr )J^1 (\Bfr_{\rm max})=J(\Bfr_{\rm max})$ and $\eta (\Bfrm
,\Bfrm ) =\eta (\Bfrm )=(\kappa_{\rm max})_{\vert J^{1} (\Bfrm
)}$. Since $\VV$ is admissible, $\VV^{\eta (\Bfrm )}$ is finite
dimensional and, as a $J(\Bfrm )$-module, decomposes as a finite sum
of irreducibles submodules. By lemma (XI.1.1), these irreducible representations have the
form $\kmax\otimes \tau$, where $\tau$ is an irreducible
representation of $J(\Bfrm )/J^1 (\Bfrm )\simeq U(\Bfrm )/U^1 (\Bfrm
)$. Moreover by [SZ] (see the discussion preceeding Lemma 2, page
176), for such a $\tau$, we have:
\medskip

(1) \ \ \  ${\rm Hom}_{J(\Bfrm )}\, (\kmax\otimes \tau , {\cal V})={\rm Hom}_{U(\Bfrm )/U^1
(\Bfrm )}\, (\tau , \VV (\Bfrm ))\ .$
\medskip

Recall that $\VV (\Bfrm )$ is the $U(\Bfrm )/U^1
(\Bfrm )$-module ${\rm Hom}_{J^1 (\Bfrm )}(\kmax ,\VV )$.

 By considering $L_{\Bfr_0}=U(\Bfr_0 )/U^1 (\Bfr_0 )$ as a Levi
 subgroup of
${\bar G}=U(\Bfrm )/\allowbreak U^1 (\Bfrm )$, we may form the
 {\it generalized Steinberg representation} $\St (\Bfrm ,\rho )$ with
 cuspidal support $(L_{\Bfr_0}, \rho )$. It may be defined in several
 ways. In particular it is the unique generic sub-$\bar G$-module of
 the representation of $\bar G$ parabolically induced from
 $(L_{\Bfr_{0}}, \rho )$. We then have the following crucial result.
\medskip

{\bf  (XI.1.2) Lemma}. ([SZ], Proposition 6, page 179.) {\it As a
$\bar G$-module, $\VV (\Bfrm )$ is isomorphic to $\St (\Bfrm ,\rho )$.}
\medskip

It follows from (1) and the previous lemma that the space ${\rm
Hom}_{J(\Bfrm )}\, (\kmax\otimes \tau ,\VV )$ is zero except when
$\tau \simeq \St (\Bfrm ,\rho )$ where it is $1$-dimensional. We have
proved the following result.
\medskip

{\bf (XI.1.3) Lemma}. {\it We have an isomorphism of $J(\Bfrm
)$-modules:}
$$
\VV^{\eta (\Bfrm )}\simeq \kmax\otimes \St (\Bfrm ,\rho )\ .
$$
\medskip

 Similarly, as a $U(\Bfr )J^1 (\Bfrm )$-module, $\VV^{\eta (\Bfr
 ,\Bfrm )}$ is a finite sum of irreducible submodules of the form
 $\kmax\otimes\tau$, where $\tau$ is an irreducible representation of
 $U(\Bfr )/U^1 (\Bfr )$. For such a $\tau$ we have:
\medskip

$
\matrix{
\scriptstyle{{\rm Hom}_{U(\Bfr )J^1 (\Bfrm )}\, (\kmax\otimes \tau ,\VV)}
& = & \scriptstyle{{\rm Hom}_{U(\Bfr )J^1 (\Bfrm )}(\tau ,{\rm Hom}_{U^1 (\Bfr )J^1
(\Bfrm )}\, (\kmax ,\VV ))}\cr
    &=& \scriptstyle{{\rm Hom}_{U(\Bfr )J^1 (\Bfrm )}(\tau ,\VV (\Bfrm )^{U^1 (\Bfr )J^1
(\Bfrm )})}\cr
 & = &   \scriptstyle{{\rm Hom}_{U(\Bfr )J^1 (\Bfrm )}(\tau ,\VV (\Bfrm
 )^{\Udss_\Bfr})}
}
$
\medskip

where $\VV (\Bfrm )^{\Udss_\Bfr}$ is the Jacquet module of $\VV (\Bfrm
)$ with respect to the unipotent radical $\Udss_\Bfr$ of the parabolic
subgroup $\Pdss_\Bfr$  of $\bar G$ given by $U(\Bfr )J^1 (\Bfrm
)/\allowbreak J^1 (\Bfrm )$.
Hence we have:
\medskip

$
\matrix{
{\rm Hom}_{U(\Bfr )J^1 (\Bfrm )}\, (\kmax\otimes \tau ,\VV)
& = &  {\rm Hom}_{U(\Bfr )J^1 (\Bfrm )}(\tau ,\VV (\Bfrm
 )^{\Udss_\Bfr})\cr
& = & {\rm Hom}_{\Pdss_\Bfr}\, (\tau , \VV (\Bfrm )^{\Udss_\Bfr})\cr
& = & {\rm Hom}_{\Ldss_\Bfr}\, (\tau ,\St (\Bfrm ,\rho )^{\Udss_\Bfr})
}
$
\medskip

 Denote by $\St (\Bfr ,\rho )$ the generalized Steinberg
 representation of $\Ldss_\Bfr$ with cuspidal support
 $(L_{\Bfr_0}, \rho )$.  It is classical that
$$
\St (\Bfrm ,\rho )^{\Udss_\Bfr}\simeq \St (\Bfr ,\rho )
$$
as $\Ldss_\Bfr$-modules. It follows that
$$
{\rm Dim}\,  {\rm Hom}_{U(\Bfr )J^1 (\Bfrm )}\, (\kmax\otimes \tau ,\VV)
=\left\{
\matrix{
0 & {\rm if}\ \tau \not\simeq \St (\Bfr ,\rho )\cr
1 & {\rm if} \ \tau \simeq \St (\Bfr ,\rho )
}
\right.
$$
As a consequence we have an isomorphism of $U(\Bfr )J^1 (\Bfrm
)$-modules:
$$
\VV^{\eta (\Bfr ,\Bfrm )} \simeq \kmax\otimes \St (\Bfr ,\rho )\ .
$$
\medskip

{\bf (XI.1.4) Proposition}. (i) {\it The  $U(\Afr )$-intertwing  of
$\kmax\otimes {\rm St}(\Bfr , \rho )$ is equal to $U(\Bfr )J^1 (\Bfrm )$.}

(ii) {\it The representation of $U(\Afr )$ given by
$$
\lambda (\Afr ):= {\rm Ind}_{U(\Bfr )J^1 (\Bfrm )}^{U(\Afr
)}\,  \kmax \otimes {\rm St}(\Bfr ,\rho )
$$
is irreducible.}

(iii) {\it We have
$$
\VV [\sigma_{\Bfr}]=\VV^{\lambda (\Afr )}\simeq \lambda (\Afr ),
$$
where the isomorphism is an isomorphism of $U(\Afr )$-modules.}
\medskip

{\it Proof}. The restriction of $\kmax\otimes \St (\Bfr ,\rho )$ to
$U^1 (\Bfr )J^1 (\Bfrm )$ is a multiple of $\eta (\Bfr ,\Bfrm )$, so
by Proposition (III.1.1)(v), we have
$$
I_G (\kmax\otimes \St (\Bfr , \rho ))\subset J^1 (\Bfrm )B^\times J^1
(\Bfrm )\ .
$$
In particular we have
$$
I_{U(\Afr )} (\kmax\otimes \St (\Bfr , \rho )) = U(\Bfr )J^1 (\Bfrm ) \ ,
$$
and point (i) follows.  Point (ii) is a consequence of Mackey
irreducibility criterion and point (iii) of  Lemma (V.2).

\bigskip

{\bf XI.2 Proof of conjecture (X.4.1) for irreducible discrete
 series representations.}

\bigskip

 Let $\Cfr$ be as before and $x$ be an element of $G$. Write $q={\rm
 Dim}\, \sigma_{\Cfr}$. Let $\Sigma$ be the $\Jm$-orbit $\Jm
 x\sigma_{\Cfr}$. We must prove that $K_\Sigma^{\lm} =0$.

 Recall that we have the exact sequence of $\Jm$-modules
$$
0 \lra K_\Sigma^{\lm} \lra C_q (\Sigma ,{\cal C})^{\lm} \lra \VV^{\lm}\lra
0\ ,
$$
if $x\in \Jm G_L U(\Afr )$, and
$$
0\lra  K_\Sigma^{\lm} \lra C_q (\Sigma ,{\cal C})^{\lm} \lra 0
\ ,
$$
if $x\not\in \Jm G_L U(\Afr )$.
\medskip

 Since $(\pi ,\VV )$ is a discrete series representation, $\lm$ occurs
 in $\VV$ with multiplicity $1$, so that $\VV^{\lm}\simeq \lm$ (see
 e.g. the discussion in [SZ] following the proof of Lemma 4, page 178).
 So we are reduced to proving the following
 result.
\medskip

{\bf (XI.2.1) Proposition}. {\it We have}
$$
{\rm Dim}\, {\rm Hom}_{\Jm}\, (\lm ,C_q (\Sigma , {\cal C}))\leq
\left\{
\matrix{
1 & {\rm  if }\  x\in \Jm G_L U(\Afr )\cr
0 & {\rm  otherwise}
}
\right.
$$

The rest of this section will be devoted to the proof of this
proposition.  Recall that
$$
C_q (\Sigma ,{\cal C}) =\coprod_{j\in \Jm /\Jm\cap U(\Afr )^x} jxC_q
(\sigma_\Cfr ,{\cal C}) = {\rm Ind}_{\Jm\cap U(\Afr )^x}^{\Jm}\,  xC_q
(\sigma_\Cfr ,{\cal C})\ .
$$
Using Proposition (XI.1.4), we obtain:
$$
\matrix{
xC_q (\sigma_\Cfr ,{\cal C}) & = &
{\rm Ind}_{\Jm\cap U(\Afr )^x}^{\Jm}\, x {\rm Ind}_{U(\Bfr )J^1 (\Bfrm
)}^{U(\Afr )}\kmax \otimes \St (\Bfr ,\rho )\cr
 &=& {\rm Ind}_{\Jm\cap U(\Afr )^x}^{\Jm}\, {\rm Ind}_{U(\Bfr )^x J^1
 (\Bfrm )^x}^{U(\Afr )^x}\, \kmax^x \otimes \St (\Bfr , \rho )^{x}\ .
}
$$
Mackey's restriction formula gives
$$
{\big(  {\rm Ind}_{U(\Bfr )^x J^1
 (\Bfrm )^x}^{U(\Afr )^x}\, \kmax^x \otimes \St (\Bfr , \rho )^{x}
\big)_{\vert J_{\max}\cap U(\Afr )^x}}
$$
$$
\matrix{
&=& {\bigoplus_{u\in U} {\rm Ind}_{\Jm \cap U(\Afr )^x \cap U(\Bfr
)^{ux} J^1 (\Bfrm )^{ux}}^{\Jm \cap U(\Afr )^x}\,
\kmax^{ux}\otimes \St (\Bfr ,\rho )^{ux}}\cr
& & \cr
&=&
{\bigoplus_{u\in U}  {\rm Ind}_{\Jm \cap U(\Bfr
)^{ux} J^1 (\Bfrm )^{ux}}^{\Jm \cap U(\Afr )^x}\,
\kmax^{ux}\otimes \St (\Bfr ,\rho )^{ux}}\ .
}
$$
where $U$ is the double coset set
$$
U= \Jm \cap U(\Afr )^x \backslash U(\Afr )^x /U(\Bfr )^x J^1 (\Bfrm
)^x\ .
$$
By Frobenius reciprocity we have:
$$
{\rm Hom}_{\Jm}\, (\lm ,C_q (\Sigma ,{\cal C}))
$$
$$
= \bigoplus_{u\in U} {\rm Hom}_{\Jm\cap U(\Bfr )^{ux}J^1 (\Bfrm
)^{ux}}\, (\lm , \kmax^{ux}\otimes \St (\Bfr ,\rho )^{ux})
$$
By definition of the cuspidal support of a representation of $U(\Bfr
)/U^1 (\Bfr )$, we have that $\kmax\otimes \St (\Bfr ,\rho )$ embeds
in
$$
{\rm Ind}_{U(\Bfr_0 )J^1 (\Bfrm )}^{U(\Bfr )J^1 (\Bfrm
)}\, \kmax\otimes \rho ={\rm Ind}_{\Jm}^{U(\Bfr )J^1 (\Bfrm
)}\, \lm
$$
as a $U(\Bfr )J^1 (\Bfrm )$-module. It follows that ${\rm
Hom}_{\Jm} \, (\lm ,C_q (\Sigma ,{\cal C}))$ embeds in the
$\Cdss$-vector space
$$
\bigoplus_{u\in U}{\rm Hom}_{\Jm\cap U(\Bfr )^{ux}J^1 (\Bfrm )^{ux}} \,
(\lm , {\rm Ind}_{\Jm^{ux}}^{U(\Bfr )^{ux}J^1 (\Bfr )^{ux}}\,  \lm^{ux})
$$
Using Mackey's restriction formula again, we obtain:
$$
\big(  {\rm Ind}_{\Jm^{ux}}^{U(\Bfr )^{ux}J^1 (\Bfr
)^{ux}}\,  \lm^{ux})\big)_{\vert \Jm \cap U(\Bfr )^{ux}J^1 (\Bfrm
)^{ux}}
$$
$$
=\bigoplus_{v\in V_u} {\rm Ind}_{\Jm\cap U(\Bfr )^{ux}J^1 (\Bfrm
)^{ux}\cap \Jm^{vux}}^{\Jm\cap U(\bfr )^{ux}J^1 (\Bfrm
)^{ux}}\, \lm^{vux}
$$
$$
=\bigoplus_{v\in V_u} {\rm Ind}_{\Jm \cap \Jm^{vux}}^{\Jm\cap U(\Bfr
)^{ux} J^1 (\Bfrm
)^{ux}}\, \lm^{vux}\ ,
$$
where
$$
V_u = \Jm\cap U(\Bfr )^{ux}J^1 (\Bfrm )^{ux}\backslash U(\Bfr )^{ux}
J^1 (\Bfrm )^{ux} / \Jm^{ux}\ .
$$

Hence it follows by Frobenius reciprocity that ${\rm Hom}_{\Jm}\, (\lm
,C_q (\Sigma ,{\cal C}))$ embeds in
$$
\bigoplus_{u\in U}\bigoplus_{v\in V_v} {\rm Hom}_{\Jm\cap \Jm^{vux}}\,
(\lm ,\lm^{vux})\ .
$$
asz a $\Cdss$-vector space. As a consequence, if  ${\rm Hom}_{\Jm}\, (\lm
,C_q (\Sigma ,{\cal C}))$ is non-zero, there exist $u\in U$, $v\in V_u$ such that
$vux$ intertwines $\Jm$, that is
$$
vux\in \Jm G_L \Jm\ .
$$

For such $u$ and $v$, we have $u\in xU(\Afr )x^{-1}$ and
$$
v\in ux U(\Bfr )J^1 (\Bfrm )x^{-1}u^{-1}
$$
so that
$$
vux\in uxU(\Bfr )J^1 (\Bfrm )\subset xU(\Afr )U(\Bfr ) J^1 (\Bfrm
)=xU(\Afr )\ .
$$
Hence we have $xU(\Afr )\cap \Jm G_L \Jm\not= \emptyset$, that is
$x\in \Jm G_L U(\Afr )$. As a consequence Proposition (XI.2.1) holds
when $x\not\in \Jm G_L U(\Afr )$.
\medskip

 Now let us assume that $x\in \Jm G_L U(\Afr )$. Writing $x=jx_L u$,
 $j\in \Jm$, $x_L \in G_L$ and $u\in U(\Afr )$, we have that
$$
\Sigma =\Jm x\sigma_\Cfr =\Jm x_L \sigma_\Cfr
$$
so that we may as well assume that $x\in G_L$.
\medskip

 Assume that for some $u\in xU(\Afr )x^{-1}$, we have
$$
{\rm Hom}_{\Jm\cap U(\Bfr )^{ux} J^{1}(\Bfrm )^{ux}}\, (\lm ,\kmax^{ux}\otimes \St (\Bfr ,\rho )^{ux})\not= 0\ .
$$

Then by the preceeding discussion, there exists $u\in ux U(\Bfr )J¹ (\Bfrm ) (ux)^{-1} $ such that
$$
vux \in ux U(\Bfr )J^1 (\Bfrm )\cap \Jm G_L \Jm\ .
$$
This implies that
$$
uxU(\Bfr )J^1 (\Bfrm )x^{-1} \cap \Jm G_L \Jm^{x} \not=
\emptyset
$$
that is $u\in \Jm G_L U(\Bfr )^x J^1 (\Bfrm )^x$. So without changing
the double class $\bar u$ of $u$ in $U$, we may as well assume
that $u\in \Jm G_L $. Let us write $u=jg_L$, $j\in \Jm$, $g_L \in G_L$. Since $u\in
xU(\Afr )x^{-1}$, we have
$$
u(x\sigma_\Cfr )=x\sigma_\Cfr =j(g_L x\sigma_\Cfr )\ .
$$
So $x\sigma_\Cfr$ and $g_L x\sigma_\Cfr$ are simplices of $X_L$ conjugated under the action
of $\Jm\subset U(\Afr_{\rm min})$. By Lemma (X.4.4), there exists $i\in U(\Cfr_{\rm min})$ such that
$x\sigma_\Cfr =ig_L x\sigma_\Cfr$. Hence $ig_L\in U(\Afr )^x \cap G_L =U(\Cfr )^x$ and as a consequence
$g_L \in U(\Cfr_{\rm min})U(\Cfr )^x$. It follows that $u\in \Jm U(\Cfr_{\rm min})U(\Cfr )^x =\Jm U(\Cfr )^x$, and
$u\in (\Jm\cap U(\Afr )^x ).\allowbreak U(\Cfr )$. But this implies that the image
$\bar u$ of $u$ in
$$
U=\Jm\cap U(\Afr )^x \backslash U(\Afr )^x /U(\Bfr )^x J^1 (\Bfrm
)^x
$$
 is $\bar 1$. We have proved  the following:
\medskip

{\bf  (XI.2.2) Lemma}. {\it For all $x\in G_L$, we have}
$$
{\rm Hom}_{\Jm}\, (\lm ,C_q (\Sigma ,{\cal C}))
={\rm Hom}_{\Jm \cap U(\Bfr )^x J^1 (\Bfrm )^x}\, (\lm ,\kmax^x \otimes \St (\Bfr ,\rho^x )
$$
$$
={\rm Hom}_{U(\Bfr_0 )J^1 (\Bfrm )\cap U (\Bfr )^x J^1 (\Bfrm )^x} \, (\kmax\otimes \rho ,
\kmax^x \otimes\St (\Bfr ,\rho )^x )\ .
$$

We next prove:

{\bf  (XI.2.3) Lemma}. {\it For all $x\in G_L$, we have}

$$
{\rm Dim}\, {\rm Hom}_{U(\Bfr_0 )J^1 (\Bfrm )\cap U (\Bfr )^x J^1 (\Bfrm )^x} \, (\kmax\otimes \rho ,
\kmax^x \otimes\St (\Bfr ,\rho )^x )
$$
$$
={\rm Dim}\, {\rm Hom}_{U(\Bfr_0 )\cap U(\Bfr )^x}\, (\rho , \St (\Bfr ,\rho )^x )
\ .
$$

{\it Proof}. It is inspired from that of [BK](5.3.2), page 176. Abreviate
$\rho_\Bfr =\St (\rho ,\Bfr )$. Write $Y$ (resp. $X_0$, $X$)
for the space of $\kmax$ (resp. $\rho$, $\rho_\Bfr$). Let
$\varphi\in {\rm Hom}_\Cdss \, (Y\otimes X_0 , Y\otimes X)$ $=$ ${\rm End}_\Cdss (Y)\otimes
{\rm Hom}_\Cdss (X_0 ,X)$ and write
$$
\varphi =\sum_{i\in I}  S_i \otimes T_i
$$
where $S_i \in {\rm End}_\Cdss \, (Y)$, $T_i \in {\rm Hom}_\Cdss \, (X_0 ,X)$, and where the $T_i$
are linearly independent. Then $\varphi$ intertwine $\kmax\otimes\rho$
and $\kmax^x \otimes \rho_\Bfr^x$ if and only if
$$
\sum_{i\in}\, (S_i \circ \kmax (u))\otimes (T_i \circ \rho (u))
=  \sum_{i\in I} (\kmax^x (u)\circ S_i )\otimes (\rho_\Bfr^x (u)\circ T_i )
$$
for all $u\in U(\Bfr_0 ) J^1 (\Bfrm )\cap U(\Bfr )^x J^1 (\Bfrm )^x$. In particular if
 $\varphi$ intertwines these representations, for
$u\in  J^1 (\Bfrm )\cap J^1 (\Bfrm )^x$, we must have
$$
\sum_{i\in I} (S_i \circ \kmax (u)-\kmax^x (u)\circ S_i )\otimes T_i =0
$$
Since the $T_i$ are linearly independent, we obtain
$$
S_i \in {\rm Hom}_{J^1 (\Bfrm )\cap J^1 (\Bfrm )^x}\, (\kmax ,\kmax^x )
$$
$$
={\rm Hom}_{J^1 (\Bfrm )\cap J^1 (\Bfrm )^x}\, (\eta_{\rm max} ,\eta_{\rm max}^x )
\ .
$$
By [BK](5.1.8) and (5.2.7), the spaces
$$
{\rm Hom}_{J (\Bfrm )\cap J(\Bfrm )^x}\, (\kmax ,\kmax^x )
$$
and
$$
{\rm Hom}_{J^1 (\Bfrm )\cap J^1 (\Bfrm )^x}\, (\eta_{\rm max} ,\eta_{\rm max}^x )
$$
are equal and $1$-dimensional. It follows that any $\varphi$ in
$$
{\rm Hom}_{\Jm \cap U(\Bfr )^x J^1 (\Bfr )^x }\, (\lm
,\kmax^x \otimes \rho_\Bfr^x )
$$
writes $\varphi =S\otimes T$, where $S\in {\rm Hom}_{J(\Bfrm )\cap J(\Bfrm )^x} ,
(\kmax ,\kmax^x )$ and $T\in {\rm Hom}_\Cdss\, (X_0 ,X)$. Writing that such
a $S\otimes T$ does intertwine the representations, we easily obtain that
$$
T\in {\rm Hom}_{U(\Bfr_0 )\cap U(\Bfr )^x}\, (\rho ,\rho_{\Bfr})\ .
$$
It follows that we have a canonical isomorphism of $\Cdss$-vector spaces:
$$
{\rm Hom}_{\Jm \cap U(\Bfr )^x J^1 (\Bfr )^x }\, (\lm ,\kmax^x \otimes \rho_\Bfr^x )
$$
$$
={\rm Hom}_{J(\Bfrm )\cap J(\Bfrm )^x}\, (\kmax ,\kmax^x )\otimes {\rm Hom}_{U(\Bfr_0 )\cap
U(\Bfr )^x}\, (\rho ,\rho_\Bfr^x )\ ,
$$
with
$$
{\rm Dim}\, {\rm Hom}_{J(\Bfrm )\cap J(\Bfrm )^x}\, (\kmax ,\kmax^x ) =1
$$
and the lemma follows.

To obtain Proposition (XI.2.1), we are now reduced to proving the following result.
\medskip

{\bf  (XI.2.3) Lemma}. {\it For all $x\in G_L$, we have}
$$
{\rm Dim}_{\Cdss}\, {\rm Hom}_{U(\Bfr_0 )\cap U(\Bfr )^x}\, (\rho ,
\St (\Bfr ,\rho )^{x})\leq 1\ .
$$

\bigskip

 Fix a level $0$ discrete series representation $(\pi_0 ,\VV_0 )$
 of $G_E$ belonging to the
Bernstein   component of $G_E$ defined by the type $(U(\Bfr_0 ), \rho )$. Applying the results
 of section (XI.1) to $(\pi_0 ,\VV_0 )$, we have that $\rho \simeq \VV_0^{U^1 (\Bfr_0 )}$ as
 $U(\Bfr_0 )$-modules and $\St (\Bfr ,\rho )^x \simeq \VV_0^{U^1 (\Bfr )^x}$ as
$U(\Bfr )^x$-modules.  Hence the statment of the lemma rewrites:

{\it Let $\Bfr_1$ be a hereditary order lying in the image of the canonical map
$ {\rm Her}(C)\lra {\rm Her}(B)$. Then
$$
{\rm Dim}_{\Cdss} \, {\rm Hom}_{U(\Bfr_0 )\cap U(\Bfr_1 )}\, (\VV_0^{U^1 (\Bfr_0 )}, \VV_0^{U^1 (\Bfr_1 )})
\leq 1\ .
$$}

 We may write this in the language of simplicial complexes. For $\sigma$ a simplex of $X_E$,
 write $U_\sigma$ for the parahoric subgroup of $G_E$ fixing $\sigma$ and $U_\sigma^1$ for its
 pro-unipotent radical. Write $\sigma_0 =\sigma_{\Bfr_0}$. Then our lemma is equivalent to:
\medskip

{\bf  (XI.2.4) Lemma}. {\it For all simplex $\tau$ lying in the image of the canonical
simplicial map $X_L\lra X_E$, we have
$$
{\rm Dim}_{\Cdss}\, {\rm Hom}_{U_{\sigma_0}\cap U_\tau }(\VV_0^{U^1_{\sigma_0}},\VV_0^{U^1_\tau})\leq 1\ .
$$}

Fix an apartment $\AL$ of $X_L$ containing $\sigma_0$ and $\tau$ (we see $X_L\lra X_E$ as an
inclusion). According to [Br] Lemma 4, there exists a unique chamber $\sigma$ of $\AL$ such
that we have the containments
$$
E[\sigma_0 ,\tau ]\supset \sigma \supset \tau
$$
where $E[\sigma ,\tau ]$ is the {\it enclos} of $\sigma\cup\tau$ in the sense of
[BTI] Definition (2.4.1). Moreover by [Br] Lemma 5, we have that the simplex $\sigma$ lies
between $\sigma_0$ and $\tau$ in the sense of [SS] {\S}2. This means that there exists points
$x_{\sigma_0}$ in $\vert \sigma_0\vert^{\rm o}$, $x_\sigma$ in $\vert\sigma\vert^{\rm o}$, and $x_\tau$
 in $\vert\tau\vert^{\rm o}$ such that $x_\sigma$ belongs to the geometric segment
$[x_{\sigma_0},x_\tau ]$. Since the embedding $X_L\lra X_E$ is simplicial and affine, we have that,
 as a simplex of $X_E$, $\sigma$ lies between $\sigma_0$ and $\tau$.   We may then apply
Proposition (2.5) of [SS]:
\medskip

{\bf  (XI.2.5) Lemma}. {\it The image of $U^1_\sigma = U^1_{\sigma}\cap U_\tau$ in $U_\tau /U^1_\tau$ is
contained in the image of $U^1_{\sigma_0}  \cap U_\tau$ in $U_\tau /U^1_\tau$.}
\medskip

 Next fix an $L$-basis $(v_1 , ..., v_e )$ of $V$ corresponding to the apartment $\AL$. Moreover
fix a basis $(\zeta_1 ,...,\zeta_r )$ of the $\ofr_E$-module $\ofr_L$. Set
$$
V_i =Lv_i = {\rm Vect}_{E} \langle  \zeta_j v_i \ ; \ j=1,...,r\rangle\ , \ i=1,...e\ ,
$$
and write $M$ for the Levi subgroup of $G_E$ corresponding to the decomposition
$V=V_1 \oplus \cdots \oplus V_e$.
\medskip

{\bf (XI.2.6) Lemma}. {\it Let $\theta$ be a simplex of $\AL$ (seen as a simplex of $X_E$).}

(i) {\it The intersection $U_\theta \cap M$ does not depends on $\theta$. We denote it
 by $M^0$. It is given by
$$
\prod_{i=1,..,e}{\rm GL}(r,\ofr_E )
$$
where the $i$th copy of ${\rm GL}(r,\ofr_E )$ is the maximal compact subgroup of
${\rm Aut}_E \, (V_i )$ which is standard in the basis $(\zeta_j v_i )_j$ of $V_i$.}

(ii) {\it Assume moreover that $\theta$ is a chamber of $\AL$. Then we have the Iwahori
 decomposition
$$
U_\theta =(U_\theta \cap M)U_{\theta}^1\ .
$$}

{\it Proof}. This is an easy exercice in lattice chain theory and we
only sketch the proofs.

 The simplex $\theta$ corresponds to a certain $\ofr_L$-lattice chain
 ${\cal N}=(N_k )_{k\in \Zdss}$ in $V$. The fact that $\theta$ lies in
 ${\cal A}_L$ exactly means that the chain $\cal L$ is split by the
 decomposition $V=\bigoplus V_i$, i.e. for $k$ in $\Zdss$ we have:
$$
N_k =\bigoplus N_k^i ,  \ N_k^i =N_k \cap V_i , \ i=1,...e.
$$

Let $g\in U_\theta \cap M$ that we write $g=\bigoplus g_i$, $g_i \in
{\rm End}_E \, V_i$, $i=1,...,e$. Then we get $g_i N_k^i =N_k^i$,
$k\in \Zdss$, $i=1,...,e$, that is $g_i \ofr_L v_i =\ofr_L v_i$,
$i=1,...,e$, and point (i) follows easily.

 Assume moreover that $\theta$ is a chamber in ${\cal A}_L$. Since the
 identity of (ii) is invariant under the action of the affine Weyl
 group of this apartment (since it stabilizes $M$), we may as well
 assume that $\theta$ is the
 standard chamber attached to the lattice chain $\cal N$ defined by
$N_k =\bigoplus N_k^i$ as above and, for $k=0,...,e-1$,
$N_k^i = \ofr_L v_i$, if $i\in \{ 0,...,e-1-k\}$, $N_k^i =\pfr_L v_i$,
 $i\in \{ e-1-k+1 ,...,e-1\}$. Then by a straightforward computation,
 we obtain that an element $g\in {\rm End}_E\, V$, with a block matrix
$g=(g_{uv})_{u,v=1,...,e}$ in the decomposition $V=\bigoplus V_i$, lies
 in $U_\theta$ if and only if we have $g_{uu}\in {\rm GL}(r,\ofr_E )$,
 $u=1,...,e$, $g_{uv}\in {\rm M}(r,\ofr_E )$, if $v>u$, and $g_{uv}\in
 \pfr_E {\rm M}(r,\ofr_E )$, if $u>v$. It is then classical that such
 a matrix has an Iwahori decomposition as described by the identity of
 (ii).

 We have $M^0 \subset U_{\sigma_0}\cap U_\tau$ and, as a $M^0$-module,
 $\VV_0^{U^1_{\sigma_0}}$ is isomorphic to the irreducible
 representation  $\rho^{\otimes e}$. So in order to prove
 Lemma (XI.2.3), it suffices, by Schur Lemma, to show that
 ${\rm Im}\, \varphi =\varphi (\VV_0^{U^1_{\sigma_0}})$ is independent of the
choice of a non-zero intertwining operator $\varphi$ in
${\rm Hom}_{U_{\sigma_0}\cap U_\tau}\, (\VV_0^{U^1_{\sigma_0}}, \VV_0^{U^1_\tau})$.

 Let $\varphi$ such a non-zero  intertwining operator.
Then ${\cal W}:=\varphi (\VV_0^{U^1_{\sigma_0}})$
is a sub-$M^0$-module of $\VV_0^{U^1_\tau}$ equivalent to $\rho^{\otimes e}$. The groups
$U^1_{\sigma_0}\cap U_\tau$ and $U^1_{\tau}$ act trivially on $\cal W$. Moreover by Lemma
(XI.2.5), $U^1_\sigma \subset (U^1_{\sigma_0}\cap U_\tau )U^1_\tau$. It follows that $U^1_\sigma$ acts
trivially on $\cal W$ and that the action of $U_\sigma /U^1_\tau$ on that space is the action
 inflated from the representation $\rho^{\otimes e}$ of $U_\tau /U^1_\tau$; write
$\rho_\sigma$ for the corresponding representation.

 As a $U_\sigma$-module, $\VV_0^{U^1_\tau}$ is equivalent to the generalized Steinberg representation
with cuspidal support $({\bar M}_0 ,\rho^{\otimes e})$, where ${\bar M}_0$ is the image of $M_0$
in the quotient $U_\sigma /U^1_\sigma$. It is well known that this Steinberg representation
occurs with multiplicity $1$ in the parabolically induced representation ${\rm Ind}_{U_\sigma}^{U_\tau}
\, \rho_\sigma$. It follows by Frobenius reciprocity that $\rho_\sigma$ occurs in $\VV_0^{U^1_\tau}$ with
multiplicity $1$. It follows that $\cal W$ is the unique sub-$U_\sigma$-module of $\VV_0^{U^1_\tau}$
isomorphic to $\rho_\sigma$ and Lemma (XI.2.4) is proved.

 We have unconditionnaly proved the following result.

{\bf (XI.2.7) Theorem}. {\it Assume that $(\pi ,\VV )$ is an unramified twist of an irreducible
unitary  discrete series representation lying in the Bernstein block $\RR_{(J,\lambda )}$.
Then the augmented chain complex (IX.1) is exact.}

\bigskip

\centerline{\bf XII. Explicit pseudo-coefficients for discrete series
representations.}
\bigskip

Let $(\pi ,\VV )$ be an irreducible (unitary) discrete series
representation of $G$. In this section, following [SS2]{\S}II.4, we
show that Theorem (XI.2.7)   leads to an explicit pseudo-coefficient
$\varphi_\pi$ for $\pi$.  We then show how to derive an explicit
formula for the value of the Harish-Chandra character of $\pi$ at an
elliptic regular element.
\medskip

{\bf XII.1 The coefficient system ${\cal C}(\pi )$.}
\medskip

 Recall that, with the notation of {\S}XI, the coefficient system
 ${\cal C}={\cal C}(\pi )$ canonically attached to $\pi$ is given on
a part of $X_L$ by $V[\sigma_\Cfr ] =\VV^{\lambda (\Afr )} \simeq \lambda (\Afr
 )$, where
\medskip

 -- $\sigma =\sigma_\Cfr$ is any simplex of $X_L$ satisfying $\Cfr_{\rm
    min}\subset \Cfr \subset \Cfr_{\rm max}$,

 -- $\Afr =\Afr (\Cfr )$,

 -- $\lambda (\Afr )={\rm Ind}_{U(\Bfr )J^1 (\Bfrm )}^{U(\Afr )} \kappa_{\rm
    max}\otimes \St (\Bfr ,\rho )$.
\medskip

Since the coefficient system $\cal C$ is $G$-equivariant, for any
order $\Afr$ as above, the representation $\lambda (\Afr )$ extends to
a representation of ${\cal K}(\Afr )={\cal N}_G (\sigma )$ that we
still denote by $\lambda (\Afr )$. In the sequel we shall also write
$$
{\cal N}_G (\sigma )={\cal K}_\sigma\ {\rm and}\ \lambda (\Afr
)=\lambda_\sigma
\ .
$$

By equivariance, we may define an irreducible smooth representation
$\lambda_\sigma$ of ${\cal K}_\sigma$ for {\it any} simplex $\sigma$
of $X[L]$, and by equivariance of $\cal C$ we have $\VV [\sigma
]=\VV^{\lambda_\sigma}\simeq \lambda_\sigma$.
\medskip

{\bf XII.2 Euler-Poincar\'e functions.}
\medskip

Let $\chi$ be the central character of $\pi$. All representations that  we
consider will lie in the category ${\cal S}_\chi (G)$ of those smooth
representations admitting a central character equal to $\chi$.  If
${\cal V}', {\cal V}''\in {\cal S}_\chi (G)$ with ${\cal V}'$  of finite length and ${\cal V}''$
admissible, we define the Euler-Poincar\'e characteristic :
$$
{\rm EP}\, ({\cal V}',{\cal V}'')=\sum_{q\geq 0} (-1)^q\, {\rm dim}\, {\rm
Ext}^q_{{\cal S}_\chi (G)}\, ({\cal V}' ,{\cal V}'')\ .
$$

We denote by $Z$ the center of $G$ and fix a Haar measure $\mu_{G/Z}$
on $G/Z$. We denote by $\HH_\chi (G)$ the convolution Hecke algebra of
locally constant functions $f$~: $G\lra \Cdss$ satisfying
\medskip

 -- $f(zf)=\chi^{-1}(z)f(g)$, $z\in Z$, $g\in G$,
\smallskip

 -- $f$ has compact support modulo $Z$.
\medskip

Representations in ${\cal S}_\chi (G)$ are naturally left $\HH_\chi
(G)$-modules.  The {\it character} of an admissible representation
$(\pi ',{\cal V}')$ in ${\cal S}_\chi (G)$ is the functionnal
$$
{\rm Tr}_{{\cal V}'}~: \ \HH_\chi (G)\lra \Cdss\ , \ \psi \lra {\rm tr}\,
(\pi '(\psi))\ ,
$$
where $\pi '(\psi )$ is the endomorphism of ${\cal V}'$ attached to $\pi'$;
it is formally given by the integral
$$
\pi '(\psi )=\int_{G/Z}\psi (g)\pi '(g)\, d\mu_{G/Z}\, (\dot{g})\ .
$$

We set $d={\rm dim}\, X[L]$. For $q=0, ...,d$, we fix a set $\FF_q$ of
representatives of $G$-orbits in the set $X[L]_q$ of $q$-simplices in
$X[L]$. If $\sigma =\sigma_\Cfr$ is a simplex of  $X[L]$, we denote by
$\epsilon_\sigma$~: $\KK_\sigma\lra \{\pm 1\}$ the abelian character
defined as follows. If $g\in \KK_\sigma$, $\epsilon_\sigma (g)$ if the
sign of the permutation of the vertex set of $\sigma$ induced by the
action of $g$. Moreover for such a simplex $\sigma$, we denote by
$\tau_\sigma^{\VV}$ the character of the representation $(\KK_\sigma
,\lambda_\sigma )$. For all simplices of $X[L]$ we extend the class
functions $\epsilon_\sigma$ and $\tau_\sigma^\VV$ by zero to functions
on $G$. Following Kottwitz [Kot] and Schneider and Stuhler [SS2], we
define the {\it Euler-Poincar\'e function} attached to $(\pi ,\VV )$
by the formula:
$$
f_{\rm EP}^\VV := \sum_{q=0}^{d} \sum_{\sigma \in {\cal F}_q}
(-1)^q .\,   \mu_{G/Z} (\KK_\sigma /Z )^{-1} . {\bar \tau}_\sigma^{\VV}
. \epsilon_\sigma\ .
$$

{\bf Remark}. The Euler-Poincar\'e function $f_{\rm EP}^\VV$ does depend
on the choices of representative sets ${\cal F}_q$, $q=0, ...,d$.
\bigskip

{\bf (XII.2.1) Proposition}. {\it For all admissible representations
$(\pi ', \VV ')$ in ${\cal S}_\chi (G)$, we have
$$
{\rm Tr}_{\VV '}\, (f_{\rm EP}^\VV )={\rm EP}_{{\cal S}_\chi (G)}\,
(\VV ,\VV ')\ .
$$}
\medskip

{\it Proof}. We have the decomposition.
$$
C_c^{\rm or} (X_{(q)}, \CC (\VV ))=\bigoplus_{\sigma \in \FF_q}
C_c^{\rm or} (G.\sigma ,\CC (\VV ))\ ,
$$
where $C_c^{\rm or} (G.\sigma ,\CC (\VV ))$ denotes the $G$-space of
oriented chains with support in $G.\{ (\sigma ,o_1 ), (\sigma ,o_2 )\}$,
$o_1$ and $o_2$ denoting the two possible orientations of
$\sigma$. Since
$$
C_c^{\rm or}(X_{(\bullet )},\CC (\VV ))\augmentation \VV
$$
is a projective resolution of $\VV$ in $\SS (G)_\chi$, Lefschetz
formula gives :
$$
{\rm EP}_{\SS (G)_\chi}\, (\VV ,\VV ')=\sum_{q=0}^d
(-1)^q \sum_{\sigma \in \FF_q} \, {\rm
dim}\, {\rm Hom}_G \, (C_c^{\rm or} (G.\sigma ,\CC (\VV )), \VV ')\ .
$$
By definition of compact induction we have:
$$
C_c^{\rm or}(G.\sigma ,\CC (\VV ))={\rm c-Ind}_{\KK_\sigma}^G\, C_c^{\rm
or}(\sigma ,\CC (\VV ))\ ,
$$
where $\CC_c^{\rm or}(\sigma ,\CC (\VV ))$ denotes the
$\KK_\sigma$-space of chains withs support in\hfill\break
$\{ (\sigma ,o_1 ),
(\sigma ,o_2 )\}$. Moreoer, again by definition, we have the
isomorphism  of $\KK_\sigma$-modules:
$$
C_c^{\rm or}(\sigma ,\CC (\VV ))
= \lambda_\sigma \otimes \epsilon_\sigma\ .
$$
Using Frobenius reciprocity for compact induction, we obtain
$$
{\rm Hom}_G\, ( C_c^{\rm or}(\sigma ,\CC (\VV )),\VV ')={\rm
Hom}_{\KK_\sigma}(\lambda_\sigma\otimes \epsilon_\sigma ,\VV ')\ .
$$
Moreover ${\rm dim}\, {\rm Hom}_{\KK_\sigma} \,
(\lambda_\sigma\otimes \epsilon_\sigma ,\VV ')$ is nothing other than
the multiplicity of $\lambda_\sigma \otimes \epsilon_\sigma$ in the
isotypic component $(\VV ')^{\lambda_\sigma \otimes\epsilon_\sigma}$:
$$
{\rm dim}\, {\rm Hom}_{\KK_\sigma}\, (\lambda_\sigma\otimes\epsilon_\sigma ,\VV
')= {1\over{{\rm dim}\, \lambda_\sigma}}\,  {\rm dim}\, (\VV
')^{\lambda_\sigma\otimes \epsilon_\sigma}\ .
$$
Hence we have obtained
$$
{\rm EP}_{\SS (G)_\chi} (\VV ,\VV
')=\sum_{q=0}^{d}\sum_{\sigma\in \FF_q}{(-1)^q\over{\rm
dim}\, \lambda_\sigma}\, {\rm dim}\, (\VV
')^{\lambda_\sigma\otimes \epsilon_\sigma}\ .
$$

We need to compare this with ${\rm Tr}_{\VV '}\, (f_{\rm
EP}^{\VV})$. For this we have to compute\hfill\break
 ${\rm tr}_{\VV '}({\bar \tau}_\sigma^\VV .\epsilon_\sigma )$, for all $q$ and
$\sigma\in \FF_q$. Recall that for such $q$ and $\sigma$,
$$
E_\sigma := {1\over \mu_{G/Z}(\KK_\sigma
)}\, {\bar \tau}_\sigma^\VV \epsilon_\sigma . {\rm dim}\,
(\lambda_\sigma \epsilon_\sigma )
$$
is an idempotent of $\HH (G)_\chi$, and that $E_\sigma$ seen as an
endomorphism of $\VV '$ is the projection of the
$\lambda_\sigma \otimes \epsilon_\sigma$-isotypic component $(\VV
')^{\lambda_\sigma \otimes\epsilon_\sigma}$. Hence we have  that ${\rm
Tr}_{\VV '} (E_\sigma )={\rm dim}\, (\VV
')^{\lambda_\sigma \otimes\epsilon_\sigma}$ and the proposition
follows.

 We shall need the following result.
\bigskip

{\bf (XII.2.2) Theorem}. {\it Let $\VV '$ be an irreducible tempered representation
in $\SS_\chi (G)$. Then :}
$$
{\rm EP}_{\SS_\chi (G)}(\VV ,\VV ')=\cases{
$1$ & {\rm if} \  $\VV '\simeq \VV$ ,\cr
$0$ &  {\rm otherwise.}\cr
}
$$

{\it Proof}.  It is shown in [SZ2] Prop. 9.3 and subsequent remark (based upon a result 
of R. Meyer in [Me]) that
$$
{\rm
Ext}^\ast_{{\cal S}_\chi (G)}\, ({\cal V} ,{\cal V}') = {\rm
Ext}^\ast_{{\cal S}_\chi^{temp} (G)}\, ({\cal V} ,{\cal V}') \ ,
$$
where ${\cal S}_\chi^{temp} (G)$ denotes the category of all tempered smooth representations with central 
character $\chi$. But by a variant of [SZ2] Prop. 2.3 the representation $\VV$ is a projective object in
${\cal S}_\chi^{temp} (G)$.

\bigskip

 Recall that a function $f\in \HH (G)_\chi$ is a pseudo-coefficient of
 $(\pi ,V)$ if for any irreducible tempered representation in
 $\SS_\chi (G)$, we have
$$
{\rm Tr}_{\VV '}(f)=\cases{
$1$ & if $\VV '\simeq \VV$ ,\cr
$0$ & otherwise.\cr }
$$

As a consequence of (XII.2.1) and (XII.2.2) we have :
\bigskip

{\bf (XII.2.3) Theorem.} {\it The Euler-Poincar\'e function $f_{\rm
EP}^{\VV}$  is a pseudo-coefficient of  $(\pi ,\VV )$.}
\medskip

 In [Br2], the first author obtained pseudo-coefficients for discrete
 series representations of $G$ using a quite different approach (but
 also based on Bushnell and Kutzko type theory). Our
 pseudo-coefficients are likely to be very close to those of [Br2], but
 the comparison has yet to be done.

\bigskip

{\bf XII.3. An explicit character formula.}
\bigskip

 If $\psi\in \HH (G)_\chi$ and $h\in G$ is a regular elliptic
 element, the orbital integral
$$
\psiv (h):=\int_{G/Z}\psi (g^{-1}hg)d\mu_{G/Z}(\dot{g})
$$
is known to converge (see e.g. [SS2], page 140 in the case of a
reductive group with compact center, the non-compact case being
similar).

 Let $\Theta_\pi$ denote the Harish-Chandra character of $(\pi ,\VV
 )$. This is  a locally constant function on the set $G^{\rm reg}$ of
 regular semisimple elements of $G$. The following result relates
 values of $\Theta_\pi$ with the orbital integral of a
 pseudo-coefficient of $\pi$.

{\bf (XII.3.1) Theorem}. (Kazhdan-Badulescu) {\it Let $f_0$ be a
pseudo-coefficient of $(\pi ,\VV )$. Then for all regular elliptic
element $h$ of $G$, we have
$$
\Theta_\pi (h) = \fv_0 (h^{-1})\ .
$$}

{\it Remark}. This theorem is due to Kazhdan ([Ka], Prop. 3, page 28)
for a reductive group with compact center when $F$ has characteristic
$0$. It is due
to Badulescu ([Ba] Th\'eor\`eme (4.3)(ii), page 64) for our group $G$ without
restriction on $F$.
\bigskip

 Let $h\in G$ be a regular elliptic element. To obtain  a
 formula for $\Theta_\pi (h)$ it suffices to compute $(f_{\rm EP}^\VV )^{\vee}
 (h^{-1})$ explicitely. For this we closely follow the proof of Lemma
 (III.4.10) of [SS2] where a similar computation is done.

 If $\vert X\vert$ denotes the geometric realization of the building
 of $G$, it is known that $\vert X\vert^h$ is compact (see e.g. [SS2],
 page 141). Hence so is $\vert X[L]\vert^h$ the set of $h$-fixed
 points in the geometric realization of $X[L]$ since the subset $X[L]\subset X$ is
closed. Let us sketch the proof of this latter fact. Let $x_n =g_n .c_n$ be a converging sequence
of points in $X[L]$ with limit $x$  where,  for all $n$, $g_n$ is in $G$ and $c_n$ lies
in some fixed (closed) chambre $C_L$
of $X_L$. Then $(c_n )$ has a convergent subsequence and replacing $(x_n )$ by a subsequence
 we may assume that $c_n$ converges to some $c\in C_L$. Let $d$ be a $G$-invariant metric on $G$. we have
$d(x, g_n .c_n )=d(g_n^{-1}.x, c_n )\lra 0$. Hence $d(c_n ,G.x )\lra 0$ and $d(c,G.x)=0$. But it is an
easy exercice  in Bruhat-Tits theory (left to the reader) that the $G$-orbit of any point of $X$ is closed
in $X$. Hence $c\in G.x$, that is $x\in G.c\subset X[L]$ as required.

It follows that their
 exists a finite number of simplices $\sigma$ in $X[L]$ such that
 $h.\sigma =\sigma$. For such a $\sigma$, the intersection
 $\sigma\cap \vert X[L]\vert^h$ is non-empty. The collection
of $\sigma (h)$ where $\sigma$ runs over the simplices of $X[L]$
 globally fixed by $h$ endows the compact topological set
$\vert X[L]\vert^h$ with a simplicial
 structure. As noticed by Kottwitz ([Kot], page 635), it is an easy
 exercise to check that for all $\sigma$ in $X[L]$ fixed by $h$ we
 have
$$
\epsilon_\sigma (h)=(-1)^{{\rm dim}\, \sigma -{\rm dim}\, \sigma
 (h)} \ .
$$

{\bf (XII.3.2) Theorem}. {\it For all regular elliptic element $h$ of
$G$, we have}
$$
\Theta_\pi (h) = \sum_{q=0}^{{\rm dim}\, \vert
X[L]\vert^h} \sum_{\sigma (h)\in \vert X[L]\vert_q^h} (-1)^q \, {\rm
Tr}\, (h, \lambda_\sigma )\ ,
$$
where $\vert X[L]\vert_q^h$ denotes the set of $q$-simplices in $\vert
X[L]\vert^h$.

{\it Proof}. We have to prove that
$$
(f_{\rm EP}^{\VV})^{\vee}(h)=\sum_{q=0}^{d}\sum_{\sigma (h)\in \vert
X[L]\vert^h_q}(-1)^q {\bar \tau}_\sigma (h)\ .
$$

Let $\psi\in \HH (G)_\chi$ be any function with support in
$\KK_\sigma$, for some $q$-dimensional simplex $\sigma$ of $X[L]$,
such that $\psi_{\vert \KK_\sigma}$ is a class function. Let
$(G.\sigma )^h$ be the set of simplices in the $G$-orbit of $\sigma$
that are fixed by $h$. Finally let $G_h$ denote the centralizer of $h$ in $G$.
Following [SS2], page 141, we write

$$
\eqalign{\int_{G/Z}\psi (g^{-1}hg)d\mu_{G/Z} (\dot{g}) & =
\sum_{g\in G_h\backslash G/\KK_\sigma ,\ g^{-1}hg\in \KK_\sigma} \psi (g^{-1}hg)\mu_{G/Z}(G_h g\KK_\sigma /Z)\cr
  & = \sum_{g\sigma\in G_h\backslash (G.\sigma )^h}\psi (g^{-1}hg)\mu_{G/Z}(\KK_\sigma /Z) .
 [G_h : G_h \cap \KK_{g\sigma}]\cr
  & = \mu_{G/Z}(\KK_\sigma /Z).\sum_{g\sigma \in (G.\sigma )^h} \psi (g^{-1}hg)\ .\cr
}
$$
We then apply this to each component of our Euler-Poincar\'e function $f_{\rm EP}^\VV$:
$$
\eqalign{(f_{\rm EP}^\VV )^{\vee}(h) & =
\sum_{q=0}^{d} \sum_{\sigma \in \FF_q}(-1)^q . \sum_{g\sigma\in (G.h)^h}
({\bar \tau}_\sigma^\VV .\epsilon_\sigma )(g^{-1}hg)  \cr
 & = \sum_{q=0}^d \sum_{\sigma \in \FF_q}\sum_{g\sigma \in (G.\sigma )^h}(-1)^q . \epsilon_{g\sigma}(h)
{\bar \tau}_{g\sigma}^\VV (h)  \cr
 & = \sum_{q=0}^d\sum_{\sigma \in (X[L]_q)^h} (-1)^q.\epsilon_\sigma (h).{\bar \tau}_\sigma^\VV
(h) \cr
 & = \sum_{q=0}^{{\rm dim}\, \vert X[L]\vert^h} \sum_{\sigma(h) \in \vert X[L]\vert^h_q}(-1)^{{\rm dim}\, \sigma (h)} .
{\bar \tau}_\sigma^\VV (h)   \cr
}
$$
and we are done.

\bigskip

\centerline{\bf XII.4 The character of discrete series representations at minimal elements.}
\bigskip

In this section we prove that the character formula  of theorem (XII.3.2) takes a striking simple form
under a simple assumption on the regular elliptic element $h$.

Let $\gamma\in G$ satisfying: the algebra $K:=F[\gamma ]\subset A$ is a field (we shall
assume later that the extension $K/F$ is separable, but we do not need this hypothesis
for the moment). Let $v_K$ denote the normalized valuation of $K$. Following [BK1](1.4.14),
one says that $\gamma$ is {\it minimal over $F$} if it satifies:
\medskip

 (i) ${\rm gcd}(v_K (\gamma ), e(K/F))=1$,
\smallskip

 (ii) $\varpi_F^{-v_K(\gamma )}\gamma^{e(K/F)}+\pfr_K$ generates the extension
of residue fields $\Fdss_K /\Fdss$.
\medskip

Here $\varpi_F$ is some uniformizer of $F$ that we fix once for all.

>From [BK1], Exercice (1.5.6), page 44, we have the following result.
\medskip

{\bf (XII.4.1) Lemma}.  {\it Assume that $\gamma \in G$ is minimal over $F$ and let
$\Afr$ be a hereditary order of $A$. Then $\gamma$ normalizes $\Afr$ if, and only if,
$K^{\times}$ normalizes $\Afr$.}
\medskip

Our next  result is a more precise version of this lemma.
\medskip

{\bf (XII.4.2) Proposition}. {\it Assume that $\gamma\in G$ is minimal over $F$.}

(i) {\it We have $X^{\gamma}=X^{K^\times}$ (fixed points set in the geometric realizations). In particular
$X^\gamma$ coincides with the canonical image of $X_K$ in $X$ (cf. Theorem (I.2.1)).}

(ii) {\it In particular, if $K/F$ is a maximal subfield extension of $A$, then $X^\gamma$ reduces
to a single point $x_\gamma$, isobarycenter of simplex corresponding to the unique hereditary order $\Afr_\gamma$ normalized
by $K^\times$ (it has $\ofr_F$-period $e(K/F)$. }
\medskip

{\it Proof}. We use the lattice model of the geometric realization of $X$ given in [BL] {\S}I. Let
us describe this model. Let $L(V)$ denote the set of $\ofr_F$-lattices in $V$. Let ${\rm Latt}^1_{\ofr_F}(V)$ denote the set of functions
$\Lambda~: \ \Rdss \lra L(V)$ satisfying:
\medskip

 -- $\Lambda$ is non-increasing, that is $\Lambda (r)\subset \Lambda (s)$, if $r\geq s$,

 -- $\Lambda$ is periodical, that is $\Lambda (r+1)=\pfr_F \Lambda (r)$, $r\in \Rdss$,

 -- $\Lambda$ is left-continuous for the discrete topology on $L(V)$: for all $r\in\Rdss$, there exists
$\epsilon >0$, such that $\Lambda$ is constant on the segment $[r-\epsilon ,r]$.
\medskip

We let $G$ acts on ${\rm Latt}^1_{\ofr_F} (V)$ by
$$
(g.\Lambda ) (r)=g.\Lambda (r) ,\ g\in G , \ r\in \Rdss\ .
$$

We define the set ${\rm Latt}_{\ofr_F}(V)$ of {\it lattice functions}
 in $V$ as the quotient\break
  ${\rm Latt}_{\ofr_F}^1 (V)/\sim$
for the equivalence relation defined by $\Lambda_1 \sim \Lambda_2$, if
 there exists $s\in \Rdss$  such that $\Lambda_1 (r)
=\Lambda_2 (r+s)$, for all $r\in \Rdss$. Then ${\rm Latt}_{\ofr_F}(V)$ is a $G$-set in an obvious way.

 The point of [BL] {\S}I is that, as a $G$-set,  the geometric realization of $X$ is naturally
isomorphic to ${\rm Latt}_{\ofr_F}(V)$.

 Let $\bar \Lambda$ be a lattice function, with representative $\Lambda \in {\rm Latt}_{\ofr_F}^1 (V)$.
 Assume that $\gamma. {\bar \Lambda}={\bar \Lambda}$. We must prove
 that $\bar \Lambda$ is
fixed by $K^\times$.
 Consider the lattice chain ${\cal L}=\{ \Lambda (r)\ ; \ r\in \Rdss\}$, and let $\Afr ({\cal L} )$
 and $\sigma_{\cal L}$ be the associate hereditary order and simplex respectively.
Then by [BL] Proposition (3.1), $\bar \Lambda$ lies in the interior of the simplex $\sigma_{\cal L}$.
It follows that $\sigma_{\cal L}$ is stabilized by $\gamma$
and therefore that $\Afr ({\cal L})$ is normalized by $\gamma$. Applying Lemma (XII.4.1),
we obtain that $\Afr ({\cal L})$ is normalized by $K^\times$.
In particular it follows that $\cal L$ is a chain of $\ofr_K$-lattices in $V$, and that for all $r\in \Rdss$,
 $\Lambda (r)$ is fixed by $\ofr_K^\times$. Hence
$\bar \Lambda$ is fixed by $\ofr_K^\times$. By condition (i) in the definition
 of a minimal element, there exist  integers $r$, $s$ such that
$\varpi_K :=\varpi_F^r \gamma^s$ is a uniformizer of $K$, and it follows that
 $\bar \Lambda$ is fixed by $\varpi_K$. Hence it is fixed
by $K^\times =\langle \varpi_K \rangle \ofr_K^\times$, as required.
\bigskip

 With the notation as above, we fix an unramified twist of an irreducible  discrete series
 representation $(\pi , \VV )$ of $G$ with type $(J,\lambda )$. Its coefficient system ${\cal C}(\pi)$
has support $X[L]$. We also fix an elliptic regular element $\gamma\in G$ assumed to
be minimal over $F$. In other words, $\gamma$ is minimal
over $F$ and the field extension $K/F$ is separable and maximal.
In particuler the fixed point set $X^\gamma$ is reduced to a single point $x_\gamma$, isobarycenter
of simplex $\sigma_\gamma$ attached to a principal hereditary order $\Afr_\gamma$ with $\ofr_F$-period $e(K/F)$.
\medskip

{\bf (XII.4.3) Lemma}. {\it With the notation as above, we have that $x_\gamma\in X[L]$ if,
and only if, $f(L/F)\vert f(K/F)$ and $e(L/F)\vert e(K/F)$. }
\medskip

{\it Proof}. Using the numerical criterion (I.3.5),
 we have that $x_\gamma\in X[L]$ (that is $\Afr_\gamma$ has a $G$-conjugate normalized by $L^\times$)
if, and only if, with the notation of {\S}I, we have:

i) $f(L/F)\vert d(\Afr_\gamma )_k$, for all $k\in \Zdss$,

ii) $e(L/F)\vert e(\Afr_\gamma /\ofr_F )/p(\Afr_\gamma )$.

 But $\Afr_\gamma$ being principal with period $e(\Afr_\gamma /\ofr_F )=e(K/F)$, we easily see that
$(d(\Afr_\gamma )_k$ is constant with value $f(K/F)$ and
that $p(\Afr_\gamma )=1$. The lemma follows.
\medskip

As a straightforward consequence of the previous lemma and theorem (XII.3.2), we obtain the following simple formula
for the value of the Harish-Chandra character at a minimal element.
\medskip

{\bf (XII.4.4) Proposition}. {\it The Harish-Chandra character of the discrete series representation $(\pi , \VV )$ satisfies :
$$
\Theta_\pi (\gamma )=\cases{
{\rm Tr}\, (\gamma , \lambda_{\sigma_\gamma}) & if $f(L/F)\vert f(K/F)$ and $e(L/F)\vert f(K/F)$,  \cr
0 & otherwise. \cr
}
$$
}

In some particular cases, the same formula was obtained by the first
author in [Br2] using a different approach.

\vfill\eject

{\bf References}

\parindent=20truept

\ref{[BF]} Bushnell C., Fr\"ohlich A.: Nonabelian congruence Gauss
sums and $p$-adic simple algebras. Proc. London Math. Soc. 50, 207-264
(1985)

\ref{[BH]} Bushnell C., Henniart H. : Local tame lifting for ${\rm GL}(N)$,
 I : Simple characters. Publ. math. IHES 83, 105-233 (1996)

\ref{[BK1]} Bushnell C., Kutzko P.: The admissible dual of $GL(n)$
via compact open subgroups. Ann. Math. Studies 129. Princeton Univ.
Press 1993

\ref{[BK2]} Bushnell C., Kutzko P.: Simple types in $GL(N)$:
Computing conjugacy classes. Contemp. Math. 177, 107-135 (1994)

\ref{[BK3]} Bushnell C., Kutzko P.: Smooth representations of
reductive $p$-adic groups: structure theory via types. Proc. London
Math. Soc. 77, 582-634 (1998)

\ref{[BL]} Broussous P., Lemaire B.: Building of ${\rm GL}(m,D)$ and
centralizers.  Transformation Groups 7, 15-50 (2002)

\ref{[BT]} Bruhat F., Tits J.: Sch\'emas en groupes et immeubles des
groupes classiques sur un corps local, ${1}^{\rm re}$ partie~: le
groupe lin\'eaire g\'en\'eral. Bull. Soc. Math. France 112, 259-301
(1984)

\ref{[BTI]} Bruhat F., Tits J.:  Groupes r\'eductifs sur un corps local.
 (French) Inst. Hautes Etudes Sci. Publ. Math. No. 41 (1972), 5–251.

\ref{[Ba]} Badulescu I.: Un r\'esultat de transfert et un r\'esultat
d'int\'egrabilit\'e locale des caract\`eres en caract\'eristique non
nulle, J. reine angew. Math {\bf 565}, 101--124, 2003.

\ref{[Br]} Broussous P.: Acyclycity of Schneider and Stuhler's coefficient systems:
another apporach in the level $0$ case, Journal of Algebra 279 (2004) 737--748.

\ref{[Br2]} Broussous P.:  Transfert du pseudo-coefficient de Kottwitz et formules de caract\`ere
pour la s\'erie discr\`ete de GL(N) sur un corps local,  to appear in Canadian J. of Math. 2013.

\ref{[Ka]} Kazhdan D.: Cuspidal geometry of $p$-adic groups,
J. d'Analyse Math. {\bf 47}, 1--36 1986.

\ref{[Kot]} Kottwitz, R.: Tamagawa numbers, Ann. Math. {\bf 127}
(1988), 629--646.

\ref{[Me]} Meyer, R.: Homological algebra for Schwartz algebras of reductive
$p$-adic groups, in Noncommutative
geometry and number theory, Aspects Math., E37, pp. 263-300. Wiesbaden: Vieweg 2006

\ref{[Rei]} Reiner, I.: Maximal Orders. New York: Academic Press 1975

\ref{[SS]} Schneider P., Stuhler U.: Resolutions for smooth
representations of the general linear group over a local field. J.
reine angew. Math. 436, 19-32 (1993)

\ref{[SS2]} Schneider P., Stuhler U.: Representation theory and
sheaves on the Bruhat-Tits building, Inst. Hautes \'Etudes
Sci. Publ. Math. No. 85 (1997), 97--191.

\ref{[SZ]} Schneider P., Zink E.-W.: $K$-types for the tempered
components of a $p$-adic general linear group. J. reine angew. Math.
517, 161-208 (1999)

\ref{[SZ2]} Schneider P., Zink E.-W.: The algebraic theory of tempered representations 
of $p$-adic groups, Part II: Projective generators. Geometric And Funct. Analysis 17, 2018-2065 (2007)

Paul Broussous {\it Laboratoire de Math\'ematiques,
T\'el\'eport 2 - BP 30179,
Bd Marie et Pierre Curie,
86962 Futuroscope Chasseneuil Cedex, France} 

Email adress: Paul.Broussous{@}math.univ-poitiers.fr

Peter Schneider {\it Universit\"at M\"unster, 
Mathematisches Institut, Einsteinstr. 62, 48291 M\"unster,
Germany}

Email adress : pschnei@uni-muenster.de

\end

%% file: defps.tex

%
\magnification1200
\pretolerance=100
\tolerance=200
\hbadness=1000
\vbadness=1000
\linepenalty=10
\hyphenpenalty=50
\exhyphenpenalty=50
\binoppenalty=700
\relpenalty=500
\clubpenalty=5000
\widowpenalty=5000
\displaywidowpenalty=50
\brokenpenalty=100
\predisplaypenalty=7000
\postdisplaypenalty=0
\interlinepenalty=10
\doublehyphendemerits=10000
\finalhyphendemerits=10000
\adjdemerits=160000
\uchyph=1
\delimiterfactor=901
\hfuzz=0.1pt
\vfuzz=0.1pt
\overfullrule=5pt
\hsize=146 true mm
\vsize=8.9 true in
\maxdepth=4pt
\delimitershortfall=.5pt
\nulldelimiterspace=1.2pt
\scriptspace=.5pt
\normallineskiplimit=.5pt
\mathsurround=0pt
\parindent=20pt
\catcode`\_=11
\catcode`\_=8
\normalbaselineskip=12pt
\normallineskip=1pt plus .5 pt minus .5 pt
\parskip=6pt plus 3pt minus 3pt
\abovedisplayskip = 12pt plus 5pt minus 5pt
\abovedisplayshortskip = 1pt plus 4pt
\belowdisplayskip = 12pt plus 5pt minus 5pt
\belowdisplayshortskip = 7pt plus 5pt
\normalbaselines
\smallskipamount=\parskip
 \medskipamount=2\parskip
 \bigskipamount=3\parskip
\jot=3pt
%
%
\def\ref#1{\par\noindent\hangindent2\parindent
 \hbox to 2\parindent{#1\hfil}\ignorespaces}
%
%
\font\tenss=cmss10
\font\sevenss=cmss8 at 7pt
\font\fivess=cmss8 at 5pt
\newfam\ssfam %
\textfont\ssfam=\tenss
\scriptfont\ssfam=\sevenss
\scriptscriptfont\ssfam=\fivess
%
%
%
%
%
%
%
%
%
\catcode`\_=11
\def\suf_fix{}
\def\scaled_rm_box#1{%
 \relax
 \ifmmode
   \mathchoice
    {\hbox{\tenrm #1}}%
    {\hbox{\tenrm #1}}%
    {\hbox{\sevenrm #1}}%
    {\hbox{\fiverm #1}}%
 \else
  \hbox{\tenrm #1}%
 \fi}
\def\suf_fix_def#1#2{\expandafter\def\csname#1\suf_fix\endcsname{#2}}
\def\I_Buchstabe#1#2#3{%
 \suf_fix_def{#1}{\scaled_rm_box{I\hskip-0.#2#3em #1}}
}
\def\rule_Buchstabe#1#2#3#4{%
 \suf_fix_def{#1}{%
  \scaled_rm_box{%
   \hbox{%
    #1%
    \hskip-0.#2em%
    \lower-0.#3ex\hbox{\vrule height1.#4ex width0.07em }%
   }%
   \hskip0.50em%
  }%
 }%
}
\I_Buchstabe B22
\rule_Buchstabe C51{34}
\I_Buchstabe D22
\I_Buchstabe E22
\I_Buchstabe F22
\rule_Buchstabe G{525}{081}4
\I_Buchstabe H22
\I_Buchstabe I20
\I_Buchstabe K22
\I_Buchstabe L20
\I_Buchstabe M{20em }{I\hskip-0.35}
\I_Buchstabe N{20em }{I\hskip-0.35}
\rule_Buchstabe O{525}{095}{45}
\I_Buchstabe P20
\rule_Buchstabe Q{525}{097}{47}
\I_Buchstabe R21 
\rule_Buchstabe U{45}{02}{54}
\suf_fix_def{Z}{\scaled_rm_box{Z\hskip-0.38em Z}}
\catcode`\"=12
\newcount\math_char_code
\def\suf_fix_math_chars_def#1{%
 \ifcat#1A
  \expandafter\math_char_code\expandafter=\suf_fix_fam
  \multiply\math_char_code by 256
  \advance\math_char_code by `#1
  \expandafter\mathchardef\csname#1\suf_fix\endcsname=\math_char_code
  \let\next=\suf_fix_math_chars_def
 \else
  \let\next=\relax
 \fi
 \next}
%
%
%
%
\def\font_fam_suf_fix#1#2 #3 {%
 \def\suf_fix{#2}
 \def\suf_fix_fam{#1}
 \suf_fix_math_chars_def #3.
}
\font_fam_suf_fix
 0rm
 ABCDEFGHIJKLMNOPQRSTUVWXYZabcdefghijklmnopqrstuvwxyz
\font_fam_suf_fix
 2scr
 ABCDEFGHIJKLMNOPQRSTUVWXYZ
\font_fam_suf_fix
 \slfam sl
 ABCDEFGHIJKLMNOPQRSTUVWXYZabcdefghijklmnopqrstuvwxyz
\font_fam_suf_fix
 \bffam bf
 ABCDEFGHIJKLMNOPQRSTUVWXYZabcdefghijklmnopqrstuvwxyz
\font_fam_suf_fix
 \ttfam tt
 ABCDEFGHIJKLMNOPQRSTUVWXYZabcdefghijklmnopqrstuvwxyz
\font_fam_suf_fix
 \ssfam
 ss
 ABCDEFGHIJKLMNOPQRSTUVWXYZabcdefgijklmnopqrstuwxyz
\catcode`\_=8
\def\Cdss{{\fam\ssfam
    \mkern 4.2 mu \mathchoice%
    {\vrule height 6.5pt depth -.55pt width 1pt}%
    {\vrule height 6.5pt depth -.57pt width 1pt}%
    {\vrule height 4.55pt depth -.28pt width .8pt}%
    {\vrule height 3.25pt depth -.19pt width .6pt}%
    \mkern -6.3mu C}}%
\def\Fdss{{\fam\ssfam I\mkern -2.5mu F}}%
\def\Gdss{{\fam\ssfam
    \mkern 3.8 mu \mathchoice%
    {\vrule height 6.5pt depth -.62pt width 1pt}%
    {\vrule height 6.5pt depth -.65pt width 1pt}%
    {\vrule height 4.55pt depth -.44pt width .8pt}%
    {\vrule height 3.25pt depth -.30pt width .6pt}%
    \mkern -5.9mu G}}%
\def\Ldss{{\fam\ssfam I\mkern -2.5mu L}}%
\def\Pdss{{\fam\ssfam I\mkern -2.5mu P}}%
\def\Rdss{{\fam\ssfam I\mkern -2.5mu R}}%
\def\Udss{{\fam\ssfam U\mkern-10mu U}}%
\def\Zdss{{\fam\ssfam Z\mkern-8.1mu Z}}%
%
%
%
%
\font\teneuf=eufm10 
\font\seveneuf=eufm7
\font\fiveeuf=eufm5
\newfam\euffam \def\euf{\fam\euffam\teneuf} 
\textfont\euffam=\teneuf \scriptfont\euffam=\seveneuf
\scriptscriptfont\euffam=\fiveeuf
\def\Afr{{\euf A}}       
\def\Bfr{{\euf B}}       \def\bfr{{\euf b}}
\def\Cfr{{\euf C}}

       \def\ofr{{\euf o}}
\def\Pfr{{\euf P}}       \def\pfr{{\euf p}}